\documentclass[12pt]{book}
\usepackage{amsmath,amssymb,bbm,amsthm,graphicx,url}

\newtheorem{theorem}{Theorem}[chapter]
\newtheorem{corollary}[theorem]{Corollary}
\newtheorem{lemma}[theorem]{Lemma}
\newtheorem{proposition}[theorem]{Proposition}
\newtheorem{remark}[theorem]{Remark}

\theoremstyle{definition}
\newtheorem{definition}[theorem]{Definition}
\newtheorem{example}[theorem]{Example}
\newtheorem{problem}{Problem}

\newcommand{\C}{{\mathbb C}}
\newcommand{\Cd}{{\C^d}}
\newcommand{\D}{{\mathbb D}}
\newcommand{\N}{{\mathbb N}}
\newcommand{\R}{{\mathbb R}}
\newcommand{\Rd}{{\R^d}}
\newcommand{\T}{{\mathbb T}}
\newcommand{\Td}{{\T^d}}
\newcommand{\Z}{{\mathbb Z}}
\newcommand{\Zd}{{\Z^d}}
\newcommand{\pr}{\frac{1}{2\pi}}

\newcommand{\Aside}{\noindent \emph{Aside.}\ }
\newcommand{\Fejer}{Fej\'{e}r\ }
\newcommand{\Proof}{\emph{Proof.}\ }

\newcommand{\charfn}{\mathbbm{1}}
\newcommand{\e}{{\varepsilon}}
\renewcommand{\Re}{\operatorname{Re}}

\newcommand{\Pe}{\operatorname{Pe}}
\newcommand{\side}{\operatorname{side}}
\newcommand{\sign}{\operatorname{sign}}
\newcommand{\sinc}{\operatorname{sinc}}

\newcommand{\fhat}{\widehat{f}}
\newcommand{\fhatj}{\fhat(j)}
\newcommand{\fhatn}{\fhat(n)}
\newcommand{\fhatxi}{\fhat(\xi)}
\newcommand{\ghat}{\widehat{g}}

\newcommand{\ghatn}{\ghat(n)}
\newcommand{\ghatxi}{\ghat(\xi)}

\newcommand{\gcheck}{\check{g}}

\newcommand{\pp}{\frac{1}{2\pi}}
\newcommand{\intp}{\frac{1}{2\pi} \! \int}

\newcommand{\la}{\langle}
\newcommand{\ra}{\rangle}

\newcommand{\lv}{\lVert}
\newcommand{\rv}{\rVert}

\newcommand{\li}{{L^\infty}}

\newcommand{\lird}{{L^\infty(\Rd)}}
\newcommand{\lit}{{L^\infty(\T)}}

\newcommand{\lp}{{L^p}}
\newcommand{\lpr}{{L^p(\R)}}
\newcommand{\lprd}{{L^p(\Rd)}}
\newcommand{\lpt}{{L^p(\T)}}
\newcommand{\lptd}{{L^p(\Td)}}
\newcommand{\lo}{{L^1}}
\renewcommand{\lor}{{L^1(\R)}}
\newcommand{\lord}{{L^1(\Rd)}}
\newcommand{\lot}{{L^1(\T)}}
\newcommand{\lotd}{{L^1(\Td)}}
\newcommand{\lt}{{L^2}}
\newcommand{\ltr}{{L^2(\R)}}
\newcommand{\ltrd}{{L^2(\Rd)}}
\newcommand{\ltt}{{L^2(\T)}}

\title{\textbf{Harmonic Analysis} \\ \textbf{Lecture Notes} \\ \ \\ University of Illinois \\ at Urbana--Champaign}

\author{Richard S. Laugesen
\footnote{Copyright \copyright\ 2017, Richard S. Laugesen
(Laugesen@illinois.edu). This work is licensed under the Creative
Commons Attribution-NonCommercial-NoDerivs 4.0 International
License. To view a copy of this license, visit
\protect\url{https://creativecommons.org/licenses/by-nc-nd/4.0/}.} }

\begin{document}

\maketitle

\section*{Preface}

A \emph{textbook} presents more than any professor can
cover in class. In contrast, these \emph{lecture notes} present
exactly\footnote{modulo some improvements after the fact} what I
covered in Harmonic Analysis (Math 545) at the University
of Illinois, Urbana--Champaign, in Fall 2008.

The first part of the course emphasizes Fourier series, since so
many aspects of harmonic analysis arise already in that classical
context. The Hilbert transform is treated on the circle, for
example, where it is used to prove $L^p$ convergence of Fourier
series. Maximal functions and Calder\'{o}n--Zygmund decompositions
are treated in $\Rd$, so that they can be applied again in the
second part of the course, where the Fourier transform is studied.

Real methods are used throughout. In particular, complex methods
such as Poisson integrals and conjugate functions are \emph{not}
used to prove boundedness of the Hilbert transform.

Distribution functions and interpolation are covered in the
Appendices. I inserted these topics at the appropriate places in my
lectures (after Chapters~\ref{ch:l1} and \ref{ch:hp}, respectively).

The references at the beginning of each chapter provide guidance to
students who wish to delve more deeply, or roam more widely, in the
subject. Those references do not necessarily contain all the
material in the chapter.

Finally, a word on personal taste: while I appreciate a good
counterexample, I prefer spending class time on positive results.
Thus I do not supply proofs of some prominent counterexamples (such
as Kolmogorov's integrable function whose Fourier series diverges at
every point).

I am grateful to Noel DeJarnette, Eunmi Kim, Aleksandra Kwiatkowska,
Kostya Slutsky, Khang Tran and Ping Xu for {\TeX}ing parts of the
document, and to Alexander Tumanov for pointing out a number of typos. 

\vspace{6pt} Please email me with corrections, and with suggested
improvements of any kind.

\vspace{24pt} \noindent Richard S. Laugesen \qquad \qquad \qquad
\textsc{Email:} Laugesen\@@illinois.edu
\\ Department of Mathematics \\ University of Illinois \\ Urbana, IL 61801 \\ U.S.A.

\newpage
\section*{Introduction}

Harmonic analysis began with Fourier's effort to analyze (extract
information from) and synthesize (reconstruct) the solutions of the
heat and wave equations, in terms of harmonics. Specifically, the
computation of Fourier coefficients is analysis, while writing down
the Fourier series is synthesis, and the harmonics in one dimension are $\sin(nt)$
and $\cos(nt)$. Immediately one asks: does the
Fourier series converge? to the original function? In what sense does it converge:
pointwise? mean-square? $L^p$? Do analogous results
hold on $\Rd$ for the Fourier transform?

We will answer these classical qualitative questions (and more!)
using modern quantitative estimates, involving tools such as
summability methods (convolution), maximal operators, singular
integrals and interpolation. These topics, which we address for both
Fourier series and transforms, constitute the theoretical core of
the course. We further cover the sampling theorem, Poisson summation
formula and uncertainty principles.

\vspace{6pt} This graduate course is \emph{theoretical} in nature.
Students who are intrigued by the fascinating \emph{applications} of
Fourier series and transforms are advised to browse \cite{DM},
\cite{Ko} and \cite{SS}, which are all wonderfully engaging books.

\vspace{6pt} If more time (or a second semester) were available, I
might cover additional topics such as: Littlewood--Paley theory for
Fourier series and integrals, Fourier analysis on locally compact
abelian groups \cite{R} (especially Bochner's theorem on Fourier
transforms of nonnegative functions), short-time Fourier transforms
\cite{Gr}, discrete Fourier transforms, the Schwartz class and
tempered distributions and applications in Fourier analysis
\cite{St}, Fourier integral operators (including solutions of the
wave and Schr\"{o}dinger equations), Radon transforms, and some
topics related to signal processing, such as maximum entropy,
spectral estimation and prediction \cite{B}. I might also cover
multiplier theorems, ergodic theorems, and almost periodic
functions.

\tableofcontents

\part{Fourier series}

\chapter{Fourier coefficients: basic properties} \label{ch:fc}

\subsubsection*{Goal}

Derive basic properties of Fourier coefficients

\subsubsection*{Reference} \cite{K} Section I.1

\subsubsection*{Notation}

\noindent $\T = \R / 2\pi \Z$ is the one dimensional \textbf{torus}

\noindent $\lpt = \{ \text{complex-valued, $p$-th power integrable,
$2\pi$-periodic functions} \}$

\noindent $\lv f \rv_\lpt = \big( \intp_\T |f(t)|^p \, dt
\big)^{1/p}$ where $\int_\T$ can be taken over any interval of
length $2\pi$

\noindent Nesting of $L^p$-spaces: $\lit \subset \ltt \subset \lot$

\noindent $C(\T) = \{ \text{complex-valued, continuous,
$2\pi$-periodic functions} \}$, Banach space with norm $\lv \cdot
\rv_\lit$

\noindent Trigonometric polynomial $P(t)=\sum_{n=-N}^N a_n e^{int}$

\noindent Translation $f_\tau(t) = f(t-\tau)$

\begin{definition}
For $f \in \lot$ and $n \in \Z$, define
\begin{align}
\fhatn & = \text{$n$-th \emph{Fourier coefficient} of $f$} \notag \\
& = \intp_\T f(t) e^{-int} \, dt . \label{eq:fhatdef}
\end{align}
The formal series $S[f] = \sum \fhatn e^{int}$ is the \emph{Fourier
series} of $f$.
\end{definition}

\Aside For $f \in \ltt$, note $\fhatn = \la f , e^{int} \ra$ where
$\la f , g \ra = \intp_\T f(t) \overline{g(t)} \, dt$ is that $\lt$
inner product. Thus $\fhatn=$amplitude of $f$ in direction of
$e^{int}$. See Chapter~\ref{ch:l2}.

\begin{theorem}[Basic properties] \label{th:bp}
Let $f,g \in \lot, j,n \in \Z, c \in \C, \tau \in \T$.

\noindent Linearity $\widehat{(f+g)}(n)=\fhatn+\widehat{g}(n)$ and
$\widehat{(cf)}(n)=c\fhatn$

\noindent Conjugation
$\widehat{\overline{f}}(n)=\overline{\fhat(-n)}$

\noindent Trigonometric polynomial $P(t)=\sum_{n=-N}^N a_n e^{int}$
has $\widehat{P}(n) = a_n$ for $|n| \leq N$ and $\widehat{P}(n) = 0$
for $|n|>N$

\noindent $\widehat{\ }$ \ takes translation to modulation,
$\widehat{\, f_\tau \,}(n) = e^{-in\tau} \fhatn$

\noindent $\widehat{\ }$ \ takes modulation to translation,
$[f(t)e^{ijt}]\widehat{\ }(n) = \fhat(n-j)$

\noindent $\widehat{\ } : \lot \to \ell^\infty(\Z)$ is bounded, with
$|\fhatn| \leq \lv f \rv_\lot$

\noindent Hence if $f_m \to f$ in $\lot$ then $\widehat{f_m}(n) \to
\fhatn$ (uniformly in $n$) as $m \to \infty$.
\end{theorem}
\begin{proof}
Exercise.
\end{proof}

\begin{lemma}[Difference formula] \label{le:diffT} For $n \neq 0$,
\[
\fhatn = \frac{1}{4\pi} \int_\T [f(t)-f(t-\pi/n)] \, e^{-int} \, dt
.
\]
\end{lemma}
\begin{proof}
\begin{align}
\fhatn & = - \intp_\T f(t) e^{-in(t+\pi/n)} \, dt && \text{since
$e^{-i\pi}=-1$} \notag \\
& = - \intp_\T f(t-\pi/n) e^{-int} \, dt \label{eq:diffT1}
\end{align}
by $t \mapsto t-\pi/n$ and periodicity. By \eqref{eq:diffT1} and the
definition \eqref{eq:fhatdef},
\[
\fhatn = \frac{1}{2} \fhatn + \frac{1}{2} \fhatn = \frac{1}{4\pi}
\int_\T [f(t)-f(t-\pi/n)] \, e^{-int} \, dt .
\]
\end{proof}
\begin{lemma}[Continuity of translation] \label{le:ct} Fix $f \in \lpt, 1 \leq p < \infty$. The map
\begin{align*}
\phi : \T & \to \lpt \\
\tau & \mapsto f_\tau
\end{align*}
is continuous.
\end{lemma}
\begin{proof}
Let $\tau_0 \in \T$. Take $g \in C(\T)$ and observe
\begin{align*}
\lv f_\tau - f_{\tau_0} \rv_\lpt & \leq \lv f_\tau - g_\tau \rv_\lpt
+ \lv g_\tau - g_{\tau_0} \rv_\lpt + \lv g_{\tau_0} - f_{\tau_0} \rv_\lpt \\
& = 2 \lv f - g \rv_\lpt + \lv g_\tau - g_{\tau_0} \rv_\lpt \\
& \to 2 \lv f - g \rv_\lpt
\end{align*}
as $\tau \to \tau_0$, by uniform continuity of $g$. By density of
continuous functions in $\lpt, 1 \leq p < \infty$, the difference $f
- g$ can be made arbitrarily small. Hence $\limsup_{\tau \to \tau_0}
\lv f_\tau - f_{\tau_0} \rv_\lpt = 0$, as desired.
\end{proof}
\begin{corollary}[Riemann--Lebesgue lemma] \label{co:rl} $\fhatn \to 0$ as $|n| \to
\infty$.
\end{corollary}
\begin{proof}
Lemma~\ref{le:diffT} implies
\[
|\fhatn| \leq \frac{1}{2} \lv f - f_{\pi/n} \rv_\lot ,
\]
which tends to zero as $|n| \to \infty$ by the $\lo$-continuity of
translation in Lemma~\ref{le:ct}, since $f=f_0$.
\end{proof}

\subsubsection*{Smoothness and decay}

The Riemann--Lebesgue lemma says $\fhatn=o(1)$, with $\fhatn =O(1)$
explicitly by Theorem~\ref{th:bp}. We show the smoother $f$ is, the
faster its Fourier coefficients decay.

\begin{theorem}[Less than one derivative] \label{th:rd1}
If $f \in C^\alpha(\T), 0<\alpha \leq 1$, then
$\fhatn=O(1/|n|^\alpha)$.
\end{theorem}

Here $C^\alpha(\T)$ denotes the H\"older continuous functions: $f
\in C^\alpha(\T)$ if $f \in C(\T)$ and there exists $A>0$ such that
$|f(t)-f(\tau)|\leq A|t-\tau|^\alpha$ whenever $|t-\tau| \leq 2\pi$.

\begin{proof}
\[
\fhatn=\frac{1}{4\pi} \int_\T [f(t)-f(t-\pi/n)] e^{-int} \, dt
\]
by the Difference Formula in Lemma~\ref{le:diffT}. Therefore
\[
|\fhatn| \leq \frac{1}{4\pi} A \left| \frac{\pi}{n} \right|^\alpha
\, 2\pi = \frac{\text{const.}}{|n|^\alpha}.
\]
\end{proof}

\begin{theorem}[One derivative] \label{th:rd2}
If $f$ is $2\pi$-periodic and absolutely continuous ($f \in
W^{1,1}(\T)$) then $\fhatn=o(1/n)$ and $|\fhatn| \leq \lv f^\prime
\rv_\lot/|n|$.
\end{theorem}

\begin{proof}
Absolute continuity of $f$ says
\[
f(t)=f(0)+ \int_0^t f^\prime(\tau) \, d\tau,
\]
where $f^\prime \in \lot$. Integrating by parts gives
\[
\fhatn= \intp_\T f(t)e^{-int} \, dt= \intp_\T \frac{e^{-int}}{in}
f^\prime(t) \, dt.
\]
By Riemann-Lebesgue applied to $f^\prime$,
\[
\fhatn=\frac{1}{in}\widehat{f^\prime}(n)=\frac{1}{in}o(1)=o(\frac{1}{n})
,
\]
with
\[
|\fhatn| = \frac{1}{|n|}|\widehat{f^\prime}(n)| \leq
\frac{1}{|n|} \lv f^\prime \rv_\lot.
\]
\end{proof}

\begin{theorem}[Higher derivatives] \label{th:rd3}
If $f$ is $2\pi$-periodic and $k$ times differentiable ($f \in
W^{k,1}(\T)$) then $\fhatn=o(1/|n|^k)$ and $|\fhatn| \leq \lv
f^{(k)} \rv_\lot/|n|^k$.
\end{theorem}

\begin{proof}
Integrate by parts $k$ times.
\end{proof}

\begin{remark}\rm
Similar decay results hold for functions of bounded variation,
provided one integrates by parts using the Lebesgue--Stieltjes
measure $df(t)$ instead of $f^\prime(t) \, dt$.
\end{remark}

\subsubsection*{Convolution}

\begin{definition} Given $f,g \in \lot$, define their
\emph{convolution}
\[
(f*g)(t) = \intp_\T f(t-\tau) g(\tau) \, d\tau , \qquad t \in \T.
\]
\end{definition}
\begin{theorem}[Convolution and Fourier coefficients] \label{th:cfc}
If $f \in \lpt, 1 \leq p \leq \infty$, and $g \in \lot$, then $f * g \in \lpt$ with
\begin{equation} \label{eq:convL1}
\lv f*g \rv_\lpt \leq \lv f \rv_\lpt \lv g \rv_\lot
\end{equation}
and
\[
\widehat{(f*g)}(n) = \fhatn \ghatn , \qquad n \in \Z .
\]
Further, if $f \in C(\T)$ and $g \in \lot$ then $f * g \in C(\T)$.
\end{theorem}
Thus $\widehat{\ }$ takes convolution to multiplication.
\begin{proof}
That $f * g \in \lpt$ satisfies \eqref{eq:convL1} is exactly Young's Theorem~\ref{th:yt}. Then by Fubini's theorem,
\begin{align*}
\widehat{(f*g)}(n) & = \intp_\T \big( \intp_\T f(t-\tau) g(\tau) \,
d\tau \big) e^{-int} \, dt \\
& = \intp_\T \big( \intp_\T f(t-\tau) e^{-in(t-\tau)} \, dt
\big) g(\tau) e^{-in\tau} \, d\tau \\
& = \fhatn \ghatn .
\end{align*}
Finally, if $f \in C(\T)$ and $g \in \lot$ then $f*g$ is continuous because $(f * g)(t+\delta) \to (f*g)(t)$ as $\delta \to 0$ by uniform continuity of $f$.
\end{proof}
\paragraph*{Convolution facts} \cite[Section~I.1.8]{K}

\vspace{6pt} \noindent 1. Convolution is commutative:
\begin{align*}
(f*g)(t) & = \intp_\T f(t-\tau) g(\tau) \, d\tau \\
& = \intp_\T f(\theta) g(t-\theta) \, d\theta && \text{where
$\tau = t - \theta, d\tau = -d\theta$} \\
& = (g*f)(t) .
\end{align*}
Convolution is also associative, and linear with respect to $f$ and
$g$.

\vspace{6pt} \noindent 2. Convolution is continuous on $\lpt$: if
$f_m \to f \in \lpt, 1 \leq p \leq \infty$, and $g \in \lot$ then $f_m*g \to
f*g$ in $\lpt$.

\noindent \Proof Use linearity and \eqref{eq:convL1}, to prove $f_m*g \to
f*g$ in $\lpt$.

\vspace{6pt} \noindent 3. Convolution with a trigonometric
polynomial gives a trigonometric polynomial: if $f \in \lot$ and
$P(t) = \sum_{j=-n}^n a_j e^{ijt}$ then
\begin{equation} \label{eq:convtrig}
(P*f)(t) = \sum_{j=-n}^n a_j \fhatj e^{ijt} .
\end{equation}
\noindent \Proof
\begin{align*}
(P*f)(t) & = \sum_{j=-n}^n a_j \intp_\T e^{ij(t-\tau)} f(\tau) \,
d\tau \\
& = \sum_{j=-n}^n a_j e^{ijt} \fhatj .
\end{align*}
[Sanity check: $\widehat{(P*f)}(j) = a_j \fhatj = \widehat{P}(j)
\fhatj$ as expected.]

\vspace{6pt} \noindent More generally, \eqref{eq:convtrig} holds for
$P(t)=\sum_{j=-\infty}^\infty a_j e^{ijt}$ provided $\{ a_j \} \in
\ell^1(\Z)$.

\chapter{Fourier series: summability in norm} \label{ch:fs}

\subsubsection*{Goal}

Prove summability (averaged convergence) in norm of Fourier series

\subsubsection*{Reference} \cite{K} Section I.2

\vspace{18pt} Write
\begin{align*}
(S_n f)(t) & = \sum_{j=-n}^n \fhatj e^{ijt} \\
& = \text{$n$-th partial sum of Fourier series of $f$.}
\end{align*}
In Chapter~\ref{ch:cn} we prove norm convergence of Fourier series:
$S_n(f) \to f$ in $\lpt$, when $1<p<\infty$. In this chapter we
prove \textbf{summability} of Fourier series, meaning $\sigma_n(f)
\to f$ in $\lpt$ when $1 \leq p < \infty$, where
\begin{align*}
\sigma_n(f) & = \frac{S_0(f) + \cdots + S_n(f)}{n+1} = \sum_{j=-n}^n
\left( 1 - \frac{|j|}{n+1} \right) \fhatj e^{ijt} \\
& = \text{arithmetic mean of partial sums.}
\end{align*}

\Aside Norm convergence is stronger than summability. Indeed, if a
sequence $\{ s_n \}$ in a Banach space converges to $s$, then the
arithmetic means $(s_0 + \cdots + s_n)/(n+1)$ also converge to $s$
(Exercise).

\begin{definition}
A \emph{summability kernel} is a sequence $\{ k_n \}$ in $\lot$
satisfying:
\begin{align}
\intp_\T k_n(t) \,  dt & = 1 && \textsc{(Normalization)} \tag{S1} \label{eq:S1} \\
\sup_n \intp_\T |k_n(t)| \, dt & < \infty && \textsc{($L^1$ bound)} \tag{S2} \label{eq:S2} \\
\lim_{n \to \infty} \int_{\{ \delta < |t| < \pi \}} |k_n(t)| \,
dt & = 0 && \textsc{($L^1$ concentration)} \tag{S3} \label{eq:S3} \\
& && \text{for each $\delta \in (0,\pi)$.} \notag
\end{align}
Some kernels satisfy a stronger concentration property:
\begin{align}
\lim_{n \to \infty} \sup_{\delta < |t| < \pi} |k_n(t)| & = 0 && \textsc{($L^\infty$ concentration)} \tag{S4} \label{eq:S4} \\
& && \text{for each $\delta \in (0,\pi)$.}  \notag
\end{align}
Call the kernel \emph{positive} if $k_n \geq 0$ for each $n$.
\end{definition}

\begin{example}
Define the \textbf{Dirichlet kernel}
\begin{align}
D_n(t) & = \sum_{j=-n}^n e^{ijt} \label{eq:D1} \\
& = \frac{e^{i(n+1)t}-e^{-int}}{e^{it}-1} && \text{by geometric
series} \label{eq:D2} \\
& = \frac{\sin \big( (n+\frac{1}{2})t \big)}{\sin \big( \frac{1}{2}t
\big)}
\end{align}
\begin{figure}
\begin{center}
  \includegraphics[scale=0.8]{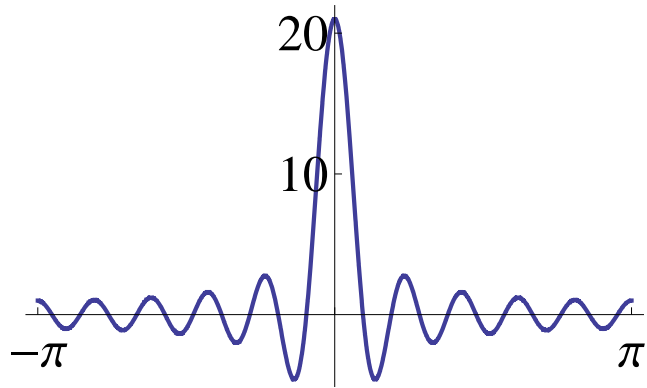}
  \caption{\label{Dnfig}
    Dirichlet kernel with $n=10$}
\end{center}
\end{figure}
\eqref{eq:S1} holds by \eqref{eq:D1}. You can show (optional
exercise) that $\lv D_n \rv_\lot \sim (\text{const.}) \log n$ as $n
\to \infty$, so that \eqref{eq:S2} fails.

$\therefore \ \{ D_n \}$ is \textbf{not} a summability kernel.
\end{example}

\begin{figure}
\begin{center}
  \includegraphics[scale=0.8]{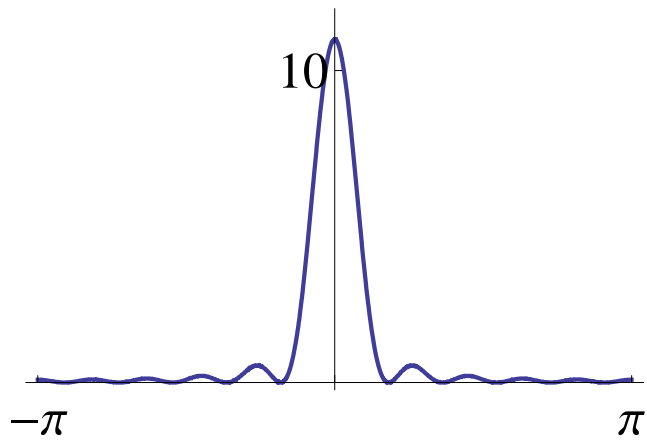}
  \caption{\label{Fnfig}
    \Fejer kernel with $n=10$}
\end{center}
\end{figure}

\begin{example} Define the \textbf{\Fejer kernel}
\begin{align}
F_n(t) & = \frac{D_0(t) + \cdots + D_n(t)}{n+1} \label{eq:F1} \\
& = \sum_{j=-n}^n \left( 1 - \frac{|j|}{n+1} \right) e^{ijt} &&
\text{by \eqref{eq:F1} and \eqref{eq:D1}} \label{eq:F2} \\
& = \frac{1}{n+1} \left( \frac{\sin \big( \frac{n+1}{2}t \big)}{\sin
\big( \frac{1}{2}t \big)} \right)^{\! \! \! 2} \label{eq:F3} &&
\text{by \eqref{eq:F1}, \eqref{eq:D2} and geometric series}
\end{align}
\eqref{eq:S1} holds by \eqref{eq:F2}, and $F_n \geq 0$ so that
\eqref{eq:S2} holds also. For \eqref{eq:S4},
\begin{align*}
\sup_{\delta < |t| < \pi} |F_n(t)| & \leq \frac{1}{n+1}
\frac{1}{\sin^2 \big( \frac{1}{2}\delta
\big)} && \text{by \eqref{eq:F3}} \\
& \to 0 && \text{as $n \to \infty$.}
\end{align*}
$\therefore \ \{ F_n \}$ is a positive summability kernel.
\end{example}

\begin{example} Define the \textbf{Poisson kernel}
\begin{align}
P_r(t) & = 1 + 2 \sum_{j=1}^\infty r^j \cos(jt) \label{eq:P1} \\
& = \sum_{j=-\infty}^\infty r^{|j|} e^{ijt} \label{eq:P2} \\
& = \frac{1-r^2}{1-2r \cos t + r^2} \label{eq:P3}
\end{align}
by summing two geometric series ($j<0$ and $j \geq 0$) in
\eqref{eq:P2} and simplifying.

The Poisson kernel is indexed by $r \in (0,1)$, with limiting
process $r \nearrow 1$. After suitably modifying the definition of
summability kernel, we see \eqref{eq:S1} holds by \eqref{eq:P1}, and
$P_r \geq 0$ by \eqref{eq:P3} so that \eqref{eq:S2} holds also. For
\eqref{eq:S4},
\begin{align*}
\sup_{\delta < |t| < \pi} |P_r(t)| & \leq \frac{1-r^2}{1-2r \cos \delta + r^2} && \text{by \eqref{eq:P3}} \\
& \to 0 && \text{as $r \nearrow 1$.}
\end{align*}
$\therefore \ \{ P_r \}$ is a positive summability kernel.
\end{example}

\begin{example} Define the \textbf{Gauss kernel}
\begin{align}
G_s(t) & = \sum_{j=-\infty}^\infty e^{-j^2 s} e^{ijt} \label{eq:G1} \\
& = \frac{2\pi}{\sqrt{4\pi s}} \sum_{n=-\infty}^\infty e^{-(t+2\pi
n)^2/4s} \label{eq:G2}
\end{align}
by Example~\ref{ex:pg} later in the course.
\begin{figure}
\begin{center}
  \includegraphics[scale=0.8]{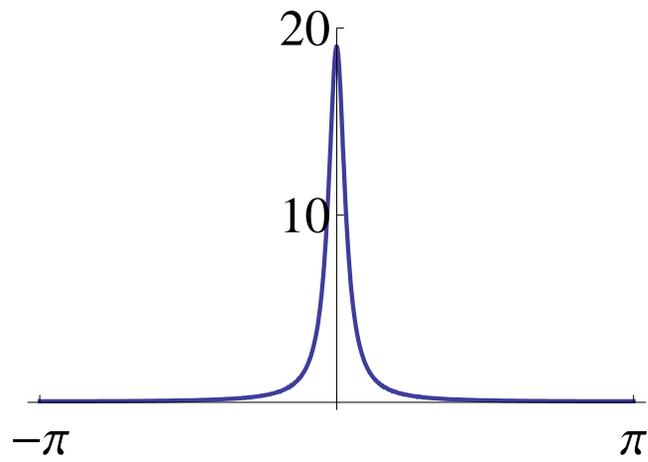}
  \caption{\label{Prfig}
    Poisson kernel with $r=0.9$}
\end{center}
\end{figure}

\begin{figure}
\begin{center}
  \includegraphics[scale=0.8]{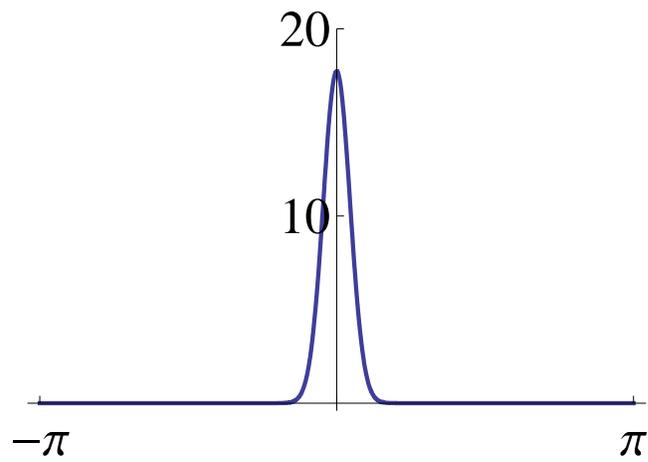}
  \caption{\label{Gsfig}
    Gauss kernel with $s=0.01$}
\end{center}
\end{figure}

The Gauss kernel is indexed by $s \in (0,\infty)$, with limiting
process $s \searrow 0$. The analogue of \eqref{eq:S1} holds by
\eqref{eq:G1}, and $G_s \geq 0$ by \eqref{eq:G2} so that
\eqref{eq:S2} holds also. For \eqref{eq:S4},
\begin{align*}
\sup_{\delta < |t| < \pi} |G_s(t)| & \leq \frac{2\pi}{\sqrt{4\pi s}} \left[ e^{-\delta^2/4s} + \sum_{n \neq 0} e^{-(\pi n)^2/4s} \right] && \text{by \eqref{eq:G2}} \\
& \to 0 && \text{as $s \searrow 0$.}
\end{align*}
$\therefore \ \{ G_s \}$ is a positive summability kernel.
\end{example}

\subsubsection*{Connection to Fourier series}

\[
\fbox{$S_n(f) = D_n * f$}
\]
\noindent \Proof $D_n(t) \overset{\eqref{eq:D1}}{=} \sum_{j=-n}^n 1
e^{ijt}$ implies
\[
(D_n * f)(t)=\sum_{j=-n}^n 1 \fhatj e^{ijt} = S_n(f)
\]
by Convolution Fact \eqref{eq:convtrig}.

\[
\fbox{$\sigma_n(f) = F_n * f$}
\]
\noindent \Proof $F_n(t) \overset{\eqref{eq:F1}}{=} \sum_{j=-n}^n
(1-\frac{|j|}{n+1}) e^{ijt}$ implies
\[
(F_n * f)(t)=\sum_{j=-n}^n \big( 1-\frac{|j|}{n+1} \big) \fhatj
e^{ijt} = \sigma_n(f)
\]
by Convolution Fact \eqref{eq:convtrig}. Alternatively, use that
$\sigma_n(f)=[S_0(f)+\cdots+S_n(f)]/(n+1)$ and
$F_n=[D_0+\cdots+D_n]/(n+1)$.
\[
\text{Thus for summability of Fourier series, we want $F_n * f \to
f$.}
\]

\[
\fbox{$\text{Abel mean of $S[f]$} = P_r * f$}
\]
\noindent \Proof $P_r(t) \overset{\eqref{eq:P2}}{=}
\sum_{j=-\infty}^\infty r^{|j|} e^{ijt}$ implies
\begin{equation} \label{eq:pabel}
(P_r * f)(t)=\sum_{j=-\infty}^\infty r^{|j|} \fhatj e^{ijt}
\end{equation}
by Convolution Fact \eqref{eq:convtrig} (with the series converging absolutely and uniformly), and this last expression is
the Abel mean of $S[f]$.

\subsubsection*{Summability in norm}

\begin{theorem}[Summability in $\lpt$ and $C(\T)$] \label{th:sk}
If $\{ k_n \}$ is a summability kernel and $f \in \lpt, 1 \leq p < \infty$, then
\[
k_n * f \to f \qquad \text{in $\lpt$, \quad as $n \to \infty$.}
\]
Similarly, if $f \in C(\T)$ then $k_n * f \to f$ in $C(\T)$.
\end{theorem}
\begin{proof}
Let $\e > 0$. By \eqref{eq:S2} and continuity of translation on $\lpt$ (Lemma~\ref{le:ct}), we can
choose $0<\delta<\pi$ such that
\begin{equation} \label{eq:sd1}
\max_{|\tau| \leq \delta} \lv f_\tau - f \rv_\lpt \cdot
\sup_n \lv k_n \rv_\lot <  \e .
\end{equation}
Then
\begin{align*}
& \big\lv (k_n * f)(t) - f(t) \big\rv_\lpt
\\
& = \big\lv \intp_\T k_n(\tau) [f_\tau(t) - f(t)] \, d\tau
\big\rv_\lpt && \text{by \eqref{eq:S1}} \\
& \leq \intp_\T |k_n(\tau)| \lv f_\tau - f \rv_\lpt \, d\tau \\
& \qquad \qquad \qquad \text{by Minkowski's Integral Inequality, Theorem~\ref{th:mii},} \\
& = \frac{1}{2\pi} \left( \int_{-\delta}^\delta + \int_{\{ \delta < |\tau| < \pi\}} \right) |k_n(\tau)| \lv f_\tau - f \rv_\lpt \, d\tau \\
& \leq \max_{|\tau| \leq \delta} \lv f_\tau - f \rv_\lpt \frac{1}{2\pi} \int_{-\delta}^\delta |k_n(\tau)| \, d\tau \\
& \quad + \max_{|\tau| \leq \pi} \lv f_\tau - f \rv_\lpt \frac{1}{2\pi} \int_{\{ \delta < |\tau| < \pi\}} |k_n(\tau)| \, d\tau \\
& < \e + \e
\end{align*}
by \eqref{eq:sd1} and \eqref{eq:S3}, for all large $n$.

If $f \in C(\T)$ then repeat the argument with $p=\infty$, using
uniform continuity of $f$ to get that $f_\tau \to f$ in $\lit$.
\end{proof}

\subsubsection*{Consequences}

\noindent $\bullet$ Summability of Fourier series in $C(\T),
\lpt, 1 \leq p <  \infty$:
\[
\sigma_n(f) \to f
\]
in norm.

\noindent \Proof Choose $k_n=F_n=$\,\Fejer kernel. Then
$\sigma_n(f)=F_n * f \to f$ in norm by Theorem~\ref{th:sk}.

\vspace{6pt} \noindent $\bullet$ Trigonometric polynomials are dense
in $C(\T), \lpt, 1 \leq p <  \infty$.

\noindent \Proof $\sigma_n(f)$ is a trigonometric polynomial
arbitrarily close to $f$.

\vspace{6pt} \Aside Density of trigonometric polynomials in $C(\T)$
proves the Weierstrass Trigonometric Approximation Theorem.

\vspace{6pt} \noindent $\bullet$ Uniqueness theorem:
\begin{equation} \label{eq:ut}
\text{if $f,g \in \lot$ with $\fhatn=\ghatn$ for all $n$, then $f=g$
in $\lot$.}
\end{equation}
In other words, the map $\widehat{\ } : \lot \to \ell^\infty(\Z)$ is
injective.

\noindent \Proof $F_n * f = F_n * g$ by Convolution Fact
\eqref{eq:convtrig}, since $\fhat=\ghat$. Letting $n \to \infty$
gives $f=g$.

\subsubsection*{Connection to PDEs}

To finish the section, we connect our summability kernels to some
important partial differential equations. Fix $f \in \lot$.

\vspace{6pt} \noindent 1. The Poisson kernel solves Laplace's
equation in a disk:
\[
v(re^{it}) = (P_r * f)(t)  = \frac{1}{2\pi} \int_\T
\frac{1-r^2}{1-2r \cos(t-\tau)+r^2} f(\tau) \, d\tau
\]
solves
\[
\Delta v = v_{rr}+r^{-1} v_r+r^{-2}v_{tt}=0
\]
on the unit disk $\{ r<1 \}$, with boundary value $v(1,t)=f(t)$ in the sense of
Theorem~\ref{th:sk}.

That is, $v$ is the harmonic extension of $f$ from the boundary
circle to the disk.

\vspace{3pt} \noindent \Proof Differentiate through formula
\eqref{eq:pabel} for $P_r*f$ and note that
\[
\Big( \frac{\partial^2\ }{\partial r^2} + \frac{1}{r}
\frac{\partial\ }{\partial r} + \frac{1}{r^2} \frac{\partial^2\
}{\partial t^2} \Big) (r^{|j|} e^{ijt}) = 0 .
\]

\vspace{6pt} \noindent 2. The Gauss kernel solves the diffusion
(heat) equation:
\[
w(s,t) = (G_s * f)(t)
\]
solves
\[
w_s=w_{tt}
\]
for $(s,t) \in (0,\infty) \times \T$, with initial value $w(0,t)=f(t)$ in the sense of Theorem~\ref{th:sk}.

\vspace{3pt} \noindent \Proof $G_s(t) \overset{\eqref{eq:G1}}{=}
\sum_{j=-\infty}^\infty e^{-j^2 s} e^{ijt}$ implies
\[
(G_s * f)(t)=\sum_{j=-\infty}^\infty e^{-j^2 s} \fhatj e^{ijt}
\]
by Convolution Fact \eqref{eq:convtrig}. Now differentiate through
the sum and use that
\[
\left( \frac{\partial\ }{\partial s} - \frac{\partial^2\ }{\partial
t^2} \right) (e^{-j^2 s} e^{ijt}) = 0 .
\]

\chapter{Fourier series: summability at a point} \label{ch:su}

\subsubsection*{Goal}

Prove a sufficient condition for summability at a point

\subsubsection*{Reference}

\cite{K} Section I.3

\vspace{18pt}

By Chapter~\ref{ch:fs}, if $f$ is continuous then $\sigma_n(f)\to f$
in $C(\T)$. That is, $\sigma_n(f) \to f$ uniformly. In particular,
$\sigma_n(f)(t)\to f(t)$, for each $t\in \T$.

But what if $f$ is merely continuous at a \emph{point}?

\begin{theorem}[Summability at a point] \label{th:sp}
Assume $\{ k_n \}$ is a summability kernel,$f \in \lot$ and $t_0 \in
\T$. Suppose either $\{ k_n \}$ satisfies the $L^\infty$
concentration hypothesis \eqref{eq:S4}, or else $f \in \lit$.

(a) If $f$ is continuous at $t_0$ then $(k_n * f)(t_0) \to f(t_0)$
as $n \to \infty$.

(b) If in addition the summability kernel is even ($k_n(-t)=k_n(t)$)
and
\[
L= \lim_{h \to 0} \frac{f(t_0+h)+f(t_0-h)}{2}
\]
exists (or equals $\pm \infty$), then
\[
(k_n * f)(t_0) \to L \qquad \text{as $n \to \infty$.}
\]
\end{theorem}
Note if $f$ has limits from the left and right at $t_0$, then the
quantity $L$ equals the average of those limits.

The \Fejer and Poisson kernels satisfy \eqref{eq:SR4}, and so
Theorem~\ref{th:sp} applies in particular to summability at a point
for $\sigma_n(f)=F_n * f$ and for the Abel mean $P_r * f$.
\begin{proof} (a) Let $\e>0$ and choose $0<\delta<\pi$ such that
\begin{equation} \label{eq:sp1}
\sup_n \lv k_n \rv_\lot \cdot \max_{|\tau| \leq \delta} |f(t_0 -
\tau) - f(t_0)| < \e ,
\end{equation}
using here \eqref{eq:S2} and continuity of $f$ at $t_0$. Then as $n
\to \infty$,
\begin{align*}
& |(k_n * f)(t_0)-f(t_0)| \\
& = \Big| \int_{\{ |\tau|<\delta \}} k_n(\tau) [f(t_0-\tau)-f(t_0)] \, d\tau \\
& \qquad - \int_{\{ \delta<|\tau|<\pi \}} k_n(\tau) \, d\tau \cdot
f(t_0) + \int_{\{ \delta<|\tau|<\pi \}} k_n(\tau) f(t_0-\tau) \,
d\tau \Big| \qquad \text{using \eqref{eq:S1}} \\
& <
\begin{cases}
\e + o(1) + o(1) \cdot \lv f \rv_\lot &
\text{by \eqref{eq:sp1}, \eqref{eq:S3} and \eqref{eq:S4}, or else} \\
\e + o(1) + o(1) \cdot \lv f \rv_\lit & \text{by \eqref{eq:sp1},
\eqref{eq:S3} and \eqref{eq:S1},}
\end{cases} \\
& < \e
\end{align*}
for all large $n$.

(b) The proof is similar to (a), but uses symmetry of the kernel.
\end{proof}

\begin{remark}\rm \

\noindent 1. How does the proof of summability at a point, in
Theorem~\ref{th:sp}(a), differ from the proof of summability in norm, in Theorem~\ref{th:sk}?

\noindent 2. Theorem~\ref{th:sp} treats summability at a single
point $t_0$ at which $f$ is continuous. Chapter~\ref{ch:sp} will
prove $k_n * f \to f$ at almost every point, for each integrable
$f$.
\end{remark}

\chapter{Fourier coefficients in $\ell^1(\Z)$ (or, $f \in A(\T)$)} \label{ch:l1}

\subsubsection*{Goal}

Establish the algebra structure of $A(\T)$

\subsubsection*{Reference}

\cite{K} Section I.6

\vspace{18pt}

Define
\begin{align*}
A(\T) & = \{ f \in \lot : \sum_{n \in \Z} |\fhatn| < \infty \} \\
& = \text{functions with Fourier coefficients in $\ell^1(\Z)$.}
\end{align*}
The map $\widehat{\ } : A(\T)\to \ell^1(\Z)$ is a linear bijection.

\Proof Injectivity follows from the uniqueness result \eqref{eq:ut}.
To prove surjectivity, let $\{c_n \} \in \ell^1(\Z)$ and define
$g(t)=\sum_{n \in \Z} c_n e^{int}$. The series for $g$ converges
uniformly since
\[
\sup_{t \in \T} \big| \sum_{|n| > N} c_n e^{int} \big| \leq
\sum_{|n|
> N} |c_n| \to 0
\]
as $N \to \infty$. (Hence $g$ is continuous.) We have $\ghat(m)=c_m$
for every $m$, and so $\ghat=\{ c_m \}$ as desired.

\vspace{6pt} Our proof has shown each $f \in A(\T)$ is represented
by its Fourier series:
\begin{equation} \label{eq:arep}
f(t) = \sum_{n \in \Z} \fhatn e^{int} \quad \text{a.e.}
\end{equation}
so that $f$ is continuous (after redefinition on a set of measure
zero). This Fourier series converges absolutely and uniformly.

\begin{definition}\rm
Define a norm on $A(\T)$ by
\[
\lv f \rv_{A(\T)} = \lv \fhat \rv_{\ell^1(\Z)} = \sum_n |\fhatn| .
\]
$A(\T)$ is a Banach space under this norm (because $\ell^1(\Z)$ is
one).
\end{definition}

Define the convolution of sequences $a,b\in \ell^1(\Z)$ by
\[
(a*b)(n)=\sum_{m \in \Z} a(m)b(n-m).
\]
Clearly
\begin{equation}
\lv a*b \rv_{\ell^1(\Z)} \leq \lv a \rv_{\ell^1(\Z)} \lv b
\rv_{\ell^1(\Z)} \label{eq:seqconv}
\end{equation}
because
\[
\sum_n |(a*b)(n)| \leq \sum_m |a(m)| \sum_n |b(n-m)| = \lv a
\rv_{\ell^1(\Z)} \lv b \rv_{\ell^1(\Z)} .
\]

\begin{theorem}[ $\widehat{\ }$ takes multiplication to convolution]
$A(\T)$ is an algebra, meaning that if $f,g\in A(\T)$ then $fg\in
A(\T)$. Indeed
\[
\widehat{fg}=\fhat*\ghat
\]
and $\lv fg \rv_{A(\T)}\leq \lv f \rv _{A(\T)} \lv g \rv_{A(\T)}$.
\end{theorem}

\begin{proof}
$fg$ is continuous, and hence integrable, with
\begin{align*}
\widehat{(fg)}(n) & = \intp_{\T} f(t)g(t) e^{-int}dt \\
& = \sum_m \fhat(m) \intp_{\T} g(t) e^{-i(n-m)t} \,  dt && \text{by \eqref{eq:arep}} \\
& = \sum_m \fhat(m) \ghat(n-m) \\
& = (\fhat*\ghat)(n) .
\end{align*}
So $\lv \widehat{(fg)} \rv_{\ell^1(\Z)} \leq \lv \fhat
\rv_{\ell^1(\Z)} \lv \ghat \rv_{\ell^1(\Z)}$ by \eqref{eq:seqconv}.
\end{proof}

Sufficient conditions for membership in $A(\T)$ are discussed in
\cite[Section I.6]{K}, for example, H\"older continuity:
$C^\alpha(\T)\subset A(\T)$ when $\alpha>\frac{1}{2}$.

\begin{theorem}[Wiener's Inversion Theorem]
If $f \in A(\T)$ and $f(t) \neq 0$ for every $t \in \T$ then $1/f
\in A(\T)$.
\end{theorem}

We omit the proof. Clearly $1/f$ is continuous, but it is \emph{not} clear
that $\widehat{(1/f)}$ belongs to $\ell^1(\Z)$.

\chapter{Fourier coefficients in $\ell^2(\Z)$ (or, $f \in \ltt$)}
\label{ch:l2}

\subsubsection*{Goal}

Study the Fourier ONB for $\ltt$, using analysis and synthesis operators

\vspace{18pt}

\subsubsection*{Notation and definitions}
Let $H$ be a Hilbert space with inner product $\la u,v \ra$ and norm
$\lv u \rv = \sqrt{\la u,u \ra}$.

Given a sequence $\{ u_n \}_{n \in \Z}$ in $H$, define the
\begin{align*}
\textbf{synthesis operator\quad} S : \ell^2(\Z) & \to H \\
\{ c_n \}_{n \in \Z} & \mapsto \sum_n c_n u_n
\end{align*}
and
\begin{align*}
\textbf{analysis operator\quad } T : H & \to \ell^2(\Z) \\
u & \mapsto \{ \la u,u_n \ra \}_{n \in \Z} .
\end{align*}

\begin{theorem}
If analysis is bounded ($\sum_n |\la u,u_n \ra|^2 \leq
(\text{const.}) \lv u \rv^2$ for all $u \in H$), then so is
synthesis, and the series $S(\{ c_n \})= \sum c_n u_n$ converges
unconditionally.
\end{theorem}

\begin{proof}
Since $T$ is bounded, the adjoint $T^* : \ell^2(\Z) \to H$ is
bounded, and for each sequence $\{c_n \}, u \in H, N \geq 1$, we
have
\begin{align*}
\la T^*( \{ c_n \}_{n=-N}^N) , u \ra
& = \la \{ c_n \}_{n=-N}^N, Tu \ra_{\ell^2} \\
& = \sum_{n=-N}^N c_n \overline{\la u,u_n \ra} \qquad \text{by definition of $Tu$} \\
& = \la \sum_{n=-N}^N c_n u_n , u \ra .
\end{align*}
Hence $T^* (\{ c_n \}_{n=-N}^N) = \sum_{n=-N}^N c_n u_n$. The limit
as $N \to \infty$ exists on the left side, and hence on the right
side; therefore $T^* (\{ c_n\})=\sum_{n=-\infty}^\infty c_n u_n$, so
that $T^*=S$. Hence $S$ is bounded.

Convergence of the synthesis series is unconditional, because if $A
\subset \Z$ then
\begin{align*}
\big\lv S(\{ c_n \}_{n \in \Z}) - S(\{ c_n \}_{n \in A}) \big\rv & =
\big\lv S(\{ c_n \}_{n \in \Z \setminus A}) \big\rv \\
& \leq \lv S \rv \lv \{ c_n \} \rv_{\ell^2(\Z \setminus A)} ,
\end{align*}
which tends to $0$ as $A$ expands to fill $\Z$, regardless of the
order in which $A$ expands.
\end{proof}

\begin{remark}\rm
The last proof shows $S=T^*$, meaning
\[
\fbox{\text{analysis and synthesis are adjoint operations.}}
\]
\end{remark}

\begin{theorem}[Fourier coefficients on $\ltt$] \label{th:fsl2}
The Fourier coefficient (or analysis) operator $\widehat{\ } : \ltt \to \ell^2(\Z)$ is an
isometry, with
\begin{align*}
\lv f \rv_\ltt & = \lv \fhat \rv_{\ell^2(\Z)} &&
\textsc{(Plancherel)} \\
\la f , g \ra_{\ltt} & = \la \fhat , \ghat \ra_{\ell^2(\Z)} &&
\textsc{(Parseval)}
\end{align*}
for all $f,g \in \ltt$.
\end{theorem}

\begin{proof}
First we prove Plancherel's identity: since $P_r * f \to f$ in $\ltt$ by Theorem~\ref{th:sk}, we have
\begin{align*}
\intp_\T |f(t)|^2 \, dt & = \lim_{r \to 1} \intp_\T f(t)
\overline{(P_r * f)(t)} \, dt \\
& = \lim_{r \to 1} \intp_\T f(t) \sum_{n \in \Z} \overline{r^{|n|} \fhatn e^{int} }
\, dt \qquad \text{by \eqref{eq:pabel} for $P_r * f$} \\
& = \lim_{r \to 1} \sum_{n \in \Z} r^{|n|} |\fhatn|^2 \\
& = \sum_{n \in \Z} |\fhatn|^2
\end{align*}
by monotone convergence.

Parseval follows from Plancherel by polarization, or by repeating
the argument for Plancherel with $\la f ,f \ra$ changed to $\la f ,
g \ra$ (and using dominated instead of monotone convergence).
\end{proof}

Since the Fourier analysis operator is bounded, so is its adjoint, the Fourier synthesis operator
\begin{align*}
\check{\ } : \ell^2(\Z) & \to \ltt \\
\{ c_n \}_{n \in \Z} & \mapsto \sum_n c_n e^{int}
\end{align*}

\begin{theorem}[Fourier ONB] \label{th:fonb} \

\noindent (a) If $f \in \ltt$ then $\sum_n \fhatn e^{int} = f$ with unconditional convergence
in $\ltt$. That is, $(\hat{f})\check{\ }=f$.

\noindent (b) If $c = \{ c_n \} \in \ell^2(\Z)$ then $(\ \sum_{n \in \Z} c_n e^{int})\widehat{\ }(j) = c_j$. That is, $(\check{c})\hat{\ }=c$.

\noindent (c) $\{ e^{int} \}_{n \in \Z}$\ is an orthonormal basis of $\ltt$.
\end{theorem}

Part (a) says Fourier series converge in $\ltt$. Parts (a) and (b) together show that Fourier analysis and synthesis are inverse operations.

\begin{proof}
Fourier analysis and synthesis are bounded operators, and analysis followed by synthesis equals the identity ($\sum_n \fhatn e^{int} = f$) on the class of trigonometric polynomials. That class is dense in $\ltt$, and so by continuity, analysis followed by synthesis equals the identity on $\ltt$.

Argue similarly for part (b), using the dense class of finite sequences in $\ell^2(\Z)$.

For orthonormality in part (c), observe
\[
\la e^{int}, e^{imt} \ra = \intp_\T e^{int} e^{-imt} \, dt \\
= \begin{cases}
1 & \text{if $n=m$,} \\
0 & \text{if $n \neq m$.}
\end{cases}
\]
The basis property follows from part (a), noting $\fhatn = \la f ,
e^{int} \ra_\ltt$.
\end{proof}

\begin{remark}\rm Fourier analysis satisfies
\begin{align*}
\widehat{\ } & : \lot \to \ell^\infty(\Z) && \text{by Theorem~\ref{th:bp},} \\
\widehat{\ }& : \ltt \to \ell^2(\Z) && \text{(isometrically) by
Theorem~\ref{th:fsl2}.}
\end{align*}
Further, $\widehat{\ } : \ltt \to \ell^2(\Z)$ is a linear bijection by Theorem~\ref{th:fonb}.

In Chapter~\ref{ch:ai} we will interpolate to show
\[
\widehat{\ }: \lpt \to \ell^{p^\prime}(\Z), \qquad \text{whenever\ }
1 \leq p \leq 2, \quad \frac{1}{p}+\frac{1}{p^\prime}=1 .
\]
\end{remark}

\chapter{Maximal functions} \label{ch:mf}

\subsubsection*{Goals}

Connect abstract maximal functions to convergence a.e.

\noindent Prove weak and strong bounds on the Hardy--Littlewood
maximal function

\noindent Prepare for summability pointwise a.e.\ in next Chapter

\subsubsection*{References}

\cite{D} Section 2.2

\noindent \cite{G} Section 2.1

\noindent \cite{S} Section 1.1

\vspace{18pt}

\begin{definition}[Weak and strong operators]
Let $(X,\mu)$ and $(Y,\nu)$ be measure spaces, and $1 \leq p,q \leq
\infty$. Suppose
\[
T : \lp(X) \to \{ \text{measurable functions on $Y$} \} .
\]
(We do not assume $T$ is linear.)

Call $T$ \textbf{strong} $\mathbf{(p,q)}$ if $T$ is bounded from
$\lp(X)$ to $L^q(Y)$, meaning a constant $C > 0$ exists such that
\[
    \lv Tf \rv_{L^q(Y)} \leq C \lv f \rv _{\lp(X)}, \qquad f \in
    \lp(X).
\]
When $q<\infty$, we call $T$ \textbf{weak} $\mathbf{(p,q)}$ if $C
> 0$ exists such that
\[
\nu \big( \{ y \in Y : |(Tf)(y)| > \lambda \} \big)^{1/q} \leq
\frac{C \lv f \rv _{\lp(X)}}{\lambda} \qquad \forall \ \lambda
> 0, \quad f \in \lp(X) .
\]
When $q=\infty$, we call $T$ \textbf{weak} $\mathbf{(p,\infty)}$ if
it is strong $(p,\infty)$:
\[
    \lv Tf \rv_{L^\infty(Y)} \leq C \lv f \rv _{\lp(X)}, \qquad f \in
    \lp(X).
\]
\end{definition}

\begin{lemma}
Strong $(p,q)$ \ $\Rightarrow$ \ weak $(p,q)$.
\end{lemma}

\begin{proof}
When $q=\infty$ the result is immediate by definition. Suppose
$q<\infty$. Write
\[
E(\lambda) = \{ y \in Y: |(Tf)(y)| > \lambda \}
\]
for the level set of $Tf$ above height $\lambda$. Then
\begin{align*}
\lambda ^q \, \nu \big( E(\lambda) \big) & = \int_{E(\lambda)} \lambda^q \, d\nu(y) \\
& \leq \int_{E(\lambda)} |(Tf)(y)|^q \, d\nu(y) && \text{since $\lambda < |Tf|$ on $E(\lambda)$} \\
& \leq \lv Tf \rv ^q_{L^q(Y)}
\end{align*}
and so
\begin{align*}
\nu \big( E(\lambda) \big)^{1/q} & \leq \frac{\lv Tf
\rv_{L^q(Y)}}{\lambda} \\
& \leq \frac{C \lv f \rv_{\lp(X)}}{\lambda}
\end{align*}
if T is strong $(p,q)$.
\end{proof}

\begin{lemma}
If $T$ is weak $(p,q)$ then $Tf \in L^r_{loc}(Y)$ for all $0<r<q$.
\end{lemma}

Thus intuitively, $T$ ``almost'' maps $L^p$ into $L^q$, locally.

\begin{proof}
Let $f \in \lp(X)$ and suppose $Z \subset Y$ with $\nu (Z) <
\infty$. We will show $Tf \in L^r(Z)$.

Write $g = Tf$. Then
\begin{align*}
& \int_Z |g(y)|^r \, d\nu(y) \\
& = \int_0^\infty r \lambda^{r-1} \nu \big( \{y \in Z : |g(y)|>\lambda \} \big) \, d\lambda && \text{by Appendix\ref{ap:df}} \\
& \leq \int_0^1 r \lambda^{r-1} \nu(Z) \, d\lambda + \int_1^\infty r
\lambda^{r-1} \Big( \frac{C \lv f \rv_{\lp(X)}}{\lambda} \Big)^{\!
\! q}
\, d\lambda && \text{by weak $(p,q)$} \\
& < \infty
\end{align*}
since $\nu(Z)<\infty$ and $\int_1^\infty \lambda^{-1-q+r} \,
d\lambda < \infty$ (using that $-q+r < 0$).
\end{proof}

\begin{theorem}[Maximal functions and convergence a.e.] \label{th:mfcae}
Assume
\[
T_n : L^p(X) \to \{ \text{measurable functions on $X$} \}
\]
for $n=1,2,3,\ldots$. Define
\[
T^* : L^p(X) \to \{ \text{measurable functions on $X$} \}
\]
by
\[
(T^*f)(x) = \sup_n |(T_n f)(x)| , \qquad x \in X .
\]

If $T^*$ is weak $(p,q)$ and each $T_n$ is linear, then the
collection
\[
{\mathcal C} = \{ f \in L^p(X) : \lim_n (T_n f)(x) = f(x) \
\text{a.e} \}
\]
is closed in $L^p(X)$.
\end{theorem}

$T^*$ is called the \textbf{maximal operator} for the family $\{ T_n
\}$. Clearly it takes values in $[0,\infty]$. Note $T^*$ is not
linear, in general.

\begin{remark}\rm
In this theorem a \emph{quantitative} hypothesis (weak $(p,q)$)
implies a \emph{qualitative} conclusion (closure of the collection
$\mathcal C$ where $T_n f \to f$ a.e.).
\end{remark}

\begin{proof}
Let $f_k \in {\mathcal C}$ with $f_k \to f$ in $L^p(X)$. We show $f
\in {\mathcal C}$.

Suppose $q<\infty$. For any $\lambda>0$,
 \begin{align*}
& \mu \big( \{x \in X: \limsup_n |(T_n f(x)-f(x)|>2\lambda \} \big) \\
& = \mu \big( \{ x \in X : \limsup_n |T_n(f-f_k)(x)-(f-f_k)(x)|>2\lambda \} \big) \\
& \qquad \qquad \text{by linearity and the pointwise convergence $T_n f_k \to f_k$ a.e.} \\
& \leq \mu \big( \{x \in X : T^*(f-f_k)(x) + |(f-f_k)(x)|>2\lambda \} \big) \qquad \text{by triangle inequality} \\
& \leq \mu \big( \{ x \in X : T^*(f-f_k)(x)>\lambda \} \big) + \mu \big( \{ x \in X : |(f-f_k)(x)|>\lambda \} \big) \\
& \leq \Big( \frac{C \lv f-f_k \rv_{L^p(X)}}{\lambda} \Big)^{\! q} +
\Big( \frac{\lv f-f_k \rv_{L^p(X)}}{\lambda} \Big)^{\! p}
\qquad \text{by weak $(p,q)$ on $T^*$} \\
& \to 0
\end{align*}
as $k \to \infty$.

Therefore $\limsup_n |(T_n f)(x) - f(x)| \leq 2\lambda$ a.e. Taking
a countable sequence of $\lambda \searrow 0$, we conclude $\limsup_n
|(T_n f)(x) - f(x)| = 0$ a.e. Therefore $\lim_n (T_nf)(x)=f(x)$
a.e.,  so that $f \in {\mathcal C}$.

The case $q=\infty$ is left to the reader.
\end{proof}

To apply maximal functions on $\Rd$ and $\T$, we will need:
\begin{lemma}[Covering] \label{le:cl}
Let $\{B_i\}_{i=1}^k$ be a finite collection of open balls in $\Rd$.
Then there exists a pairwise disjoint subcollection
$\{B_{i_j}\}_{j=1}^l$ of balls such that
\[
\big| \bigcup_{i=1}^k B_i \big| \leq 3^d \big| \bigcup_{j=1}^l
B_{i_j} \big| = 3^d \sum_{j=1}^l |B_{i_j}| .
\]
\end{lemma}

Thus the subcollection covers at least $1/3^d$ of the total volume
of the balls.

\begin{proof}
Re-label the balls in decreasing order of size: $|B_1| \geq |B_2|
\geq \cdots \geq |B_k|$. Choose $i_1=1$ and employ the following
greedy algorithm. After choosing $i_j$, choose $i_{j+1}$ to be the
smallest index $i>i_j$ such that $B_i$ is disjoint from $B_{i_1},
\ldots, B_{i_j}$. Continue until no such ball $B_i$ exists.

The $B_{i_j}$ are pairwise disjoint, by construction.

Let $i \in \{ 1,\ldots, k \}$. If $B_i$ is \emph{not} one of the
$B_{i_j}$ chosen, then $B_i$ must intersect one of the $B_{i_j}$ and
be smaller than it, so that
\[
\text{radius} \, (B_i) \leq \text{radius}(B_{i_j}) .
\]
Hence $B_i \subset 3B_{i_j}$ (where we mean the ball with the same
center and three times the radius). Thus
\begin{align*}
\big| \bigcup_{i=1}^k B_i \big| & \leq \big| \bigcup_{j=1}^l
(3B_{i_j})
\big| \\
& \leq \sum_{j=1}^l |3B_{i_j}| \\
& = 3^d \sum_{j=1}^l |B_{i_j}| \\
& = 3^d \big| \bigcup_{j=1}^l B_{i_j} \big|
\end{align*}
by disjointness of the $B_{i_j}$.
\end{proof}

\begin{definition}
The \emph{Hardy--Littlewood (H-L) maximal function} of a locally
integrable function $f$ on $\Rd$ is
\begin{align*}
(Mf)(x) & = \sup_{r>0} \frac{1}{|B(x,r)|} \int_{B(x,r)} |f(y)| \,
dy \\
& = \text{``largest local average'' of $|f|$ around $x$.}
\end{align*}
\end{definition}

\subparagraph*{Properties}
\begin{align*}
Mf & \geq 0 \\
|f| \leq |g| \ \Rightarrow \ Mf & \leq Mg \\
M(f+g) & \leq Mf+Mg \quad \text{(sub-linearity)} \\
Mc & = c \quad \text{if $c = \text{(const.)} \geq 0$}
\end{align*}

\begin{theorem}[H-L maximal operator]
$M$ is weak $(1,1)$ and strong $(p,p)$ for $1 < p \leq \infty$.
\end{theorem}

\begin{proof}
For weak $(1,1)$ we show
\begin{equation} \label{eq:hl1}
        |E(\lambda)| \leq \frac{3^d\lv f \rv_{\lord}}{\lambda}
\end{equation}
where $E(\lambda) = \{ x \in \Rd : Mf(x)>\lambda \}$. If $x \in
E(\lambda)$ then
\[
\frac{1}{|B(x,r)|} \int_{B(x,r)} |f(y)|\, dy > \lambda
\]
for some $r>0$. The same inequality holds for all $x^\prime$ close
to $x$, so that $x^\prime \in E(\lambda)$. Thus $E(\lambda)$ is open
(and measurable), and $Mf$ is lower semicontinuous (and measurable).

Let $F \subset E(\lambda)$ be compact. Each $x \in F$ is the center
of some ball $B$ such that
\begin{equation} \label{eq:hl2}
|B| < \frac{1}{\lambda} \int_B |f(y)| \, dy .
\end{equation}
By compactness, $F$ is covered by finitely many such balls, say
$B_1, \ \ldots,\ B_k$. The Covering Lemma~\ref{le:cl} yields a
subcollection $B_{i_1}, \ldots, B_{i_l}$. Then
\begin{align*}
|F| & \leq \big| \bigcup_{i=1}^k B_i \big| \\
& \leq 3^d \sum_{j=1}^l |B_{i_j}| && \text{by Covering Lemma~\ref{le:cl}} \\
& \leq \frac{3^d}{\lambda} \sum_{j=1}^l \int_{B_{i_j}} |f(y)| \, dy && \text{by \eqref{eq:hl2}} \\
& \leq \frac{3^d}{\lambda} \int_\Rd |f(y)| \, dy && \text{by disjointness} \\
& = \frac{3^d}{\lambda} \lv f \rv_{\lord} .
\end{align*}
Taking the supremum over all compact $F \subset E(\lambda)$ gives
\eqref{eq:hl1}.

For strong $(\infty,\infty)$, note $Mf(x) \leq \lv f \rv_{\lird}$
for all $x \in \Rd$, by definition of $Mf$. Hence $\lv Mf \rv_\lird
\leq \lv f \rv_\lird$.

For strong $(p,p)$ when $1<p<\infty$, let $\lambda>0$ and define
\begin{align*}
g(x) = \begin{cases} f(x) & \text{if $|f(x)| > \lambda/2$} \\
0 & \text{otherwise}
\end{cases}
\ = \text{``large'' part of $f$,} \\
h(x) = \begin{cases} f(x) & \text{if $|f(x)| \leq \lambda/2$} \\
0 & \text{otherwise}
\end{cases}
\ = \text{``small'' part of $f$.}
\end{align*}
Then $f=g+h$ and $|h| \leq \lambda/2$, so that $Mf \leq Mg +
\lambda/2$. Hence
\begin{align}
|E(\lambda)| & = \big| \{ x : Mf(x) > \lambda \} \big| \notag \\
& \leq \big| \{ x : Mg(x) > \lambda/2 \} \big| \notag \\
& \leq \frac{3^d \lv g \rv_{\lord}}{\lambda/2} && \text{by the above weak $(1,1)$ result} \notag \\
& = \frac{2 \cdot 3^d}{\lambda} \int_{ \{ x : |f(x)|>\lambda/2 \}}
|f(x)| \, dx . \label{eq:hl3}
\end{align}
Therefore
\begin{align*}
& \int_\Rd |Mf(x)|^p \, dx \\
& = \int_0^\infty p\lambda^{p-1} |E(\lambda)| \, d\lambda && \text{by Appendix~\ref{ap:df}} \\
& \leq 2\cdot 3^d p \int_0^\infty \lambda^{p-2} \int_{ \{ x : |f(x)|>\lambda/2 \} } |f(x)| \, dx d\lambda && \text{by \eqref{eq:hl3}} \\
& = \Big( 2^p 3^d \frac{p}{p-1} \Big) \int_\Rd |f(x)|^p \, dx
\end{align*}
by Lemma~\ref{le:df} with $r=1,\alpha=2$. We have proved the strong
$(p,p)$ bound.
\end{proof}

Notice the constant in the strong $(p,p)$ bound blows up as $p
\searrow 1$. As this observation suggests, the Hardy--Littlewood
maximal operator is not strong $(1,1)$. For example, the indicator
function $f=\charfn_{[-1,1]}$ in $1$ dimension has $Mf(x) \sim
c/|x|$ when $|x|$ is large, so that $Mf \notin \lor$.

The maximal function is \emph{locally} integrable provided $f \in L
\log L(\Rd)$; see Problem~\ref{pr:llogl}.

\chapter[Fourier summability pointwise a.e.]{Fourier series: summability pointwise a.e.} \label{ch:sp}

\subsubsection*{Goal} Prove summability a.e.\ using Fej\'{e}r and Poisson maximal functions

\vspace{18pt}

\begin{definition}
\begin{align*}
\text{Dirichlet maximal function} \quad (D^*f)(t) & = \sup_n |(D_n * f)(t)| = \sup_n |S_n(f)(t)| \\
\text{\Fejer maximal function} \quad (F^*f)(t) & = \sup_n |(F_n * f)(t)| = \sup_n |\sigma_n(f)(t)| \\
\text{Poisson maximal function} \quad (P^*f)(t) & = \sup_{0<r<1} |(P_r * f)(t)| \\
\text{Gauss maximal function} \quad (G^*f)(t) & = \sup_{0<s<\infty} |(G_s * f)(t)| \\
\text{Lebesgue maximal function} \quad (L^*f)(t) & = \sup_{0<h<\pi} |(L_h * f)(t)|
\end{align*}
where the Lebesgue kernel is $L_h(t) = 2\pi \frac{1}{2h} \charfn_{[-h,h]}(t)$, extended $2\pi$-periodically. Notice $(L_h * f)(t) = \frac{1}{2h} \int_{t-h}^{t+h} f(\tau) \, d\tau$ is a local average of $f$ around $t$.
\end{definition}

\begin{lemma}[Majorization] \label{le:maj}
If $k \in \lot$ is nonnegative and symmetric ($k(-t)=k(t)$), and decreasing on $[0,\pi]$, then
\[
|(k*f)(t)| \leq \lv k \rv_\lot (L^*f)(t) \qquad \text{for all $t
\in \T, \quad f \in \lot$.}
\]
\end{lemma}
Thus convolution with a symmetric decreasing kernel is majorized by the Hardy--Littlewood maximal function.
\begin{proof}
Assume $k$ is absolutely continuous, for simplicity. We first establish a ``layer cake'' decomposition of
$k$, representing it as a linear combination of kernels $L_h$:
\begin{align*}
k(t) = k(|t|) & = k(\pi) - \int_{|t|}^\pi k^\prime(h) \, dh \\
& = k(\pi) - \pr \int_0^\pi 2h L_h(t)
k^\prime(h) \, dh ,
\end{align*}
since
\[
\pr 2h L_h(t) =
\begin{cases}
1 , & \text{if $h \geq |t|$,} \\
0, & \text{if $h<|t|$.}
\end{cases}
\]
Hence
\begin{align*}
(k*f)(t) & = k(\pi) \intp_\T f(\tau) \, d\tau + \pr \int_0^\pi 2h (L_h * f)(t) \big( -k^\prime(h) \big) \, dh \\
|(k*f)(t)| & \leq k(\pi) |(L_\pi * f)(t)| + \pr \int_0^\pi 2h \big( -k^\prime(h) \big) \, dh \, (L^* f)(t) \\
& \qquad \qquad \qquad \text{using $k(\pi) \geq 0$ and $k^\prime \leq 0$} \\
& \leq \frac{2}{2\pi} \int_0^\pi k(h) \, dh \, (L^* f)(t) \qquad \text{by parts} \\
& = \lv k \rv_\lot (L^* f)(t)
\end{align*}
by symmetry of $k$.
\end{proof}

\begin{theorem}[Lebesgue dominates \Fejer and Poisson] \label{th:lfp}
For all $f \in \lot$,
\begin{align*}
F^* f & \leq 2 L^*|f| \\
P^* f & \leq \ L^* f
\end{align*}
\end{theorem}

\begin{proof}
$P_r(t)$ is nonnegative, symmetric, and decreasing on $[0,\pi]$ (exercise), with $\lv P_r \rv_\lot = 1$. Hence $|P_r * f| \leq L^* f$ by Majorization Lemma~\ref{le:maj}, so that $P^* f \leq L^* f$.

The \Fejer kernel is \emph{not} decreasing on $[0,\pi]$, but it is bounded by a symmetric decreasing kernel, as follows:
\begin{align*}
F_n(t) & = \frac{1}{n+1} \left( \frac{\sin \big( \frac{n+1}{2} t \big)}{\sin \big( \frac{1}{2} t \big)} \right)^{\! \! 2} \\
& \leq k(t) \overset{\text{def}}{=} \frac{1}{n+1}
\begin{cases}
(n+1)^2 & \text{if $|t| \leq \pi/(n+1)$,} \\
\pi^2/t^2 & \text{if $\pi/(n+1) \leq |t| \leq \pi$,}
\end{cases}
\end{align*}
since
\begin{align*}
\sin(n+1)\theta & \leq (n+1) \sin \theta , && 0 \leq \theta \leq \frac{\pi}{2} , \\
\sin \big( \frac{1}{2} t \big) & \geq \frac{t}{\pi} , && 0 \leq t \leq \pi .
\end{align*}
Note the kernel $k$ is nonnegative, symmetric, and decreasing on $[0,\pi]$, with
\[
\lv k \rv_\lot = \pr \big( 4\pi - \frac{2\pi}{n+1} \big) < 2 .
\]
Hence $|F_n * f| \leq k*|f| \leq 2L^* |f|$ by Majorization Lemma~\ref{le:maj}, so that $F^* f \leq 2 L^* |f|$.
\end{proof}

The Gauss kernel can be shown to be symmetric decreasing, so that $G^* f \leq \ L^* f$, but we omit the proof.

\begin{corollary} \label{co:kplweak}
$F^*, P^*$ and $L^*$ are weak $(1,1)$ on $\T$.
\end{corollary}

\begin{proof}
\[
\big| \{ t \in \T : (L^* f)(t) > 2\pi \lambda \} \big| \leq \big| \{ t \in \T : (L^* |f|)(t) > 2\pi \lambda \} \big| \leq \frac{3}{\lambda} \int_\T |f(t)| \, dt
\]
by repeating the weak $(1,1)$ proof for the Hardy--Littlewood maximal function. These weak $(1,1)$ estimates for $L^* f$ and $L^* |f|$ imply weak $(1,1)$ for $F^*f$, since if $(F^* f)(t)>\lambda$ then $(L^* |f|)(t) > \lambda/2$ by  Theorem~\ref{th:lfp}. Argue similarly for $P^*f$.
\end{proof}

\begin{theorem}[Summability a.e.] \label{th:sae}
If $f \in \lot$ then
\begin{align*}
\sigma_n(f) = F_n * f \to f \ \text{a.e.} & \quad \text{as $n \to \infty$} && \text{(\Fejer summability)} \\
              P_r * f \to f \ \text{a.e.} & \quad \text{as $r \nearrow 1$} && \text{(Abel summability)} \\
              L_h * f \to f \ \text{a.e.} & \quad \text{as $h \searrow 0$} && \text{(Lebesgue differentiation theorem)}
\end{align*}
\end{theorem}

\begin{proof}
By the weak $(1,1)$ estimate in Corollary~\ref{co:kplweak} and the abstract convergence result in Theorem~\ref{th:mfcae}, the set
\[
{\mathcal C} = \{ f \in \lot : \lim_n (F_n * f)(t) = f(t) \ \text{a.e.} \}
\]
is closed in $\lot$.

Obviously ${\mathcal C}$ contains the continuous functions on $\T$, since $F_n * f \to f$ uniformly when $f$ is continuous. Thus ${\mathcal C}$ is dense in $\lot$. Because ${\mathcal C}$ is also closed, it must equal $\lot$, thus proving \Fejer summability a.e.\ for each $f \in \lot$.

Argue similarly for $P_r * f$ and $L_h * f$.
\end{proof}

The result that $L_h * f \to f$ a.e.\ means
\[
\frac{1}{2h} \int_{t-h}^{t+h} f(\tau) \, d\tau \to f(t) \quad \text{a.e.,}
\]
which is the Lebesgue differentiation theorem on $\T$.

\chapter{Fourier series: convergence at a point} \label{ch:cp}

\subsubsection*{Goals}

State divergence pointwise can occur for $\lot$

\noindent Show divergence pointwise can occur for $C(\T)$

\noindent Prove convergence pointwise for $C^\alpha(\T)$ and
$BV(\T)$

\subsubsection*{References}

\cite{K} Section II.2, II.3

\noindent \cite{D} Section 1.1

\vspace{18pt}

Fourier series can behave badly for integrable functions.

\begin{theorem}[Kolmogorov] \label{th:kol}
There exists $f \in \lot$ whose Fourier series diverges
unboundedly at \emph{every} point. That is,
\[
\sup_n |S_n(f)(t)| = \infty \qquad \text{for all $t \in \T$,}
\]
so that $D^* f \equiv \infty$.
\end{theorem}
Recall $S_n(f)=D_n * f$ and $D^* f$ is the maximal function for the
Dirichlet kernel.
\begin{proof}
\cite[Section II.3]{K}.
\end{proof}

Even continuous functions can behave badly.

\begin{theorem}
There exists a continuous function whose Fourier series diverges
unboundedly at $t=0$. That is,
\[
\sup_n |S_n(f)(0)|=\infty .
\]
\end{theorem}

\begin{proof}
Define
\begin{align*}
    T_n : C(\T) & \to \C \\
     f & \mapsto S_n(f)(0)=\text{($n$th partial sum of $f$ at $t=0$).}
\end{align*}
Then $T_n$ is linear. Each $T_n$ is bounded since
\begin{align*}
|T_n (f)| & = |S_n(f)(0)| \\
& = |(D_n * f)(0)| \\
& = \big| \intp_\T D_n(\tau) f(0-\tau) \, d\tau \big| \\
& \leq \lv D_n \rv_\lot \lv f \rv_\lit .
\end{align*}
Thus $\lv T_n \rv \leq \lv D_n \rv_\lot$. We show $\lv T_n \rv = \lv
D_n \rv_\lot$. Let $\e > 0$ and choose $g \in C(\T)$ with $\lv g
\rv_\lit = 1$ and $g$ even and
\[
g(t)= \begin{cases}
1 & \text{if $D_n(t)>0$,} \\
-1 & \text{if $D_n(t)<0$,} \\
& \text{except for small intervals around the zeros of $D_n$,} \\
& \text{with total length of those intervals $<\e/(2n+1)$.}
\end{cases}
\]
Then
\begin{align*}
|T_n(g)| & = \big| \intp_\T D_n(\tau) g(\tau) \, d\tau \big| \\
& \geq \intp_{\T \setminus \{\text{intervals}\}} |D_n(\tau)| \, d\tau - \intp_{\{\text{intervals}\}} |D_n(\tau)| \, d\tau \\
& = \intp_\T |D_n(\tau)| \, d\tau - \frac{2}{2\pi} \int_{\{\text{intervals}\}} |D_n(\tau)| \, d\tau \\
& \geq \lv D_n \rv_\lot - \frac{1}{\pi} \frac{\e}{2n+1}(2n+1) \qquad \text{using definition \eqref{eq:D1} of $D_n$} \\
& = \lv D_n \rv_\lot - \frac{\e}{\pi} \\
& = \big( \lv D_n \rv_\lot - \frac{\e}{\pi} \big) \lv g \rv_\lit .
\end{align*}
Thus $\lv T_n \rv \geq \lv D_n \rv_\lot - \e/\pi$ for all $\e>0$,
and so $\lv T_n \rv = \lv D_n \rv_\lot$.

Recalling that $\lv D_n \rv_\lot \to \infty$ as $n \to \infty$ (in
fact, $\lv D_n \rv \sim c \log n$ by \cite[Ex.~II.1.1]{K}) we
conclude from the Uniform Bounded Principle (Banach--Steinhaus) that
there exists $f \in C(\T)$ with $\sup_n |T_n(f)|=\infty$, as
desired.
\end{proof}

\noindent \emph{Another proof.} \cite[Sec II.2]{K} gives an explicit
construction of $f$, proving divergence not only at $t=0$ but on a
dense set of $t$-values.

\vspace{12pt} Now we prove convergence results.

\begin{theorem}[Dini's Convergence Test] \label{th:dini}
Let $f \in \lot, t \in \T$. If
\[
\int_{-\pi}^{\pi} \left| \frac{f(t-\tau)-f(t)}{\tau} \right| d\tau
< \infty
\]
then the Fourier series of $f$ converges at $t$ to $f(t)$.
\end{theorem}

\begin{proof}
\begin{align}
S_n(f)(t)-f(t) & = \intp_{-\pi}^\pi [f(t-\tau)-f(t)] \frac{\sin \big( (n+\frac{1}{2}) \tau \big)}{\sin \big( \frac{1}{2}\tau \big)} \, d\tau \notag \\
& \qquad \qquad \qquad \text{using that $\intp_{-\pi}^\pi D_n(\tau) \, d\tau=1$} \notag \\
& = \intp_{-\pi}^\pi \left\{ \frac{f(t-\tau)-f(t)}{\tau}
\frac{\tau}{\sin \big( \frac{1}{2}\tau \big)}
\cos \big( \frac{1}{2}\tau \big) \right\} \sin(n\tau) \, d\tau \notag \\
& \qquad \qquad + \intp_{-\pi}^\pi [f(t-\tau)-f(t)] \cos(n\tau) \,
d\tau \label{eq:dini1}
\end{align}
by expanding $\sin \big( (n+\frac{1}{2}) \tau \big)$ with a
trigonometric identity.

Notice the factor $\{\cdots\}$ is integrable with respect to
$\tau$, by the Dini hypothesis. And $\tau \mapsto
[f(t-\tau)-f(t)]$ is integrable too. Hence both integrals in
\eqref{eq:dini1} tend to $0$ as $n \to \infty$, by the
Riemann--Lebesgue Corollary \ref{co:rl} (after expressing
$\sin(n\tau)$ and $\cos(n\tau)$ in terms of $e^{\pm i n\tau}$).
\end{proof}

\begin{corollary}[Convergence for H\"{o}lder continuous $f$]
If $f \in C^{\alpha}(\T), 0 < \alpha \leq 1$, then the Fourier
series of $f$ converges to $f(t)$, for every $t \in \T$.
\end{corollary}

\begin{proof}
Put H\"{o}lder into Dini:
\[
\int_{-\pi}^\pi \left| \frac{f(t-\tau)-f(t)}{\tau} \right| d\tau
\leq \int_{-\pi}^\pi \frac{(\text{const.})|\tau|^{\alpha}}{|\tau|}
\, d\tau < \infty .
\]
Now apply Dini's Theorem~\ref{th:dini}.

(Exercise. Prove the Fourier series in fact converges uniformly.)
\end{proof}

\begin{corollary}[Localization Principle]
Let $f \in \lot, t \in \T$. If $f$ vanishes on a neighborhood of
$t$, then $S_n(f)(t) \to 0$ as $n \to \infty$.
\end{corollary}

\begin{proof}
Apply Dini's Theorem~\ref{th:dini}.
\end{proof}

In particular, if two functions agree on a neighborhood of $t$ and
the Fourier series of one of them converges at $t$, then the
Fourier series of the other function converges at $t$ to the same
value. Thus Fourier series depend only on local information.

\begin{theorem}[Convergence for bounded variation f]
If $f \in BV(\T)$ then the Fourier series converges everywhere to
$\frac{1}{2}[f(t+)+f(t-)]$, and hence converges to $f(t)$ at every
point of continuity.
\end{theorem}

\begin{proof}
Let $t \in \T$. On the interval $(t-\pi, t+\pi)$, express $f$ as
the difference of two bounded increasing functions, say $f=g-h$.
It suffices to prove the theorem for $g$ and $h$ individually.

We have
\begin{align}
S_n(g)(t)-\frac{1}{2}[g(t+)+g(t-)]
& = \intp_0^\pi \big( g(t-\tau)-g(t-) \big) D_n(\tau) \, d\tau \label{eq:bv1} \\
& + \intp_0^\pi \big( g(t+\tau)-g(t+) \big) D_n(\tau) \, d\tau
\label{eq:bv2}
\end{align}
since $D_n(\tau)$ is even and $\intp_0^\pi D_n(\tau) \, d\tau =
\frac{1}{2}$.

Let $G(\tau)=g(t+\tau)-g(t+)$ for $\tau \in (0, \pi)$, so that $G$
is increasing with $G(0+)=0$. Write
\[
H_n(\tau)=\int_0^\tau D_n(\sigma) \, d\sigma
\]
so that $H_n^{\prime} = D_n$. Let $0<\delta<\pi$. Then
\begin{align*}
\eqref{eq:bv2} & = \intp_0^\delta G(\tau) H_n^\prime(\tau) \, d\tau +
\intp_\delta^\pi G(\tau) D_n(\tau) \, d\tau \\
& = \pp G(\delta)H_n(\delta) - \intp_{(0,\delta]} H_n(\tau) \,
dG(\tau) + o(1)
\end{align*}
as $n \to \infty$, by parts in the first term and by the
Localization Principle in the last term, since the function
\[
\begin{cases}
G(\tau), & \delta<\tau<\pi, \\
0 , & -\delta<\tau<\delta, \\
G(-\tau), & -\pi<\tau<-\delta,
\end{cases}
\]
vanishes near the origin. Hence
\begin{align*}
 \limsup_n |\eqref{eq:bv2}| & \leq \pp \sup_n \lv H_n \rv_\lit \left( G(\delta) + \int_{(0,\delta]} \, dG(\tau) \right) \\
 & = \pp \sup_n \lv H_n \rv_\lit \cdot 2G(\delta) \qquad \text{since $G(0+)=0$} \\
 & \to 0
\end{align*}
as $\delta \to 0$. Therefore $\eqref{eq:bv2} \to 0$ as $n \to
\infty$. Argue similarly for \eqref{eq:bv1}, and for $h$.

Thus we are done, provided we show
\[
\sup_n \lv H_n \rv_\lit < \infty .
\]
We have
\begin{align*}
|H_n(\tau)| & \leq \Big| \int_0^{\tau} \frac{\sin \big(
(n+\frac{1}{2}) \sigma \big)}{\frac{1}{2}\sigma} \, d\sigma
\Big| + \Big| \int_0^\tau \sin \big( (n+\frac{1}{2}) \sigma \big) \Big( \frac{1}{\sin \big(\frac{1}{2} \sigma \big)} - \frac{1}{\frac{1}{2} \sigma} \Big) d\sigma \Big| \\
& \leq 2 \Big| \int_0^{(n+\frac{1}{2}) \tau}
\frac{\sin \sigma}{\sigma} \, d\sigma \Big| + \int_0^\pi \text{(const.)} \frac{\sigma^3}{\sigma^2} \, d\sigma \qquad \text{by a change of variable} \\
& \leq 2 \sup_{\rho>0} \Big| \int_0^\rho \frac{\sin \sigma}{\sigma} \, d\sigma \Big| + \text{(const.)} \\
& < \infty
\end{align*}
since $\lim_{\rho \to \infty} \int_0^\rho
\frac{\sin\sigma}{\sigma} \, d\sigma$ exists.

\end{proof}

The convergence results so far in this chapter rely just on
Riemann--Lebsgue and direct estimates. A much deeper result is:

\begin{theorem}[Carleson--Hunt]
If $f \in \lpt, 1< p < \infty$ then the Fourier series of $f$
converges to $f(t)$ for almost every $t \in \T$.
\end{theorem}

For $p=1$, the result is spectacularly false by Kolmogorov's
Theorem \ref{th:kol}.

\begin{proof}
Omitted. The idea is to prove that the Dirichlet maximal operator
$(D^* f)(t) = \sup_n |(D_n * f)(t)|$ is strong $(p,p)$ for
$1<p<\infty$. Then it is weak $(p,p)$, and so convergence a.e.\
follows from Chapter~\ref{ch:mf}.

Thus one wants
\[
\big\lv \sup_n |D_n * f| \big\rv_\lpt \leq C_p \lv f \rv_\lpt
\]
for $1<p<\infty$. The next Chapters show
\[
\sup_n \lv D_n * f \rv_\lpt \leq C_p \lv f \rv_\lpt ,
\]
but that is not good enough to prove Carleson--Hunt!
\end{proof}

\chapter{Fourier series: norm convergence} \label{ch:cn}

\subsubsection*{Goals}

Characterize norm convergence in terms of uniform norm bounds

\noindent Show norm divergence can occur for $\lot$ and $C(\T)$

\noindent Show norm convergence for $\lpt$ follows from boundedness
of the Hilbert transform

\subsubsection*{Reference} \cite{K} Section II.1

\vspace{18pt}

\begin{theorem} \label{th:nc1}
Let $B$ be one of the spaces $C(\T)$ or $\lpt, 1 \leq p < \infty$.

(a) If $\sup_n \lv S_n \rv_{B \to B} < \infty$ then Fourier series
converge in $B$:
\[
\lim_{n \to \infty} \lv S_n(f) - f \rv_B = 0 \qquad \text{for each
$f \in B$.}
\]

(b) \label{div} If $\sup_n \lv S_n \rv_{B \to B} = \infty$ then
there exists $f \in B$ whose Fourier series diverges unboundedly:
$\sup_n \lv S_n(f) \rv_B = \infty$.
\end{theorem}

\begin{proof} \

(b) This part follows immediately from the Uniform Boundedness
Principle in functional analysis.

(a) The collection of trigonometric polynomials is dense in $B$ (as remarked after Theorem~\ref{th:sk}). Further, if $g$ is a trigonometric polynomial then $S_n(g)=g$ whenever $n$ exceeds the degree of
$g$. Hence the set
\[
{\mathcal C} = \{ f \in B : \lim_{n \to \infty} S_n(f) = f \
\text{in $B$} \}
\]
is dense in $B$. The set ${\mathcal
C}$ is also closed, by the following proposition, and so ${\mathcal C}=B$, which proves part (a).
\end{proof}

\begin{proposition} \label{pr:red}
Let $B$ be any Banach space and assume the $T_n : B \to B$ are
bounded linear operators.

If $\sup_n \lv T_n \rv_{B \to B} < \infty$ then
\[
{\mathcal C} = \{f \in B : \lim_{n \to \infty} T_n f = f \ \text{in
$B$} \}
\]
is closed.
\end{proposition}

\begin{proof} Let $A=\sup_n \lv T_n \rv_{B \to B}$.
Consider a sequence $f_m \in {\mathcal C}$ with $f_m\to f$. We must
show $f \in {\mathcal C}$, so that ${\mathcal C}$ is closed.

Choose $\e > 0$ and fix $m$ such that $\lv f_m - f \rv < \e/2(A +
1)$. Since $f_m \in {\mathcal C}$ there exists $N$ such that $\lv
T_n f_m - f_m \rv < \e/2$ whenever $n > N$. Then
\begin{align*}
\lv T_n f - f \rv & \leq \lv T_n f - T_n f_m \rv + \lv T_n f_m - f_m
\rv + \lv f_m - f \rv \\
& \leq (A + 1) \lv f - f_m \rv + \lv T_n f_m - f_m \rv < \e
\end{align*}
whenever $n>N$, as desired.
\end{proof}

\subsubsection*{Norm Estimates}
\[
\lv S_n \rv_{B \to B} \leq \lv D_n \rv_{\lot}
\]
when $B$ is $C(\T)$ or $\lpt, 1 \leq p < \infty$, since
\begin{align*}
\lv S_n(f) \rv_B & = \big\lv D_n * f \big\rv_B \leq \lv D_n \rv_{\lot} \lv f \rv_B .
\end{align*}
This upper estimate is not useful, since we know $\lv D_n
\rv_{\lot} \to \infty$.

\begin{example}[Divergence in $C(\T)$]\rm
For $B = C(\T)$ we have
\[
\lv S_n \rv_{C(\T) \to C(\T)} = \lv D_n \rv_{\lot} .
\]
Indeed, for each $\e > 0$ one can construct $g \in C(\T)$ that
approximates $\sign(D_n)$ (like in Chapter~\ref{ch:cp}), so that
\[
\lv S_n(g) \rv_{C(\T)} \geq |S_n(g)(0)| \geq \big( \lv D_n
\rv_{\lot} - \e \big) \lv g \rv_{C(\T)}.
\]

Therefore $\sup_n \lv S_n \rv_{C(\T) \to C(\T)} = \infty$, so that
(by Theorem~\ref{th:nc1}(b)) there exists a continuous function $f
\in C(\T)$ whose Fourier series diverges unboundedly in the uniform
norm: $\sup_n \lv S_n(f) \rv_{C(\T)} = \infty$.

Of course, this result follows already from the pointwise divergence
in Chapter~\ref{ch:cp}.
\end{example}

\begin{example}[Divergence in $\lot$]\rm
For $B = \lot$ we have
\[
\lv S_n \rv_{\lot \to \lot} = \lv D_n \rv_{\lot} .
\]
\noindent \Proof Fix $n$. Then
$S_n(F_N) = F_N * D_n \to D_n$ in $\lot$ as $N \to \infty$, and so
\begin{align*}
\lv D_n \rv_\lot & = \lim_{N \to \infty} \lv S_n(F_N) \rv_\lot \\
& \leq \lv S_n \rv_{\lot \to \lot} \lv F_N \rv_\lot \\
& = \lv S_n \rv_{\lot \to \lot} .
\end{align*}

\vspace{6pt} Therefore $\sup_n \lv S_n \rv_{\lot \to \lot} =
\infty$, so that (by Theorem~\ref{th:nc1}(b)) there exists an
integrable function $f \in \lot$ whose Fourier series diverges
unboundedly in the $\lo$ norm: $\sup_n \lv S_n(f) \rv_{\lot} =
\infty$.

\vspace{6pt} \noindent \emph{Aside.} For an explicit example of
$L^1$ divergence, see \cite[Exercise 3.5.9]{G}.
\end{example}

\subsubsection*{Convergence in $\lpt, 1<p<\infty$}

\noindent 1. We shall prove (in Chapters~\ref{ch:h2}--\ref{ch:hp})
the existence of a bounded linear operator
\[
H: \lpt \to \lpt , \qquad 1<p<\infty,
\]
called the \textbf{Hilbert transform} on $\T$, with the property
\[
\widehat{(Hf)}(n) = -i \sign(n) \fhatn .
\]
(Thus $H$ is a \emph{Fourier multiplier} operator.) That is
\[
Hf \sim \sum_{n=-\infty}^\infty (-i)\sign(n) \fhatn e^{int}.
\]

\vspace{6pt} \noindent 2. Then the \emph{Riesz projection} $P: \lpt
\to \lpt$ defined by
\[
Pf = \frac{1}{2} \fhat(0) + \frac{1}{2}(f + iHf)
\]
is also bounded, when $1<p<\infty$. (Note the constant term
$\fhat(0)$ is bounded by $\lv f \rv_\lpt$, by H\"{o}lder's
inequality.)

Observe $P$ projects onto the nonnegative frequencies:
\[
Pf \sim \sum_{n \geq 0} \fhatn e^{int}
\]
since $i(-i \sign(n))=\sign(n)$.

\vspace{6pt} \noindent 3. The following formula expresses the
Fourier partial sum operator in terms of the Riesz projection and
some modulations:
\begin{equation} \label{eq:rp}
e^{-imt}P(e^{imt}f) - e^{i(m+1)t}P(e^{-i(m+1)t}f) = S_m(f).
\end{equation}
\noindent \Proof
\begin{align*}
e^{imt}f & \sim \sum_{n=-\infty}^\infty \fhatn e^{i(m+n)t} \\
P(e^{imt}f) & \sim \sum_{n \geq -m} \fhatn e^{i(m+n)t} \\
e^{-imt}P(e^{imt}f) & \sim \sum_{n \geq -m} \fhatn e^{int} \\
e^{i(m+1)t}P(e^{-i(m+1)t}f) & \sim \sum_{n \geq m+1} \fhatn e^{int}
\end{align*}
Subtracting the last two formulas gives $S_m(f)$, on the right side,
and we conclude that the left side of \eqref{eq:rp} has the same
Fourier coefficients as $S_m(f)$. By the uniqueness result
\eqref{eq:ut}, the left side of \eqref{eq:rp} must equal $S_m(f)$.

\vspace{6pt} \noindent 4. From \eqref{eq:rp} and boundedness of the
Riesz projection it follows that
\[
\sup_m \lv S_m \rv_{\lpt \to \lpt} \leq 2 \lv P \rv_{\lpt \to \lpt}
<  \infty
\]
when $1<p<\infty$. Hence from Theorem~\ref{th:nc1} we conclude:

\begin{theorem}[Fourier series converge in $\lpt$] Let $1<p<\infty$.
Then
\[
\lim_{n \to \infty} \lv S_n(f) - f \rv_{\lpt} = 0 \qquad \text{for
each $f \in \lpt$.}
\]
\end{theorem}

\vspace{6pt} It remains to prove $L^p$ boundedness of the Hilbert
transform.

\chapter{Hilbert transform on $\ltt$} \label{ch:h2}

\subsubsection*{Goal} Obtain time and frequency representations of
the Hilbert transform

\subsubsection*{Reference} \cite{EG} Section 6.3

\begin{definition}
The \emph{Hilbert transform} on $\ltt$ is
\begin{align*}
H : \ltt & \to \ltt \\
f & \mapsto \sum_{n= -\infty}^\infty \big (-i\sign(n) \fhatn \big)
e^{int} .
\end{align*}
\end{definition}
We call $\{-i\sign(n)\}$ the \emph{multiplier sequence} of $H$.

Since $|\sign(n)| \leq 1$, the definition indeed yields $Hf \in
\ltt$, with
\[
\lv Hf \rv_\ltt^2 = \sum_{n \in \Z} |\widehat{(Hf)}(n)|^2 = \sum_{n
\neq 0} |\fhatn|^2 \leq \lv \fhat \rv_{\ell^2(\Z)}^2 = \lv f
\rv_\ltt^2
\]
by Plancherel in Chapter~\ref{ch:l2}. Hence $\lv H \rv_{\lt \to \lt
} = 1$. Observe also $H^2(f)=H(Hf)=-\sum_{n \neq 0} \fhatn e^{int} =
-f+\fhat(0)$.

\begin{lemma}[Adjoint of Hilbert transform] \label{le:aht}
$H^*=-H$
\end{lemma}

\begin{proof}
For $f,g \in \ltt$,
\begin{align*}
\la Hf , g \ra_\ltt & = \la \widehat{Hf} , \ghat \ra_{\ell^2(\Z)} \\
& = \la -i \sign(n) \fhatn , \ghatn \ra_{\ell^2(\Z)} \\
& = \la \fhatn , i \sign(n) \ghatn \ra_{\ell^2(\Z)} \\
& = \la \fhat , -\widehat{Hg} \ra_{\ell^2(\Z)} \\
& = \la f , -Hg \ra_\ltt .
\end{align*}
\end{proof}

\begin{proposition} \label{pr:hst}
If $f \in \ltt$ is $C^1$-smooth on an open interval $I \subset \T$,
then
\begin{align}
(Hf)(t) & = \intp_0^\pi [ f(t-\tau) - f(t+\tau) ]
\cot \big( \frac{\tau}{2} \big) \, d\tau \label{eq:10repr} \\
& = \lim_{\e \to 0} \intp_{\e < |\tau| < \pi} f(t-\tau) \cot \big(
\frac{\tau}{2} \big) \, d\tau \label{eq:10conv}
\end{align}
for almost every $t \in I$.
\end{proposition}

\begin{remark}\rm
Formally \eqref{eq:10conv} says that
\[
Hf = f* \cot \big( \frac{t}{2} \big) .
\]
But the convolution is ill-defined because the Hilbert kernel
$\cot(t/2)$ is not integrable. That is why \eqref{eq:10conv}
evaluates the convolution in the \emph{principal valued} sense,
taking the limit of integrals over $\T \setminus [-\e,\e]$.
\end{remark}

{\it Proof.} First, geometric series calculations show that
\begin{align}
\sum_{n=-N}^{N} \big( -i\sign(n)e^{in\tau} \big) & =
i\sum_{n=-N}^{-1}e^{in\tau} - i\sum_{n=1}^{N}e^{in\tau} \notag\\
& = i\frac{e^{-i(N+1)\tau} - e^{-i\tau}}{e^{-i\tau}-1} -
i\frac{e^{i(N+1)\tau} - e^{i\tau}}{e^{i\tau} - 1} \notag \\
& = i\frac{e^{-i(N+1/2)\tau} - e^{-i\tau/2} + e^{i(N+1/2)\tau} -
  e^{i\tau/2}}{e^{-i\tau/2} - e^{i\tau/2}} \notag \\
& = \frac{\cos\left(\frac{\tau}{2}\right) - \cos\big(
    (N+\frac{1}{2}) \tau \big)}{\sin(\frac{\tau}{2})}.
    \label{eq:hst1}
\end{align}

Second, the $N$th partial sum of $Hf$ is

\begin{align*}
& \sum_{n=-N}^{N} \big (-i\sign(n)\fhatn \big) e^{int} \\
& = \intp_\T f(\tau) \sum_{n=-N}^{N} (-i)
\sign(n)e^{in(t-\tau)} \, d\tau \\
& = \intp_{-\pi}^\pi f(t-\tau) \sum_{n=-N}^{N}
  (-i) \sign(n)e^{in\tau} \, d\tau && \text{by $\tau \mapsto t-\tau$} \\
& = \intp_0^\pi [f(t-\tau) - f(t+\tau)]
\frac{\cos\left(\frac{\tau}{2}\right) -
  \cos\left((N+\frac{1}{2})\tau\right)}{\sin(\frac{\tau}{2})} \, d\tau \\
& = \intp_0^\pi [f(t-\tau) - f(t+\tau)] \cot \left(\frac{\tau}{2}\right) \, d\tau && \text{by \eqref{eq:hst1}} \\
& \qquad - \intp_0^\pi \frac{f(t-\tau) -
f(t+\tau)}{\sin(\frac{\tau}{2})} \cos \big( (N+\frac{1}{2} ) \tau
\big) \, d\tau.
\end{align*}

If $t \in I$ then the second integrand belongs to $\lot$ since it is
bounded for $\tau$ near $0$,  by the $C^1$-smoothness of $f$. Hence
the second integral tends to $0$ as $N \to \infty$ by the
Riemann-Lebesgue Corollary~\ref{co:rl}. Formula \eqref{eq:10repr}
now follows, because the partial sum
\[
\sum_{n=-N}^N \left( -i\sign(n) \fhatn \right) e^{in\tau}
\]
converges to $Hf(t)$ in $\ltt$ and hence some subsequence of the
partial sums converges to $(Hf)(t)$ a.e.

Now write \eqref{eq:10repr} as
\[
(Hf)(t) = \lim_{\e \to 0} \intp_\e^\pi [f(t-\tau) - f(t+\tau)] \cot
\left( \frac{\tau}{2} \right) \, d\tau
  \]
and use oddness of $\cot(\tau/2)$ to obtain \eqref{eq:10conv}.

\chapter{Calder\'{o}n--Zygmund decompositions} \label{ch:cz}

\subsubsection*{Goal} Decompose a function into good and bad parts,
preparing for a weak $(1,1)$ estimate on the Hilbert transform

\subsubsection*{References}

\cite{D} Section 2.5

\noindent \cite{G} Section 4.3

\vspace{18pt}

\begin{definition}
For $k \in \Z$, let
\[
Q_k = \{ 2^{-k} \big( [0,1)^d + m \big) : m \in \Zd \} .
\]
Notice the cubes in $Q_k$ are small when $k$ is large.

Call $\cup_k Q_k$ the \emph{collection of dyadic cubes}.
\end{definition}

\subsubsection*{Facts \text{(exercise)}}
\begin{enumerate}
\item For all $x \in \Rd$ and $k \in \Z$, there exists a unique $Q \in Q_k$ such that $x \in Q$. That is, there exists a unique $m \in \Zd$ with $x \in 2^{-k} \big( [0,1)^d + m \big)$.
\item Given $Q \in Q_k$ and $j<k$, there exists a unique $\widetilde{Q} \in Q_j$ with $Q \subset \widetilde{Q}$.
\item Each cube in $Q_k$ contains exactly $2^d$ cubes in $Q_{k+1}$.
\item Given two dyadic cubes, either one of them is contained in the other, or else the cubes are disjoint.
\end{enumerate}

\begin{definition}
For $f \in L^1_{loc}(\Rd)$, let
\[
(E_k f)(x) = \sum_{Q \in Q_k} \Big( \frac{1}{|Q|} \int_Q f(y) \, dy \Big) \charfn_Q(x) .
\]
Then $E_k f$ is constant on each cube in $Q_k$ (equalling there the average of $f$ over that cube), and
\begin{equation} \label{eq:cz1}
\int_\Omega E_k f \, dx = \int_\Omega f \, dx
\end{equation}
whenever $\Omega$ is a finite union of cubes in $Q_k$.

Define the \emph{dyadic maximal function}
\begin{align*}
(M_d f)(x) & = \sup_k |(E_k f)(x)| \\
& = \sup \big\{ \Big| \frac{1}{|Q|} \int_Q f(y) \, dy \Big| : \text{$Q$ is a dyadic cube containing $x$} \big\} .
\end{align*}
\end{definition}

\begin{theorem} \label{th:dm} \

(a) $M_d$ is weak $(1,1)$.

(b) If $f \in L^1_{loc}(\Rd)$ then $\lim_{k \to \infty} (E_k f)(x) = f(x)$ a.e.
\end{theorem}

\begin{proof}
We employ a ``stopping time'' argument like in probability theory for martingales.

For part (a), let $f \in \lord, \lambda>0$. Since $M_d f \leq M_d |f|$, we can assume $f \geq 0$. Let
\begin{align*}
\Omega & = \{ x \in \Rd : (M_d f)(x) > \lambda \} , \\
\Omega_k & = \{ x \in \Rd : \text{$(E_k f)(x) > \lambda$ and $(E_j f)(x) \leq \lambda$ for all $j < k$} \} .
\end{align*}
Clearly $\Omega_k \subset \Omega$. And if $x \in \Omega$ then $(E_k f)(x) > \lambda$ for some $k$; a smallest such $k$ exists, because
\begin{align*}
\lim_{j \to -\infty} (E_j f)(x) & \leq \lim_{j \to -\infty} \frac{1}{(2^{-j})^d} \int_\Rd f(y) \, dy \\
& = 0 \\
& < \lambda .
\end{align*}
Choosing the smallest $k$ implies $(E_j f)(x) \leq \lambda$ for all $j<k$, and so $x \in \Omega_k$. Hence $\Omega = \cup_k \Omega_k$, so that
\begin{align*}
|\Omega| & = \sum_k |\Omega_k| && \text{by disjointness of the $\Omega_k$} \\
& \leq \frac{1}{\lambda} \sum_k \int_{\Omega_k} E_k f \, dx && \text{since $E_k f > \lambda$ on $\Omega_k$} \\
& = \frac{1}{\lambda} \sum_k \int_{\Omega_k} f \, dx && \text{by \eqref{eq:cz1}, since $\Omega_k$ equals a union of cubes in $Q_k$} \\
& && \text{(recall $E_k f$ is constant on each cube in $Q_k$)} \\
& \leq \frac{1}{\lambda} \int_\Rd f \, dx .
\end{align*}
Therefore $M_d$ is weak $(1,1)$.

Part (b) holds if $f$ is continuous, and hence if $f \in L^1_{loc}(\Rd)$ by Theorem~\ref{th:mfcae} (exercise), using that the dyadic maximal operator $M_d$ is weak $(1,1)$.
\end{proof}

Note we did not need a covering lemma, when proving the dyadic maximal function is weak $(1,1)$, because disjointness of the cubes is built into the construction.

\begin{theorem}[Calder\'{o}n--Zygmund decomposition at level $\lambda$] \label{th:cz}
Let $f \in \lord, \lambda>0$. Then there exists a ``good' function $g \in \lo \cap \lird$ and a ``bad'' function $b \in \lord$ such that
\begin{enumerate}
\item[i.] $f=g+b$
\item[ii.] $\lv g \rv_\lord \leq \lv f \rv_\lord, \quad \lv g \rv_\lird \leq 2^d \lambda , \quad \lv b \rv_\lord \leq 2 \lv f \rv_\lord$
\item[iii.] $b = \sum_l b_l$ where $b_l$ is supported in a dyadic cube $Q(l)$ and the $\{ Q(l) \}$ are disjoint; we do \emph{not} assume $Q(l) \in Q_l$, just $Q(l) \in Q_k$ for some $k$.
\item[iv.] $\int_{Q(l)} b_l(x) \, dx = 0$
\item[v.] $\lv b_l \rv_\lord \leq 2^{d+1} \lambda |Q(l)|$
\item[vi.] $\sum_l |Q(l)| \leq \frac{1}{\lambda} \lv f \rv_\lord$
\end{enumerate}
\end{theorem}

\begin{proof}
Apply the proof of Theorem~\ref{th:dm} to $|f|$, and decompose the disjoint sets $\Omega_k$ into dyadic cubes in $Q_k$. Together, these cubes form the collection $\{ Q(l) \}$. Property (vi) is just the weak $(1,1)$ estimate that we proved.

For (i), (iii), (iv), argue as follows. Let
\[
b_l(x) = \Big( f(x) - \frac{1}{|Q(l)|} \int_{Q(l)} f(y) \, dy \Big) \charfn_{Q(l)}(x)
\]
so that $b_l$ integrates to $0$. Define
\[
b(x) = \sum_l b_l(x) =
\begin{cases}
f(x) - \frac{1}{|Q(l)|} \int_{Q(l)} f(y) \, dy & \text{on $Q(l)$, for each $l$,} \\
0 & \text{on $\Rd \setminus \cup_l Q(l)$.}
\end{cases}
\]
Then let
\begin{align*}
g & = f-b \\
& =
\begin{cases}
\frac{1}{|Q(l)|} \int_{Q(l)} f(y) \, dy & \text{on $Q(l)$, for each $l$,} \\
f(x) & \text{on $\Rd \setminus \cup_l Q(l)$.}
\end{cases}
\end{align*}

For (ii), note $\lv g \rv_\lord \leq \lv f \rv_\lord$, since $g=f$ off $\cup_l Q(l)$ and on $Q(l)$ we have
\[
\int_{Q(l)} |g(x)| \, dx \leq \int_{Q(l)} |f(x)| \, dx .
\]
Hence $\lv b \rv_\lord = \lv f-g \rv_\lord \leq 2 \lv f \rv_\lord$.

Next we show $\lv g \rv_\lird \leq 2^d \lambda$. Suppose $x \in \Rd \setminus \cup_l Q(l)$. Then $g(x)=f(x)$. Since $x \notin \Omega_k$ for all $k$ we have $(E_k |f|)(x) \leq \lambda$ for all $k$. Hence $|f(x)| \leq \lambda$ (for almost every such $x$) by Theorem~\ref{th:dm}(b), so that $|g(x)| \leq \lambda$.

Next suppose $x \in Q(l)$ for some $l$, so that $x \in \Omega_k$ for some $k$. Then $(E_{k-1}|f|)(x) \leq \lambda$, which means
\[
\frac{1}{|Q|} \int_Q |f(y)| \, dy \leq \lambda
\]
for some cube $Q \in Q_{k-1}$ with $x \in Q(l) \subset Q$. Hence
\begin{equation} \label{eq:cz2}
\frac{1}{2^d |Q(l)|} \int_{Q(l)} |f(y)| \, dy \leq \lambda
\end{equation}
since $Q(l) \subset Q$ and $\side(Q)=2\side(Q(l))$. Therefore $|g(x)| \leq 2^d \lambda$, by definition of $g$.

For (v), just note
\begin{align*}
\int_{Q(l)} |b_l(x)| \, dx & \leq 2 \int_{Q(l)} |f(x)| \, dx \qquad \text{by definition of $b_l$} \\
& \leq 2^{d+1} \lambda |Q(l)|
\end{align*}
by \eqref{eq:cz2}.
\end{proof}

Now we adapt the theorem to $\T$. We will restrict to ``large'' $\lambda$ values, so that the dyadic intervals have length at most $2\pi$ and thus fit into $\T$.

\begin{corollary}[Calder\'{o}n--Zygmund decomposition on $\T$] \label{co:cz}
Let $f \in \lot, \lambda > \lv f \rv_\lot$. Then there exists a ``good' function $g \in \lit$ and a ``bad'' function $b \in \lot$ such that
\begin{enumerate}
\item[i.] $f=g+b$
\item[ii.] $\lv g \rv_\lot \leq \lv f \rv_\lot, \quad \lv g \rv_\lit \leq 2 \lambda , \quad \lv b \rv_\lot \leq 2 \lv f \rv_\lot$
\item[iii.] $b = \sum_l b_l$ where $b_l$ is supported in some interval $I(l)$ of the form $2\pi \cdot 2^{-k} \big( [0,1) + m \big)$ where $k \geq 1, 0 \leq m \leq 2^k - 1$, and where the $\{ I(l) \}$ are disjoint.
\item[iv.] $\int_{I(l)} b_l(t) \, dt = 0$
\item[v.] $\lv b_l \rv_\lot \leq \frac{4}{2\pi} \lambda |I(l)|$
\item[vi.] $\sum_l |I(l)| \leq \frac{2\pi}{\lambda} \lv f \rv_\lot$
\end{enumerate}
\end{corollary}

\begin{proof}
Let $d=1$. Apply the Calder\'{o}n--Zygmund Theorem~\ref{th:cz} to
\[
\widetilde{f}(t) =
\begin{cases}
f(2\pi t) , & 0 \leq t < 1 , \\
0 , & \text{otherwise,}
\end{cases}
\]
to get $\widetilde{f}=\widetilde{g}+\widetilde{b}$. Note $\Omega_k$ is empty for $k \leq 0$, since
\begin{align*}
(E_k|\widetilde{f}|)(t) & \leq \frac{1}{2^{-k}} \int_0^1 |\widetilde{f}(\tau)| \, d\tau \\
& = 2^k \pr \int_0^{2\pi} |f(\tau)| \, d\tau \\
& \leq \lv f \rv_\lot \qquad \text{since $k \leq 0$} \\
& < \lambda
\end{align*}
by assumption on $\lambda$.

Further, $\Omega_k \subset [0,1]$ for $k \geq 1$, since $E_k|\widetilde{f}|=0$ outside $[0,1]$. Thus $I(l)=2\pi Q(l)$ has the form stated in the Corollary.

The Corollary now follows from Theorem~\ref{th:cz}, with $\widetilde{f}=\widetilde{g}+\widetilde{b}$ yielding $f=g+b$.
\end{proof}

\chapter{Hilbert transform on $\lpt$} \label{ch:hp}

\subsubsection*{Goals}

Prove a weak $(1,1)$ estimate on the Hilbert transform on $\T$

\noindent Deduce strong $(p,p)$ estimates by interpolation and
duality

\subsubsection*{Reference} \cite{D} Section 3.3

\vspace{18pt}

\begin{theorem}[weak $(1,1)$ on $\ltt$] \label{th:hw11}
There exists $A>0$ such that
\[
| \{t \in \T : |(Hf)(t)| > \lambda \} | \leq \frac{A}{\lambda} \lv
f\rv_\lot
\]
for all $\lambda>0$ and $f \in \ltt$.
\end{theorem}

\begin{proof}
If $\lambda \leq \lv f \rv_\lot$ then $A=2\pi$ works. So suppose
$\lambda > \lv f \rv_\lot$. Apply the Calder\'{o}n--Zygmund
Corollary~\ref{co:cz} to get $f=g+b$. Note $g \in \lit$ and so $g
\in \ltt$, hence $Hg \in \ltt$ by Chapter~\ref{ch:h2}. And $b=f-g
\in \ltt$ so that $Hb \in \ltt$. Further, $b_l \in \ltt$ and
$b=\sum_l b_l$ with convergence in $\ltt$, using disjointness of the
supports of the $b_l$. Hence $Hb=\sum_l Hb_l$ with convergence in
$\ltt$.

Since $Hf=Hg+Hb$, we have
\begin{align*}
& | \{ t \in \T : |(Hf)(t)| > \lambda \}| \\
& \leq | \{ t \in \T :
|(Hg)(t)| > \lambda/2 \} | + | \{ t \in \T : |(Hb)(t)| > \lambda/2 \} | \\
& = \gamma+\beta,
\end{align*}
say. First, use the $\lt$ theory on $g$:
\begin{align*}
\gamma & \leq \int_\T \frac{|(Hg)(t)|^2}{(\lambda/2)^2} \, dt \\
& \leq \frac{4}{\lambda^2} \int_\T |g(t)|^2 \, dt
&& \text{since $\lv H \rv_{\ltt \to \ltt} = 1$ by Chapter~\ref{ch:h2}} \\
& \leq \frac{8}{\lambda} \int_\T |g(t)| \, dt && \text{since $\lv g \rv_\lit \leq 2\lambda$} \\
& \leq \frac{8 \cdot 2\pi}{\lambda} \lv f \rv_\lot && \text{since
$\lv g \rv_\lot \leq \lv f \rv_\lot$.}
\end{align*}
Second, use $\lo$ estimates on $b$, as follows:
\begin{align*}
\beta & \leq \big| \bigcup_l 2I(l) \big| + | \{ t \in \T \setminus \bigcup_l 2I(l) : |(Hb)(t)| > \lambda/2 \}| \\
& \leq \frac{4\pi}{\lambda} \lv f \rv_\lot + \int_{\T \setminus
\cup_l 2I(l)} \frac{|(Hb)(t)|}{\lambda/2} \, dt
\\
& \qquad \qquad \text{by the Calder\'{o}n--Zygmund Corollary~\ref{co:cz}(vi)} \\
& \leq \frac{4\pi}{\lambda} \lv f \rv_\lot + \frac{2}{\lambda}
\sum_l \int_{\T \setminus 2I(l)} |(Hb_l)(t)| \, dt
\end{align*}
since $|Hb| \leq \sum_l |Hb_l|$ a.e.

To finish the proof, we show
\begin{equation} \label{eq:hlp1}
\sum_l \int_{\T \setminus 2I(l)} |(Hb_l)(t)| \, dt \leq
(\text{const.}) \lv f \rv_\lot .
\end{equation}
By Proposition~\ref{pr:hst} on the interval $\T \setminus 2I(l)$, we
have
\begin{align*}
& \int_{\T \setminus 2I(l)} |Hb_l(t)| \, dt \\
& = \int_{\T \setminus 2I(l)} \Big| \intp_{I(l)} b_l(\tau) \cot \big( \frac{1}{2}(t-\tau) \big) \, d\tau \Big| \, dt \\
& \qquad \qquad \text{noting $t-\tau$ is bounded away from $0$,
since $\tau \in I(l)$ and $t \notin 2I(l)$,} \\
& = \int_{\T \setminus 2I(l)} \Big| \intp_{I(l)} b_l(\tau) \big[ \cot \big( \frac{1}{2}(t-\tau) \big)-\cot \big( \frac{1}{2}(t-c_l) \big) \big] \, d\tau \Big| \, dt \\
& \qquad \qquad \text{where $c_l$ is the center of $I(l)$, using here that $\int_{I(l)} b_l(\tau) \, d\tau=0$,} \\
& = \int_{\T \setminus 2I(l)} \Big| \intp_{I(l)} b_l(\tau) \frac{\sin \big(\frac{1}{2}(\tau-c_l) \big)}{\sin \big(\frac{1}{2}(t-\tau) \big) \sin \big(\frac{1}{2}(t-c_l) \big)} \, d\tau \Big| \, dt\\
& \leq (\text{const.}) \int_{I(l)} |b_l(\tau)| \int_{\R \setminus
2I(l)} \frac{|I(l)|}{|t-\tau| |t-c_l|} \, dt d\tau .
\end{align*}
Note that
\begin{align*}
|t-c_l| & \leq |t-\tau| + |\tau-c_l| \\
& \leq |t-\tau| + \frac{1}{2}|I(l)| && \text{when $\tau \in I(l)$} \\
& \leq 2|t-\tau| && \text{when $t \in \R \setminus 2I(l)$.}
\end{align*}
Hence
\begin{align*}
\int_{\R \setminus 2I(l)} \frac{|I(l)|}{|t-\tau| |t-c_l|} \, dt & \leq 2 \int_{\R \setminus 2I(l)} \frac{|I(l)|}{|t-c_l|^2 } \, dt \\
& = 4 \int_{2r}^{\infty} \frac{2r}{t^2} \, dt && \text{where $2r=|I(l)|$} \\
& = 4.
\end{align*}
Thus
\begin{align*}
\text{the left side of \eqref{eq:hlp1}}
& \leq (\text{const.}) \sum_l \int_{I(l)} |b_l(\tau)| \, d\tau \\
& = (\text{const.}) \lv b \rv_\lot \\
& \leq (\text{const.}) \lv f \rv_\lot
\end{align*}
by the Calder\'{o}n--Zygmund Corollary~\ref{co:cz}.

We have proved \eqref{eq:hlp1}, and thus the theorem.
\end{proof}

\begin{corollary} \label{co:hspp}
The Hilbert transform is strong $(p,p)$ for $1<p<\infty$, with
$\widehat{(Hf)}(n)=-i \sign(n) \fhatn$ for all $f \in \lpt,
 n \in \Z$.
\end{corollary}

\begin{proof}
$H$ is strong $(2,2)$ and linear, by definition in
Chapter~\ref{ch:h2}, and $H$ is weak $(1,1)$ on $\ltt$ (and hence on
the simple functions on $\T$) by Theorem~\ref{th:hw11}. So $H$ is
strong $(p,p)$ for $1<p<2$ by Remark~\ref{re:mi} after Marcinkiewicz
Interpolation (in Appendix~\ref{ap:ip}). That is, $H : \lpt \to
\lpt$ is bounded and linear for $1<p<2$.

For $2<p<\infty$ we will use duality and anti-selfadjointness
$H^*=-H$ on $\ltt$ (see Lemma~\ref{le:aht}) to reduce to the case
$1<p<2$. Suppose $1<p^\prime<2<p<\infty$ with $\frac{1}{p}+\frac{1}{p^\prime}=1$, and note 
\[
L^p(\T) \subset L^2(\T) \subset L^{p^\prime}(\T) .
\]
If $f \in L^p(\T)$ then
\begin{align*}
\lv Hf \rv_p & = \sup \big\{ \big| \intp_\T (Hf) \overline{g} \, dt \big| : g \in L^{p^\prime}(\T) \text{\ with norm $1$} \big\} \\
& = \sup\{ \big| \intp_\T (Hf) \overline{g} \, dt \big| : g \in L^2(\T) \text{\ with $L^{p^\prime}$-norm $1$} \} \\
& \qquad \qquad \qquad \text{by density of $L^2$ in $L^{p^\prime}$} \\
& = \sup \{ \big| \intp_\T f (\overline{Hg}) \, dt \big| : g \in L^2(\T) \text{\ with $L^{p^\prime}$-norm $1$} \}\\
& \qquad \qquad \qquad \text{since $H^*=-H$ on $\ltt$} \\
& \leq \lv f \rv_\lpt \sup \{ \lv Hg \rv_{L^{p^\prime}(\T)} : g \in L^{p^\prime}(\T) \text{\ with norm $1$} \} \quad \text{by H\"{o}lder} \\
& \leq \text{(const.)}_{p^\prime} \lv f \rv_\lpt
\end{align*}
by the strong $(p^\prime,p^\prime)$ bound proved above, using that
$1<p^\prime<2$. Thus $H$ is a bounded
operator on $\lpt$.

Finally, for $f \in \lpt, 1<p<\infty$, let $f_m \in L^p \cap
L^2(\T)$ with $f_m \to f$ in $\lpt$. Boundedness of $H$ on $L^p$
implies $Hf_m \to Hf$ in $L^p$. Hence $f_m \to f$ and $Hf_m \to Hf$
in $\lot$, and so passing to the limit in $\widehat{(Hf_m)}(n)=-i
\sign(n) \widehat{f_m}(n)$ yields $\widehat{(Hf)}(n)=-i \sign(n)
\fhatn$, as desired.
\end{proof}

\chapter{Applications of interpolation} \label{ch:ai}

\subsubsection*{Goal} Apply Marcinkiewicz and Riesz--Thorin
interpolation to the Hilbert transform, maximal operator, Fourier
analysis and convolution

\vspace{18pt}

The Marcinkiewicz and Riesz--Thorin interpolation theorems are
covered in Appendix~\ref{ap:ip}. Some important applications are:

\paragraph*{Hilbert transform.}
\[
\text{$H : \lpt \to \lpt$ is bounded, for $1<p<\infty$,}
\]
by the Marcinkiewicz interpolation and duality argument in
Corollary~\ref{co:hspp}.

\paragraph*{Hardy--Littlewood maximal operator.}

$M$ is weak $(1,1)$ and strong $(\infty,\infty)$ by
Chapter~\ref{ch:mf}, and hence $M$ is strong $(p,p)$ for
$1<p<\infty$ by the Marcinkiewicz Interpolation Theorem~\ref{th:mi}.
(Note $M$ is sublinear.)

Strong $(p,p)$ was proved directly, already, in Chapter~\ref{ch:mf}.

\paragraph*{Fourier analysis.} The \emph{Hausdorff--Young} theorem says
\[
\text{$\widehat{\ } : \lpt \to \ell^{p^\prime}(\Z)$}, \qquad 1 \leq
p \leq 2, \quad \frac{1}{p} + \frac{1}{p^\prime} = 1 .
\]
It fails for $p>2$ \cite[Section~IV.2.3]{K}.

To interpret the theorem, note $\lpt$ gets smaller as $p$ increases,
and so does $\ell^{p^\prime}(\Z)$.

\noindent \Proof The analysis operators $\widehat{\ } : \lot \to
\ell^\infty(\Z)$ and $\widehat{\ } : \ltt \to \ell^2(\Z)$ are
bounded. Observe
\[
\frac{1}{p} = \frac{1-\theta}{1}+\frac{\theta}{2} \quad
\Longleftrightarrow \quad \frac{\theta}{2} = 1 - \frac{1}{p} \quad
\Longleftrightarrow \quad \frac{1}{p^\prime} =
\frac{1-\theta}{\infty} + \frac{\theta}{2} .
\]
Now apply the Riesz--Thorin Interpolation Theorem~\ref{th:rti}.

\paragraph*{Convolution.} The \emph{Generalized Young's} theorem says
\[
\lv f*g \rv_{L^r(\Rd)} \leq \lv f \rv_\lprd \lv g \rv_{L^q(\Rd)}
\qquad \text{when} \quad \frac{1}{p}+\frac{1}{q}=\frac{1}{r}+1 ,
\quad 1 \leq p,q,r \leq \infty .
\]
\noindent \Proof Fix $g \in L^q(\Rd)$ and define $Tf=f*g$. Then $T$
is strong $(1,q)$ since
\[
\lv f*g \rv_{L^q(\Rd)} \leq \lv f \rv_\lord \lv g \rv_{L^q(\Rd)}
\]
by Young's Theorem~\ref{th:yt}, and $T$ is strong $(q^\prime,\infty)$
since
\[
\lv f*g \rv_{\lird} \leq \lv f \rv_{L^{q^\prime}(\Rd)} \lv g
\rv_{L^q(\Rd)}
\]
by H\"{o}lder's inequality. In both cases, $\lv T \rv \leq \lv g
\rv_{L^q(\Rd)}$. Observe
\[
\frac{1}{p} = \frac{1-\theta}{1}+\frac{\theta}{q^\prime} \quad
\Longleftrightarrow \quad \frac{\theta}{q} = 1 - \frac{1}{p} =
\frac{1}{q}-\frac{1}{r} \quad \Longleftrightarrow \quad \frac{1}{r}
= \frac{1-\theta}{q} + \frac{\theta}{\infty} .
\]
Now apply the Riesz--Thorin Interpolation Theorem~\ref{th:rti}.

\chapter*{Epilogue: Fourier series in higher dimensions}

We have studied Fourier series only on the one dimensional torus $\T
= \R / 2\pi \Z$. The theory extends readily to the higher
dimensional torus $\Td = \Rd / 2\pi \Zd$.

Summability kernels can be obtained by taking products of one
dimensional kernels. Thus the higher dimensional Dirichlet kernel is
\begin{align*}
{\mathcal D}_n(t)
& = D_n(t_1) \cdots D_n(t_d) \\
& = \sum_{j_1,\ldots,j_d = -n}^n e^{ijt} ,
\end{align*}
where $j=(j_1,\ldots,j_d), t = (t_1,\ldots,t_d)^\dagger$ and
$\dagger$ denotes the transpose operation.

The Dirichlet kernel corresponds to ``cubical'' partial sums of
multiple Fourier series, because
\begin{align*}
({\mathcal D}_n * f)(t) & = \frac{1}{(2\pi)^d} \int_\T \cdots
\int_\T {\mathcal D}_n(t-\tau) f(\tau) \, d\tau_1 \cdots d\tau_d \\
& = \sum_{j_1,\ldots,j_d = -n}^n \fhatj e^{ijt} .
\end{align*}
``Spherical'' partial sums of the form $\sum_{|j| \leq n} \fhatj
e^{ijt}$ can be badly behaved. For example, they can fail to
converge for $f \in \lptd$ when $p \neq 2$. See \cite{G} for this
theorem and more on Fourier series in higher dimensions.

\part{Fourier integrals}

\chapter*{Prologue: Fourier series converge to Fourier integrals}

Fourier series do not apply to a function $g \in \lor$, since $g$ is
not periodic. Instead we take a large \emph{piece} of $g$ and look
at its Fourier series: for $\rho > 0$, let
\[
f(t) = g(\rho t) , \qquad t \in [-\pi,\pi) ,
\]
and extend $f$ to be $2\pi$-periodic. Then
\begin{align*}
\fhatj & = \intp_{-\pi}^\pi g(\rho t) e^{-ijt} \, dt \\
& = \frac{1}{2\pi \rho} \int_{-\rho \pi}^{\rho \pi} g(y)
e^{-i(j/\rho)y} \, dy
\end{align*}
by changing variable. Formally, for $|x| < \rho \pi$ we have
\begin{align*}
g(x) = f(\rho^{-1} x) & = \sum_{j=-\infty}^\infty \fhatj
e^{ij(\rho^{-1} x)} \\
& = \pr \sum_{j=-\infty}^\infty \big( \int_{-\rho \pi}^{\rho \pi}
g(y) e^{-i(j/\rho)y} \, dy \big) e^{i(j/\rho)x} \cdot \frac{1}{\rho} \\
& \to \pr \int_{-\infty}^\infty \big( \int_{-\infty}^\infty g(y)
e^{-i\xi y} \, dy \big) e^{i\xi x} \, d\xi
\end{align*}
as $\rho \to \infty$, by using Riemann sums on the $\xi$-integral.

The inner integral (``Fourier transform'') is analogous to a Fourier
coefficient.

The outer integral (``Fourier inverse'') is analogous to a Fourier
series.

We aim to develop a Fourier integral theory that is analogous to the
theory of Fourier series.

\chapter{Fourier transforms: basic properties} \label{ch:ft}

\subsubsection*{Goal}

Derive basic properties of Fourier transforms

\subsubsection*{Reference} \cite{K} Section VI.1

\subsubsection*{Notation}

\noindent $\lv f \rv_\lprd = \big( \int_\Rd |f(x)|^p \, dx
\big)^{1/p}$

\noindent Nesting of $L^p$-spaces fails: $\lird \not\subset \ltrd
\not\subset \lord$ due to behavior at infinity \emph{e.g.}\
$1/(1+|x|)$ is in $\ltr$ but not $\lor$

\noindent $C_c(\Rd) = \{ \text{complex-valued, continuous functions
with compact support} \}$

\noindent $C_0(\Rd) = \{ \text{complex-valued, continuous functions
with $f(x) \to 0$ as $|x| \to \infty$} \}$, Banach space with norm
$\lv \cdot \rv_\lird$

\noindent Translation $f_y(x) = f(x-y)$

\begin{definition}
For $f \in \lord$ and $\xi \in \Rd$, define
\begin{align}
\fhatxi & = \text{\emph{Fourier transform} of $f$} \notag \\
& = \int_\Rd f(x) e^{-i\xi x} \, dx . \label{eq:fhatdefr}
\end{align}
Here $\xi$ is a row vector, $x$ is a column vector, and so $\xi x =
\xi_1 x_1 + \cdots + \xi_d x_d$ equals the dot product.
\end{definition}

\begin{theorem}[Basic properties] \label{th:bpr} Let $f,g \in \lord, \xi, \omega \in \Rd, c
\in \C, y \in \Rd, A \in GL(\R,d)$.

\noindent Linearity $\widehat{(f+g)}(\xi)=\fhatxi+\ghatxi$ and
$\widehat{(cf)}(\xi)=c\fhatxi$

\noindent Conjugation
$\widehat{\overline{f}}(\xi)=\overline{\fhat(-\xi)}$

\noindent $\widehat{\ }$ \ takes translation to modulation,
$\widehat{\, f_y \,}(\xi) = e^{-i\xi y} \fhatxi$

\noindent $\widehat{\ }$ \ takes modulation to translation,
$[f(x)e^{i\omega x}]\widehat{\ }(\xi) = \fhat(\xi - \omega)$

\noindent $\widehat{\ }$ \ takes matrix dilation to its inverse,
$[\, |\det A| f(Ax) \, ]\widehat{\ }(\xi) = \fhat(\xi A^{-1})$

\noindent $\widehat{\ } : \lord \to \lird$ is bounded, with $\lv
\fhat \rv_\lird \leq \lv f \rv_\lord$

\noindent $\fhat$ is uniformly continuous

\noindent If $f_m \to f$ in $\lord$ then $\widehat{f_m} \to \fhat$
in $\lird$.
\end{theorem}
\begin{proof}
Exercise. For continuity, observe
\begin{align*}
|\fhat(\xi+\omega)-\fhat(\xi)| & \leq \int_\Rd |f(x)| |e^{-i\xi x}|
|e^{-i\omega x}-1| \, dx \\
& \to 0
\end{align*}
as $\omega \to 0$, by dominated convergence. The convergence is
independent of $\xi$, and so $\fhat$ is uniformly continuous.
\end{proof}

\begin{corollary}[Transform of a radial function] \label{co:rf}
If $f \in \lord$ is radial then $\fhat$ is radial.
\end{corollary}
Recall that $f$ is \emph{radial} if it depends only on the distance
to the origin: $f(x)=F(|x|)$ for some function $F$. Equivalently,
$f$ is radial if $f(Ax)=f(x)$ for every $x$ and every orthogonal
(``rotation and reflection'') matrix $A$.
\begin{proof}
Suppose $A$ is orthogonal. Then $f(Ax)=f(x)$ (since $f$ is radial)
and so
\[
\fhat(\xi A^{-1}) = [\, |\det A| f(Ax) \, ]\widehat{\ }(\xi) =
\fhatxi ,
\]
using Theorem~\ref{th:bpr} and that $|\det A|=1$.
\end{proof}

\begin{lemma}[Transform of a product] \label{le:tp}
If $f_1,\ldots,f_d \in \lor$ then $f(x)=\prod_{j=1}^d f_j(x_j)$ has
transform $\fhatxi=\prod_{j=1}^d \widehat{\, f_j \,}(\xi_j)$.
\end{lemma}
\begin{proof}
Use Fubini and the homomorphism property of the exponential:
$e^{-i\xi x} = \prod_{j=1}^d e^{-i \xi_j x_j}$.
\end{proof}

\begin{lemma}[Difference formula] \label{le:diffR} For $\xi \neq 0$,
\[
\fhatxi = \frac{1}{2} \int_\Rd [f(x)-f(x-\pi \xi^\dagger/|\xi|^2)]
\, e^{-i\xi x} \, dx ,
\]
where $\xi^\dagger$ is the column vector transpose of $\xi$.
\end{lemma}
\begin{proof}
Like Lemma~\ref{le:diffT}.
\end{proof}
\begin{lemma}[Continuity of translation] \label{le:ctR} Fix $f \in \lprd, 1 \leq p < \infty$. The map
\begin{align*}
\phi : \Rd & \to \lprd \\
y & \mapsto f_y
\end{align*}
is continuous.
\end{lemma}
\begin{proof}
Like Lemma~\ref{le:ct} except using $C_c(\Rd)$, which is dense in
$\lprd$.
\end{proof}
\begin{corollary}[Riemann--Lebesgue lemma] \label{co:rlr} $\fhatxi \to 0$ as $|\xi| \to
\infty$. Thus $\fhat \in C_0(\Rd)$.
\end{corollary}
\begin{proof}
Lemma~\ref{le:diffR} implies
\[
|\fhatxi| \leq \frac{1}{2} \lv f - f_{\pi \xi^\dagger/|\xi|^2}
\rv_\lord ,
\]
which tends to zero as $|\xi| \to \infty$ by the $\lo$-continuity of
translation in Lemma~\ref{le:ctR}, since $\xi^\dagger/|\xi|^2$ has
magnitude $1/|\xi| \to 0$.
\end{proof}

\begin{example} \label{ex:ftexamples}
We compute the Fourier transforms in Table~\ref{ta:ftexamples}.

\vspace{6pt} 1. $\int_\R \charfn_{[-1,1]}(x) e^{-i\xi x} \, dx =
\int_{-1}^1 e^{-i\xi x} \, dx = 2\sin (\xi) / \xi$

\vspace{6pt} 2.$\int_\R (1-|x|) \charfn_{[-1,1]}(x) e^{-i\xi x} \,
dx = 2 \int_0^1 (1-x) \cos(\xi x) \, dx = 2 \xi^{-2} (1-\cos \xi)$,
and $1-\cos \xi=2\sin^2(\xi/2)$

\vspace{6pt} 4. Next we compute for the fourth example, the Gaussian $e^{-|x|^2/2}$, so that we can use it later for the third example $e^{-|x|}$.

For $d=1$, let $g(\xi)=\int_\R e^{-x^2/2} e^{-i\xi
x} \, dx$ be the transform we want. Note $g(0)=\sqrt{2\pi}$.
Differentiating,
\[
g^\prime(\xi) = \int_\R e^{-x^2/2} (-ix) e^{-i\xi x} \, dx ,
\]
with the differentiation through the integral justified by using
difference quotients and dominated convergence (Exercise). Hence
\begin{align*}
g^\prime(\xi) & = i \int_\R \big( e^{-x^2/2} \big)^\prime e^{-i\xi x} \, dx \\
& = -i \int_\R e^{-x^2/2} \big( e^{-i\xi x} \big)^\prime \, dx && \text{by parts} \\
& = -\xi \int_\R e^{-x^2/2} e^{-i\xi x} \, dx \\
& = -\xi g(\xi) .
\end{align*}
Solving the differential equation yields
$g(\xi)=\sqrt{2\pi}e^{-\xi^2/2}$.

For $d>1$, note the product structure $e^{-|x|^2/2} = \prod_{j=1}^d
e^{-x_j^2/2}$ and apply Lemma~\ref{le:tp}.

\begin{table}[t]
\begin{center}
\begin{tabular}{|c|c|c|}
\hline & & \\
dimension    & $f(x)$    & $\fhatxi$
  \\
\hline   & & \\ $1$ & $\charfn_{[-1,1]}(x)$ & $2 \frac{\sin
\xi}{\xi} = 2 \sinc
\xi$ \\
& & \\ $1$ & $(1-|x|) \charfn_{[-1,1]}(x)$ & $\left( \frac{\sin (
\xi/2)}{\xi/2} \right)^{\! \! 2} =\sinc^2(\xi/2)$ \\
& & \\
$d$ & $e^{-|x|}$ & $\frac{(2\pi)^d c_d}{(1+|\xi|^2)^{(d+1)/2}}$ \\
& & \\ $d$ & $e^{-|x|^2/2}$ & $(2\pi)^{d/2} e^{-|\xi|^2/2}$ \\
\hline
\end{tabular}
\vspace*{3pt} \caption{Fourier transforms from
Example~\ref{ex:ftexamples}. In the third example, $ c_d = \left. \!
\Gamma \left( \frac{d+1}{2} \right) \right/ \pi^{(d+1)/2}$,  so that
$c_1 = 1/\pi$. The fourth example says the Fourier transform of a
Gaussian is a Gaussian.} \label{ta:ftexamples}
\end{center}
\end{table}
\end{example}

\vspace{6pt} 3. For $d=1$, $\int_\R e^{-|x|} e^{-i\xi x} \, dx =
\int_0^\infty e^{-(1+i\xi)x} \, dx + \int_{-\infty}^0 e^{(1-i\xi)x}
\, dx =1/(1+i\xi)+1/(1-i\xi)$, which simplifies to the desired
result.

To handle $d>1$, we need a calculus lemma that expresses a decaying exponential as a superposition of Gaussians.

\begin{lemma} \label{le:exp}
For $b>0$,
\[
e^{-b} = \frac{1}{\sqrt{2\pi}} \int_0^\infty \frac{e^{-a/2}}{\sqrt{a}} e^{-b^2/2a} \, da .
\]
\end{lemma}

\begin{proof}
\begin{align*}
& \ e^b \frac{1}{\sqrt{2\pi}} \int_0^\infty \frac{e^{-a/2}}{\sqrt{a}} e^{-b^2/2a} \, da \\
& =
\frac{2\sqrt{b}}{\sqrt{2\pi}} \int_0^\infty e^{-b(c-1/c)^2/2} \, dc && \text{by letting $a=bc^2$} \\
& =
\frac{2\sqrt{b}}{\sqrt{2\pi}} \int_0^\infty e^{-b(c-1/c)^2/2} c^{-2} \, dc && \text{by $c \mapsto 1/c$} \\
& =
\frac{\sqrt{b}}{\sqrt{2\pi}} \int_0^\infty e^{-b(c-1/c)^2/2} (1+c^{-2}) \, dc && \text{by averaging the last two formulas} \\
& =
\frac{\sqrt{b}}{\sqrt{2\pi}} \int_{-\infty}^\infty e^{-bu^2/2} \, du && \text{where $u=c-1/c$} \\
& = 1 .
\end{align*}
\end{proof}

Now we compute the Fourier transform of $e^{-|x|}$ as
\begin{align*}
& \ \int_\Rd e^{-|x|} e^{-i\xi x} \, dx \\
& = \frac{1}{\sqrt{2\pi}} \int_0^\infty \frac{e^{-a/2}}{\sqrt{a}} \int_\Rd e^{-|x|^2/2} e^{-i(\xi \sqrt{a}) x} \, dx \, a^{d/2} \, da \qquad \text{by Lemma~\ref{le:exp} and $x \mapsto \sqrt{a} x$} \\
& = \frac{1}{\sqrt{2\pi}} \int_0^\infty a^{(d-1)/2} e^{-a/2} (2\pi)^{d/2} e^{-|\xi \sqrt{a}|^2/2} \, da
\qquad \text{by the Gaussian in Table~\ref{ta:ftexamples}} \\
& = (2\pi)^{(d-1)/2} \big( (1+|\xi|^2)/2 \big)^{-(d+1)/2} \int_0^\infty u^{(d-1)/2} e^{-u} \, du
\end{align*}
where $u=a(1+|\xi|^2)/2$. The last integral is $\Gamma((d+1)/2)$, so that the transform equals $(2\pi)^d c_d (1+|\xi|^2)^{-(d+1)/2}$ as claimed in the Table.

\subsubsection*{Smoothness and decay}

\begin{theorem}[Differentiation and Fourier transforms] \label{th:dft} \

(a) If $f \in C^1_c(\Rd)$ (or more generally, $f \in W^{1,1}(\Rd)$)
then
\[
\widehat{(\partial_j f)}(\xi) = i\xi_j \fhatxi ,
\]
where $\partial_j = \partial/\partial x_j$ for $j=1,\ldots,d$. Thus:
\[
\text{$\widehat{\ }$ \ takes differentiation to multiplication by
$i\xi_j$.}
\]

(b) If $(1+|x|)f(x) \in \lord$ then $\fhat$ is continuously
differentiable, with
\[
\widehat{(-ix_j f)}(\xi) = (\partial_j \fhat)(\xi) ,
\]
where $\partial_j = \partial/\partial \xi_j$ for $j=1,\ldots,d$.
Thus:
\[
\text{$\widehat{\ }$ \ takes multiplication by $-ix_j$ to
differentiation.}
\]
\end{theorem}

\begin{proof} For (a)
\begin{align*}
\int_\Rd (\partial_j f)(x) e^{-i\xi x} \, dx & = \int_\Rd f(x)
(i\xi_j)e^{-i\xi x} \, dx && \text{by parts} \\
& = i\xi_j \fhatxi .
\end{align*}

For (b) we compute a difference quotient, with $\delta \in \R$ and
$e_j=$ unit vector in the $j$-th direction:
\begin{align*}
\frac{\fhat(\xi+\delta e_j)-\fhatxi}{\delta} & = \int_\Rd f(x)
e^{-i\xi x}
\frac{e^{-i\delta x_j}-1}{\delta} \, dx \\
& \to \int_\Rd f(x) e^{-i\xi x} (-ix_j) \, dx = \widehat{(-ix_j
f)}(\xi)
\end{align*}
as $\delta \to 0$, by dominated convergence with dominating function
$f(x)|x| \in \lord$. Hence $\fhatxi$ has partial derivative
$\widehat{(-ix_j f)}(\xi)$, which is continuous by
Theorem~\ref{th:bpr}.
\end{proof}

\begin{theorem}[Smoothness of $f$ and decay of $\fhat$] \label{th:sdr} \

(a) If $f \in \lord$ then $\fhatxi = o(1)$ as $|\xi| \to \infty$,
and
\[
|\fhatxi| \leq \lv f \rv_\lord = O(1) .
\]

(b) If $f \in C^1_c(\Rd)$ then $\fhatxi = o(1/|\xi|)$ as $|\xi| \to
\infty$, and
\[
|\fhatxi| \leq \frac{d \max_j \lv \partial_j f \rv_\lord}{|\xi|} =
O(1/|\xi|) .
\]
\end{theorem}

\begin{proof}

(a) Use Riemann--Lebesgue (Corollary~\ref{co:rlr}) and
Theorem~\ref{th:bpr}.

(b) For each $\xi$ there exists $j$ such that $|\xi_j| \geq |\xi|/d$
(since $|\xi_1| + \cdots + |\xi_d| \geq |\xi|$). Then
\begin{align*}
|\fhatxi| = \left| \frac{\widehat{(\partial_j f)}(\xi)}{i\xi_j}
\right| & \leq \frac{|\widehat{(\partial_j f)}(\xi)|}{|\xi|/d} \\
& \leq \frac{d \max_j |\widehat{(\partial_j f)}(\xi)|}{|\xi|} \\
& = o(1/|\xi|) && \text{by Riemann--Lebesgue} \\
& \leq \frac{d \max_j \lv \partial_j f \rv_\lord}{|\xi|} && \text{by
Theorem~\ref{th:bpr}} \\
& = O(1/|\xi|) .
\end{align*}
Or one could argue more directly using the gradient vector:
\begin{align*}
|\fhatxi| & = \frac{|\widehat{(\nabla f)}(\xi)|}{|i\xi|}
= o(1/|\xi|) && \text{by Riemann--Lebesgue} \\
& \leq \frac{\lv \nabla f \rv_\lord}{|\xi|} && \text{by
Theorem~\ref{th:bpr}} \\
& = O(1/|\xi|) .
\end{align*}
\end{proof}

\subsubsection*{Convolution}

\begin{definition} Given $f,g \in \lord$, define their
\emph{convolution}
\[
(f*g)(x) = \int_\Rd f(x-y) g(y) \, dy , \qquad x \in \Rd.
\]
\end{definition}
\begin{theorem}[Convolution and Fourier transforms] \label{th:cft}
If $f,g \in \lord$ then $f * g \in \lord$ with
\[
\lv f*g \rv_\lord \leq \lv f \rv_\lord \lv g \rv_\lord
\]
and
\[
\widehat{(f*g)}(\xi) = \fhatxi \ghatxi , \qquad \xi \in \Rd.
\]
\end{theorem}
Thus the Fourier transform takes convolution to multiplication.
\begin{proof}
Like Theorem~\ref{th:cfc}.
\end{proof}

\begin{example}
Let $f=\charfn_{[-1/2,1/2]}$, so that
$(f*f)(x)=(1-|x|)\charfn_{[-1,1]}(x)$ by direct calculation. We find
$\fhatxi=\sinc(\xi/2)$ like example 1 of Table~\ref{ta:ftexamples},
and $\widehat{(f*f)}(\xi)=\sinc^2(\xi/2)$ by example 2 of
Table~\ref{ta:ftexamples}.

Hence $\widehat{(f*f)}=(\fhat)^2$, as Theorem~\ref{th:cft} predicts.

As this example illustrates, \textbf{convolution is a smoothing
operation}, and hence improves the decay of the transform:
$\sinc(\xi/2)$ decays like $1/\xi$ while $\sinc^2(\xi/2)$ decays
like $1/\xi^2$.
\end{example}

\paragraph*{Convolution facts} (similar to Chapter~\ref{ch:fs})

\vspace{6pt} \noindent 1. Convolution is commutative: $f*g=g*f$. It
is also associative, and linear with respect to $f$ and $g$.

\vspace{6pt} \noindent 2. If $f \in \lprd, 1 \leq p \leq \infty$,
and $g \in \lord$, then $f*g \in \lprd$ with
\[
\lv f*g \rv_\lprd \leq \lv f \rv_\lprd \lv g \rv_\lord .
\]
Further, if $f \in C_0(\Rd)$ and $g \in \lord$ then $f*g \in C_0(\Rd)$.

\noindent \Proof For the first claim, use Young's Theorem~\ref{th:yt}. For the second, if $f \in C_0(\Rd)$ and $g \in \lord$ then $f*g$ is continuous because $(f * g)(x+z) \to (f*g)(x)$ as $z \to 0$ by uniform continuity of $f$ (exercise). And $(f*g)(x) \to 0$ as $|x| \to \infty$ by dominated convergence, since $f(x-y) \to 0$ as $|x| \to \infty$.

\vspace{6pt} \noindent 3. Convolution is continuous on $\lprd$: if
$f_m \to f$ in $\lprd, 1 \leq p \leq \infty$, and $g \in \lord$,
then $f_m*g \to f*g$ in $\lprd$.

\noindent \Proof Use linearity and Fact 2.

\vspace{6pt} \noindent 4. If $f \in \lord$ and $P(x) = \int_\Rd
Q(\xi) e^{i\xi x} \,  d\xi$ for some $Q \in \lord$, then
\begin{equation} \label{eq:convreal}
(P*f)(x) = \int_\Rd Q(\xi) \fhatxi e^{i\xi x} \, d\xi .
\end{equation}
\noindent \Proof
\begin{align*}
(P*f)(x) & = \int_\Rd Q(\xi) \int_\Rd e^{i\xi (x-y)} f(y) \, dy d\xi && \text{by Fubini} \\
& = \int_\Rd Q(\xi) e^{i\xi x} \fhatxi \, d\xi .
\end{align*}

\chapter{Fourier integrals: summability in norm} \label{ch:fis}

\subsubsection*{Goal}

Develop summability kernels in $\lprd$

\subsubsection*{Reference} \cite{K} Section VI.1

\begin{definition}
A \textbf{summability kernel} on $\Rd$ is a family $\{ k_\omega \}$
of integrable functions such that
\begin{align*}
\int_\Rd k_\omega(x) \, dx & = 1 && \textsc{(Normalization)} \label{eq:SR1} \tag{SR1} \\
\sup_\omega \int_\Rd |k_\omega(x)| \, dx & < \infty &&
\textsc{($L^1$ bound)} \label{eq:SR2}
\tag{SR2} \\
\lim_{\omega \to \infty} \int_{\{ x : |x|>\delta \}} |k_\omega(x)|
\, dx & = 0 && \textsc{($L^1$ concentration)} \label{eq:SR3} \tag{SR3}\\
& && \text{for each $\delta>0$.}
\end{align*}
Some kernels further satisfy
\begin{align*}
\lim_{\omega \to \infty} \sup_{|x|>\delta} |k_\omega(x)| & = 0 && \textsc{($L^\infty$ concentration)} \label{eq:SR4} \tag{SR4}\\
& && \text{for each $\delta>0$.}
\end{align*}
(\emph{Notation.} Here $k_\omega(x)$ does \emph{not} mean the
translation $k(x-\omega)$.)
\end{definition}

\begin{example}
Suppose $k \in \lord$ is continuous with $\int_\Rd k(x) \, dx = 1$.
Put
\[
k_\omega(x) = \omega^d k(\omega x)
\]
for $\omega>0$. Then $\{ k_\omega \}$ is a summability kernel.

\noindent \Proof Show \eqref{eq:SR1} and \eqref{eq:SR2} by changing
variable with $y=\omega x, dy=\omega^d dx$. For \eqref{eq:SR3},
\begin{align*}
\int_{\{ x : |x|>\delta \}} |k_\omega(x)| \, dx & = \int_{\{ y :
|y|>\omega \delta \}} |k(y)| \, dy \\
& \to 0
\end{align*}
as $\omega \to \infty$, by dominated convergence.
\end{example}

\begin{example}
For $d=1$, let
\begin{align}
D(x) & = \pr \int_\R \charfn_{[-1,1]}(\xi) e^{i\xi x} \, d\xi
\label{eq:DR1} \\
& = \frac{\sin x}{\pi x} = \frac{1}{\pi} \sinc x . \label{eq:DR2}
\end{align}
\begin{figure}
\begin{center}
  \includegraphics[scale=0.8]{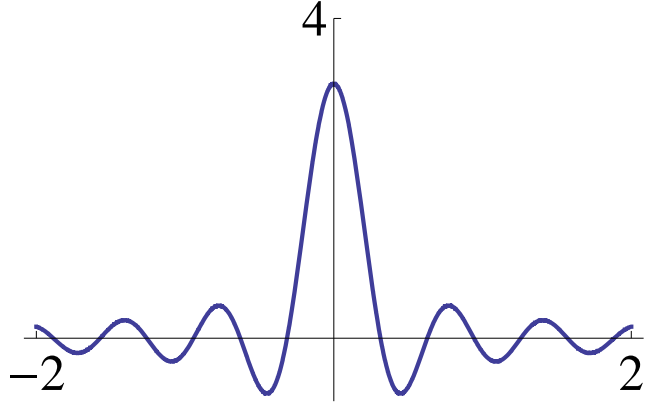}
  \caption{\label{DRfig}
    Dirichlet kernel with $\omega=10$}
\end{center}
\end{figure}
The \textbf{Dirichlet kernel} is
\begin{align}
D_\omega(x) = \omega D(\omega x) & = \pr \int_{-\omega}^\omega
e^{i\xi x} \, d\xi \label{eq:DR3} \\
& = \frac{\sin (\omega x)}{\pi x} . \label{eq:DR3b}
\end{align}
See Figure~\ref{DRfig}. $D$ is not integrable since $|D(x)| \sim
|x|^{-1}$ at infinity.

\noindent $\therefore \{ D_\omega \}$ is not a summability kernel.

In higher dimensions the Dirichlet function is $\prod_{j=1}^d
D(x_j)$, with associated kernel $D_\omega(x)=\prod_{j=1}^d
D_\omega(x_j)$.
\end{example}

\begin{example}
For $d=1$, let
\begin{align}
F(x) & = \pr \int_\R (1-|\xi|) \charfn_{[-1,1]}(\xi) e^{i\xi x} \,
d\xi
\label{eq:FR1} \\
& = \pr \left( \frac{\sin \big( \frac{1}{2} x \big)}{\frac{1}{2} x}
\right)^{\! \! \!2} && \text{by Table~\ref{ta:ftexamples}.}
\label{eq:FR2}
\end{align}
The \textbf{\Fejer kernel} is
\begin{align}
F_\omega(x) = \omega F(\omega x) & = \pr \int_{-\omega}^\omega
(1-|\xi|/\omega) e^{i\xi x} \, d\xi \label{eq:FR3} \\
& = \frac{\omega}{2\pi} \left( \frac{\sin \big( \frac{1}{2} \omega
x \big)}{\frac{1}{2} \omega x} \right)^{\! \! \!2} .
\label{eq:FR3b}
\end{align}
See Figure~\ref{FRfig}. $F$ is integrable since $F(x) \sim x^{-2}$
at infinity. And
\begin{figure}
\begin{center}
  \includegraphics[scale=0.8]{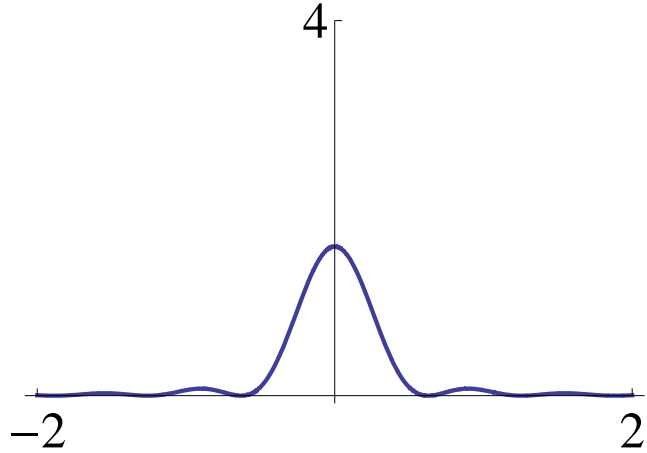}
  \caption{\label{FRfig}
    \Fejer kernel with $\omega=10$}
\end{center}
\end{figure}
\begin{align*}
\int_\R F(x) \, dx & = \frac{2}{\pi} \lim_{\rho \to \infty}
\int_{-\rho}^\rho \frac{\sin^2(x/2)}{x^2} \, dx \\
& = \frac{2}{\pi} \lim_{\rho \to \infty} \int_{-\rho}^\rho
\frac{2\sin(x/2)\cos(x/2)\cdot
(1/2)}{x} \, dx && \text{by parts} \\
& = \frac{1}{\pi} \lim_{\rho \to \infty}
\int_{-\rho}^\rho \frac{\sin x}{x} \, dx \\
& = 1 .
\end{align*}
$\therefore \{ F_\omega \}$ is a summability kernel.

In higher dimensions the \Fejer function is $\prod_{j=1}^d F(x_j)$,
with associated kernel $F_\omega(x)=\prod_{j=1}^d F_\omega(x_j)$.

\vspace{6pt} The \Fejer kernel is an arithmetic mean of Dirichlet
kernels; for example, $F(x)=\int_0^1 D_\omega(x) \,  d\omega$ in
$1$ dimension, by integrating \eqref{eq:DR3}.
\end{example}

\begin{example}
\begin{align}
P(x) & = \frac{1}{(2\pi)^d} \int_\Rd e^{-|\xi|} e^{i\xi x} \, d\xi
\label{eq:PR1} \\
& = \frac{c_d}{(1+|x|^2)^{(d+1)/2}} && \text{by
Table~\ref{ta:ftexamples}.} \label{eq:PR2}
\end{align}
The \textbf{Poisson kernel} is
\begin{figure}
\begin{center}
  \includegraphics[scale=0.8]{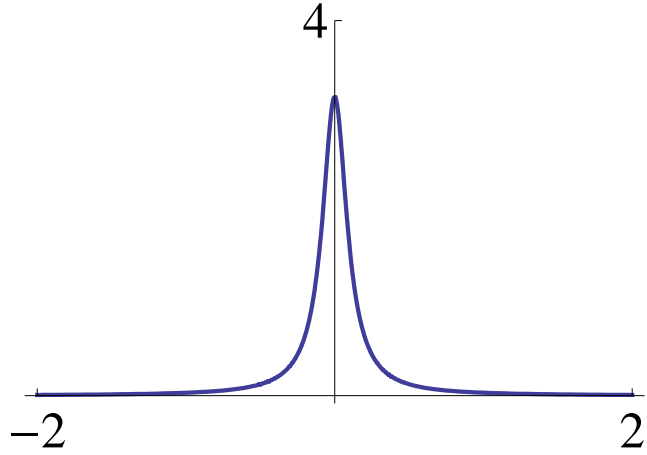}
  \caption{\label{PRfig}
    Poisson kernel with $\omega=10$}
\end{center}
\end{figure}
\begin{align}
P_\omega(x) = \omega^d P(\omega x) & = \frac{1}{(2\pi)^d}
\int_\Rd e^{-|\xi|/\omega} e^{i\xi x} \, d\xi \label{eq:PR3} \\
& = c_d \frac{\omega^{-1}}{\big( |x|^2+\omega^{-2}
\big)^{(d+1)/2}} . \label{eq:PR3b}
\end{align}
See Figure~\ref{PRfig}. $P$ is integrable since $P(x) \sim
|x|^{-(d+1)}$ at infinity. And $\int_\Rd P(x) \, dx =
\widehat{P}(0)=1$ because $\widehat{P}(\xi)=e^{-|\xi|}$ by
Example~\ref{ex:ftkernels} below; alternatively, one can integrate
\eqref{eq:PR2} directly (see \cite[p.~9]{SW} for $d>1$).

\noindent $\therefore \{ P_\omega \}$ is a summability kernel.
\end{example}

\begin{example}
\begin{align}
G(x) & = \frac{1}{(2\pi)^d} \int_\Rd e^{-|\xi|^2/2} e^{i\xi x} \,
d\xi
\label{eq:GR1} \\
& = (2\pi)^{-d/2} e^{-|x|^2/2} && \text{by
Table~\ref{ta:ftexamples}.} \label{eq:GR2}
\end{align}
The \textbf{Gauss kernel} is
\begin{align}
G_\omega(x) = \omega^d G(\omega x) & = \frac{1}{(2\pi)^d}
\int_\Rd e^{-|\xi/\omega|^2/2} e^{i\xi x} \,
d\xi \label{eq:GR3} \\
& = \frac{\omega^d}{(2\pi)^{d/2}} e^{-|\omega x|^2/2} .
\label{eq:GR3b}
\end{align}
See Figure~\ref{GRfig}. $G$ is clearly integrable, and $\int_\Rd
G(x) \, dx = 1$ from \eqref{eq:GR2}.

\begin{figure}[h]
\begin{center}
  \includegraphics[scale=0.8]{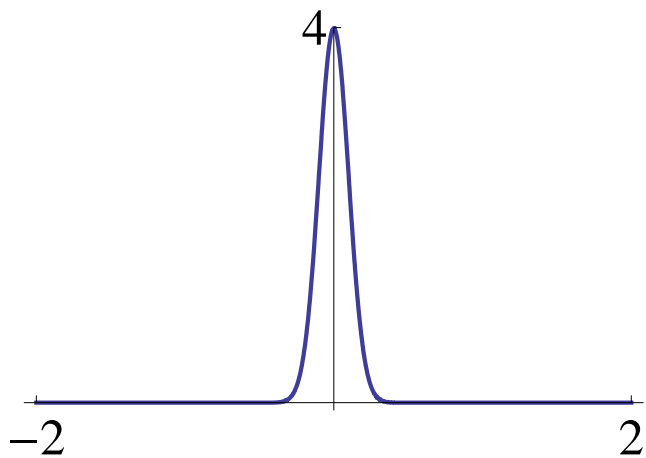}
  \caption{\label{GRfig}
    Gauss kernel with $\omega=10$}
\end{center}
\end{figure}
\noindent $\therefore \{ G_\omega \}$ is a summability kernel.
\end{example}

\subsubsection*{Connection to Fourier integrals}

For $f \in \lord$:
\begin{align}
(D_\omega * f)(x) & = \frac{1}{(2\pi)^d}
\int_{[-\omega,\omega]^d} \fhatxi
e^{i\xi x} \, d\xi \label{eq:DR4} \\
(F_\omega * f)(x) & = \frac{1}{(2\pi)^d}
\int_{[-\omega,\omega]^d} \big( \prod_{j=1}^d (1-|\xi_j|/\omega)
\big) \, \fhatxi e^{i\xi x} \, d\xi \label{eq:FR4} \\
(P_\omega * f)(x) & = \frac{1}{(2\pi)^d} \int_\Rd
e^{-|\xi|/\omega} \fhatxi e^{i\xi x} \, d\xi \label{eq:PR4} \\
(G_\omega * f)(x) & = \frac{1}{(2\pi)^d} \int_\Rd
e^{-|\xi/\omega|^2/2} \fhatxi e^{i\xi x} \, d\xi \label{eq:GR4}
\end{align}
\noindent \Proof Use Convolution Fact \eqref{eq:convreal} and
definitions \eqref{eq:DR1}, \eqref{eq:FR1}, \eqref{eq:PR1},
\eqref{eq:GR1}, respectively.

\emph{Caution.} The left sides of the above formulas make sense for
$f \in \lprd$, but the right side does not: so far we have defined
the Fourier transform only for $f \in \lord$.

\subsubsection*{Summability in norm}

\begin{theorem}[Summability in $\lprd$ and $C_0(\Rd)$] \label{th:srd}
Assume $\{ k_\omega \}$ is a summability kernel.

(a) If $f \in \lprd, 1 \leq p < \infty$, then $k_\omega * f \to f$
in $\lprd$ as $\omega \to \infty$.

(b) If $f \in C_0(\Rd)$ then $k_\omega * f \to f$ in $C_0(\Rd)$ as
$\omega \to \infty$.
\end{theorem}
Recall that $C_0(\Rd)$ uses the $\li$ norm.
\begin{proof}
Argue as for Theorem~\ref{th:sk}. Use that if $f \in C_0(\Rd)$ then $f$ is uniformly
continuous.
\end{proof}

\subsubsection*{Consequences}

\noindent $\bullet$ \Fejer summability for $f \in \lord$:
\begin{equation} \label{eq:fsr}
\frac{1}{(2\pi)^d} \int_{[-\omega,\omega]^d} \big( \prod_{j=1}^d
(1-|\xi_j|/\omega) \big) \, \fhatxi e^{i\xi x} \, d\xi \to f(x)
\qquad \text{in $\lord$.}
\end{equation}
Similarly for Poisson and Gauss summability.

\vspace{3pt} \noindent \Proof Use Theorem~\ref{th:srd} and formulas
\eqref{eq:FR4}--\eqref{eq:GR4}.

\vspace{6pt} \noindent $\bullet$ Uniqueness theorem:
\begin{equation} \label{eq:ur}
\text{if $f,g \in \lord$ with $\fhat=\ghat$ then $f=g$.}
\end{equation}
That is, the Fourier transform $\widehat{\ }: \lord \to \lird$ is
injective.

\vspace{3pt} \noindent \Proof Use \Fejer summability \eqref{eq:fsr}
on $f$ and $g$.

\subsubsection*{Connection to PDEs}

Fix $f \in \lord$.

\vspace{6pt} \noindent 1. The Poisson kernel solves Laplace's
equation in a half-space:
\begin{align*}
v(x,x_{d+1}) & = (P_{1/x_{d+1}}*f)(x) \\
& = c_d \int_\Rd \frac{x_{d+1}}{\big( |x-y|^2+x_{d+1}^2
\big)^{(d+1)/2}} f(y) \, dy
\end{align*}
solves
\[
(\partial_1^2 + \cdots + \partial_d^2 + \partial_{d+1}^2)
v=0
\]
on $\Rd \times (0,\infty)$, with boundary value $v(x,0)=f(x)$
in the sense of Theorem~\ref{th:srd}.

That is, $v$ is the harmonic extension of $f$ from $\Rd$ to the
halfspace $\Rd \times (0,\infty)$.

\vspace{3pt} \noindent \Proof Take $\omega = 1/x_{d+1}$ in
\eqref{eq:PR4} and differentiate through the integral, using
\begin{align*}
\sum_{j=1}^{d+1} \frac{\partial^2\ }{\partial x_j^2}
(e^{-|\xi|x_{d+1}} e^{i\xi x}) & = \big( (i\xi_1)^2 + \cdots +
(i\xi_d)^2 + (-|\xi|)^2 \big)
e^{-|\xi|x_{d+1}} e^{i\xi x} \\
& = 0 .
\end{align*}
For the boundary value, note $\omega = 1/x_{d+1} \to \infty$ as
$x_{d+1} \to 0$.

\vspace{6pt} \noindent 2. The Gauss kernel solves the diffusion
(heat) equation:
\begin{align*}
w(t,x) & = (G_{1/\sqrt{2t}}*f)(x) \\
& = \frac{1}{(4\pi t)^{d/2}} \int_\Rd e^{-|x-y|^2/4t} f(y) \, dy
\end{align*}
solves
\[
w_t=\Delta w
\]
for $(t,x) \in (0,\infty) \times \Rd$, with initial value $w(0,x)=f(x)$ in the sense of Theorem~\ref{th:srd}.
(Here $\Delta = \partial_1^2 + \cdots + \partial_d^2$.)

\vspace{3pt} \noindent \Proof Take $\omega = 1/\sqrt{2t}$ in
\eqref{eq:GR4} and differentiate through the integral, using
\begin{align*}
\left( \frac{\partial\ }{\partial t} - \sum_{j=1}^d
\frac{\partial^2\ }{\partial x_j^2} \right) (e^{-|\xi|^2 t} e^{i\xi
x}) & = \big( -|\xi|^2 - (i\xi_1)^2 - \cdots - (i\xi_d)^2 \big)
e^{-|\xi|^2 t} e^{i\xi x} \\
& = 0 .
\end{align*}
For the boundary value, note $\omega = 1/\sqrt{2t} \to \infty$ as
$t \to 0$.

\chapter[Fourier inversion when $\fhat \in \lord$]{Fourier transforms in $\lord$, and Fourier inversion} \label{ch:ftfi}

\subsubsection*{Goal} Fourier inversion when $\fhat$ is integrable

\subsubsection*{Reference} \cite[Section~VI.1]{K}

\vspace{18pt}
\begin{definition} \label{de:fin}
Define
\begin{align*}
\gcheck(x) & = \frac{1}{(2\pi)^d} \int_\Rd g(\xi) e^{i\xi x} \,
d\xi \\
& = \frac{1}{(2\pi)^d} \ghat(-x) .
\end{align*}
We call $\check{\ }$ the \emph{inverse Fourier transform}, in view
of the next theorem.
\end{definition}

\begin{theorem}(Fourier inversion) \label{th:fi} \

(a) If $f, \fhat \in \lord$ then $f$ is continuous and
\[
f(x) = \frac{1}{(2\pi)^d} \int_\Rd \fhatxi e^{i\xi x} \, d\xi ,
\qquad x \in \Rd .
\]

(b) If $g, \gcheck \in \lord$ then $g$ is continuous and
\[
g(\xi) = \int_\Rd \gcheck(x) e^{-i\xi x} \, dx , \qquad \xi \in \Rd
.
\]
\end{theorem}

The theorem says \fbox{$(\hat{f})\check{\ }=f$} and
\fbox{$(\gcheck)\hat{\ }=g$}.

\begin{proof}
(a) The $L^1$ convergence in \Fejer summability \eqref{eq:fsr}
implies pointwise convergence a.e.\ for some subsequence of
$\omega$-values:
\begin{align*}
f(x) & = \lim_{\omega \to \infty} \frac{1}{(2\pi)^d} \int_\Rd
\charfn_{[-\omega,\omega]^d}(\xi) \big( \prod_{j=1}^d
(1-|\xi_j|/\omega) \big) \fhatxi e^{i\xi x} \, d\xi \\
& = \frac{1}{(2\pi)^d} \int_\Rd \fhatxi e^{i\xi x} \, d\xi
\end{align*}
by dominated convergence, using that $\fhat \in \lord$.

(b) Apply part (a) to $g$, change $\xi \mapsto -\xi$, and then swap
$x$ and $\xi$.
\end{proof}

\begin{example} \label{ex:ftkernels}
The Fourier transforms of the Fej\'{e}r, Poisson and Gauss functions
can be computed by Fourier Inversion Theorem~\ref{th:fi}(b), because
definitions \eqref{eq:FR1}, \eqref{eq:PR1} and \eqref{eq:GR1}
express those kernels as inverse Fourier transforms. For example, if
we choose $g(\xi) = e^{-|\xi|^2/2}$ then definition \eqref{eq:GR1}
says $G(x) = \check{g}(x)$, so that $\widehat{G}=g$ by
Theorem~\ref{th:fi}(b).

Table~\ref{ta:ftkernels} displays the results.
\begin{table}[t]
\begin{center}
\begin{tabular}{|c|l|l|}
\hline & & \\
dimension    & $f(x)$    & $\fhatxi$
  \\
\hline & & \\ $d$ & $F(x) = \frac{1}{(2\pi)^d } \prod_{j=1}^d \left(
\frac{\sin (
x_j/2)}{x_j/2} \right)^{\! \! 2}$ & $\widehat{F}(\xi) = \charfn_{[-1,1]^d}(\xi) \prod_{j=1}^d (1-|\xi_j|)$ \\
& & \\
$d$ & $P(x) = \frac{c_d}{(1+|x|^2)^{(d+1)/2}}$ & $\widehat{P}(\xi) = e^{-|\xi|}$ \\
& & \\ $d$ & $G(x) = (2\pi)^{-d/2} e^{-|x|^2/2}$ & $\widehat{G}(\xi) = e^{-|\xi|^2/2}$ \\
\hline
\end{tabular}
\vspace*{3pt} \caption{Fourier transforms of the Fej\'{e}r, Poisson
and Gauss functions, from Example~\ref{ex:ftkernels}.}
\label{ta:ftkernels}
\end{center}
\end{table}
\end{example}

\chapter{Fourier transforms in $\ltrd$} \label{ch:ftl2}

\subsubsection*{Goal} Extend the Fourier transform to an isometric
bijection of $\ltrd$ to itself

\subsubsection*{Reference} \cite{K} Section VI.3

\subsubsection*{Notation} Inner product on $\ltrd$ is $\la f , g \ra = \int_\Rd f(x) \overline{g(x)} \,
dx$.

\begin{theorem}[Fourier transform on $\ltrd$] \label{th:ftl2}
The Fourier transform $\widehat{\ } : \ltrd \to \ltrd$ is a
bijective isometry (up to a constant factor) with
\begin{align*}
\lv f \rv_\ltrd & = (2\pi)^{-d/2} \lv \fhat \rv_\ltrd &&
\textsc{(Plancherel)} \\
\la f , g \ra & = (2\pi)^{-d} \la \fhat , \ghat \ra &&
\textsc{(Parseval)} \\
(\hat{f})\check{\ }=f, & \quad (\gcheck)\hat{\ }=g &&
\textsc{(Inversion)}
\end{align*}
for all $f,g \in \ltrd$.
\end{theorem}
The proof will show $\widehat{\ } : L^1 \cap L^2(\Rd) \to \ltrd$ is
bounded with respect to the $L^2$ norm. Then by density of $L^1 \cap
L^2$ in $L^2$, we conclude the Fourier transform extends to a
bounded operator from $L^2$ to itself.

\begin{proof}
For $f \in L^1 \cap L^2(\Rd)$,
\begin{align}
\lv f \rv_\ltrd^2 & = \lim_{\omega \to \infty} \int_\Rd f(x)
\overline{(G_\omega * f)(x)} \, dx \notag \\
& \qquad \qquad \qquad \text{since $G_\omega * f
\to f$ in $\ltrd$ by Theorem~\ref{th:srd}} \notag \\
& = \lim_{\omega \to \infty} \frac{1}{(2\pi)^d} \int_\Rd \int_\Rd
f(x) e^{-i\xi x} \overline{\fhatxi} e^{-|\xi/\omega|^2/2}
\, d\xi dx \qquad \text{by \eqref{eq:GR4}} \notag \\
& = \lim_{\omega \to \infty} \frac{1}{(2\pi)^d} \int_\Rd
|\fhatxi|^2 e^{-|\xi/\omega|^2/2} \, d\xi \quad \text{by Fubini, using $\fhat \in \lird$,} \notag \\
& = \frac{1}{(2\pi)^d} \int_\Rd |\fhatxi|^2 \, d\xi \quad \qquad \qquad \qquad \text{by monotone convergence} \notag \\
& = \frac{1}{(2\pi)^d} \lv \fhat \rv_\ltrd^2 . \label{eq:ftl2}
\end{align}
By density of $L^1 \cap L^2$ in $L^2$, the Fourier transform
$\widehat{\ }$ extends to a bounded operator from $\ltrd$ to itself.
Plancherel follows from \eqref{eq:ftl2} by density. Thus the Fourier
transform is an isometry, up to a constant factor.

Parseval follows from Plancherel by polarization, or by repeating
the argument for Plancherel with $\la f ,f \ra$ changed to $\la f ,
g \ra$ (and using dominated instead of monotone convergence).

For Inversion, note $\check{\ } : \ltrd \to \ltrd$ is bounded by
Definition~\ref{de:fin}, since the Fourier transform is bounded. If
$f$ is smooth with compact support then $\fhat$ is bounded and
decays rapidly at infinity, by repeated use of Theorem~\ref{th:sdr}.
Hence $\fhat \in \lord$, with $(\hat{f})\check{\ }=f$ by Inversion
Theorem~\ref{th:fi}. So the Fourier transform followed by the
inverse transform gives the identity on the dense set $L^1 \cap
L^2(\Rd)$, and hence on all of $\ltrd$ by continuity. Similarly
$(\gcheck)\hat{\ }=g$ for all $g \in \ltrd$.

Finally, the Fourier transform is injective by Plancherel, and
surjective by Inversion.
\end{proof}

\begin{example} \label{ex:ftdirichlet}
In $1$ dimension, the Dirichlet function
\[
D(x) = \frac{\sin x}{\pi x}
\]
belongs to $\ltr$ and has
\[
\widehat{D}(\xi) = \charfn_{[-1,1]}(\xi) .
\]
\noindent \Proof $D = (\charfn_{[-1,1]})\check{\ }$ by definition in
\eqref{eq:DR1}, and so $\widehat{D} = \charfn_{[-1,1]}$ by
Theorem~\ref{th:ftl2} Inversion.

\begin{table}[t]
\begin{center}
\begin{tabular}{|c|l|l|}
\hline & & \\
dimension    & $f(x)$    & $\fhatxi$
  \\
\hline & & \\ $d$ & $D(x) = \frac{1}{\pi^d } \prod_{j=1}^d \frac{\sin x_j}{x_j}$ & $\widehat{D}(\xi) = \charfn_{[-1,1]^d}(\xi)$ \\
& & \\
\hline
\end{tabular}
\vspace*{3pt} \caption{Fourier transform of the Dirichlet function,
from Example~\ref{ex:ftdirichlet}.} \label{ta:ftdirichlet}
\end{center}
\end{table}
\end{example}

\begin{remark}\rm
If $f \in \ltrd$ then $f \charfn_{B(n)} \in L^1 \cap L^2(\Rd)$ and
$f \charfn_{B(n)} \to f$ in $\ltrd$. Hence
\begin{align*}
\fhatxi & = \lim_{n \to \infty} \widehat{(f \charfn_{B(n)})}(\xi) &&
\text{in $\ltrd$, by Theorem~\ref{th:ftl2},} \\
& = \lim_{n \to \infty} \int_{B(n)} f(x) e^{-i\xi x} \, dx .
\end{align*}
How can this limit exist, when $f$ need not be integrable? The
answer must be that oscillations of $e^{-i\xi x}$ yield cancelations
that allow $f(x) e^{-i\xi x}$ to be integrated improperly, as above,
for almost every $\xi$.
\end{remark}

\begin{theorem}[Hausdorff--Young for Fourier transform] \label{th:hyft}
The Fourier transform
\[
\widehat{\ } : \lprd \to L^{p^\prime}(\Rd)
\]
is bounded for $1 \leq p \leq 2$, where $\frac{1}{p} +
\frac{1}{p^\prime} = 1$.
\end{theorem}
\begin{proof}
Apply the Riesz--Thorin Interpolation Theorem~\ref{th:rti}, using
boundedness of
\begin{align*}
\widehat{\ } & : \lord \to \lird && \text{in Theorem~\ref{th:bpr}, \ and} \\
\widehat{\ } & : \ltrd \to \ltrd && \text{in Theorem~\ref{th:ftl2}.}
\end{align*}
Note the Fourier transform is well defined on $L^1+L^2(\Rd)$, since
the $L^1$ and $L^2$ Fourier transforms agree on $L^1 \cap L^2(\Rd)$.
\end{proof}

\begin{remark}\rm \label{re:bpp}
The first five Basic Properties in Theorem~\ref{th:bpr} still hold
for the Fourier transform on $\lprd, 1 \leq p \leq 2$, and so do
Corollary~\ref{co:rf} (radial functions) and Lemma~\ref{le:tp}
(product functions) and \eqref{eq:DR4}--\eqref{eq:GR4} (connection
to Fourier integrals).

\noindent \Proof Given $f \in \lprd$, take $f_m \in L^1 \cap
L^p(\Rd)$ with $f_m \to f$ in $\lprd$. Then $\widehat{f_m} \to
\fhat$ in $L^{p^\prime}(\Rd)$ by the Hausdorff--Young
Theorem~\ref{th:hyft}. Here $\widehat{f_m}$ is the usual Fourier
transform of $f_m \in \lord$, so that Theorem~\ref{th:bpr},
Corollary~\ref{co:rf}, Lemma~\ref{le:tp} and
\eqref{eq:DR4}--\eqref{eq:GR4} all apply to $f_m$. Now let $m \to
\infty$ in those results.
\end{remark}

\begin{corollary}[Convolution and Fourier transforms] \label{co:cftr}
If $f \in \lord, g \in \lprd, 1 \leq p \leq 2$, then $f*g \in \lprd$
and
\[
\widehat{(f*g)} = \fhat \ \ghat .
\]
\end{corollary}

\begin{proof}
Take $g_m \in L^1 \cap L^p(\Rd)$ with $g_m \to g$ in $\lprd$. Then
$\widehat{(f*g_m)} = \fhat \, \widehat{g_m}$ by
Theorem~\ref{th:ftl2}. Let $m \to \infty$ and use the
Hausdorff--Young Theorem~\ref{th:hyft}, noting $\fhat$ is bounded.
\end{proof}

\subsubsection*{Consequence}

Analogue of Weierstrass trigonometric approximation: functions with
compactly supported Fourier transform are dense in $\lprd, 1 \leq p
\leq 2$.

\noindent \Proof $F_\omega * f \to f$ in $\lprd$ by
Theorem~\ref{th:srd}, and $\widehat{(F_\omega * f)} =
\widehat{F_\omega} \fhat$ has compact support (because
$\widehat{F_\omega}$ has compact support by
Table~\ref{ta:ftkernels}).

\chapter[Fourier integrals: summability a.e.]{Fourier integrals: summability pointwise} \label{ch:fisp}

\subsubsection*{Goal}

Prove sufficient conditions for summability at a single point, and
a.e.

\subsubsection*{Reference}

\cite{G} Sections 2.1b, 3.3b

\vspace{18pt}

If $f \in C_0(\Rd)$ then $k_\omega * f \to f$ uniformly by
Theorem~\ref{th:srd}(b), and hence convergence holds at every $x$.
But what if $f$ is merely continuous at a \emph{point}?

\begin{theorem}[Summability at a point] \label{th:spr}
Assume $\{ k_\omega \}$ is a summability kernel. Suppose either $f
\in \lord$ and $\{ k_\omega \}$ satisfies the $L^\infty$
concentration hypothesis \eqref{eq:SR4}, or else $f \in \lird$.

If $f$ is continuous at $x_0 \in \Rd$ then $(k_\omega * f)(x_0) \to
f(x_0)$ as $\omega \to \infty$.
\end{theorem}
\begin{proof}
Adapt the corresponding result on the torus, Theorem~\ref{th:sp}(a).
\end{proof}

The Poisson and Gauss kernels satisfy \eqref{eq:SR4}, and so does
the Fej\'{e}r kernel in $1$ dimension. More generally, if $k(x) =
o(1/|x|^d)$ as $|x| \to \infty$ then $k_\omega(x)=\omega^d
k(\omega x)$ satisfies \eqref{eq:SR4} (Exercise).

\vspace{6pt} Next we aim at summability a.e., by using maximal
functions like we did for Fourier series in Chapter~\ref{ch:sp}.

\begin{definition}
Define the
\begin{align*}
\text{\emph{Dirichlet} maximal function} \quad (D^* f)(x) & =
\sup_\omega
|(D_\omega * f)(x)| \\
\text{\emph{\Fejer} maximal function} \quad (F^* f)(x) & =
\sup_\omega
|(F_\omega * f)(x)| \\
\text{\emph{Poisson} maximal function} \quad (P^* f)(x) & =
\sup_\omega
|(P_\omega * f)(x)| \\
\text{\emph{Gauss} maximal function} \quad (G^* f)(x) & =
\sup_\omega
|(G_\omega * f)(x)| \\
\text{\emph{Lebesgue} maximal function} \quad (L^* f)(x) & =
\sup_\omega |(L_\omega * f)(x)|
\end{align*}
where
\[
L(x) = \frac{1}{|B(1)|} \charfn_{B(1)}(x)
\]
is the normalized indicator function of the unit ball.
\end{definition}

\begin{lemma} \label{le:lm}
\[
L^* f \leq L^* |f| = Mf
\]
where $M$ is the Hardy--Littlewood maximal operator from
Chapter~\ref{ch:mf}.
\end{lemma}
\begin{proof}
First,
\begin{equation} \label{eq:Lform}
L_{1/\omega}(y) = (1/\omega)^d L(y/\omega) =
\frac{1}{|B(\omega)|} \charfn_{B(\omega)}(y) .
\end{equation}
Hence
\begin{align*}
|(L_{1/\omega} * f)(x)| \leq (L_{1/\omega} * |f|)(x) & = \frac{1}{|B(\omega)|} \int_{B(\omega)} |f(x-y)| \, dy ,
\end{align*}
and taking the supremum over $\omega$ gives $(Mf)(x)$. 
\end{proof}

\begin{lemma}[Majorization] \label{le:majr}
If $k \in \lord$ is nonnegative and radially symmetric decreasing,
then
\[
|(k*f)(x)| \leq \lv k \rv_\lord (L^*f)(x) \qquad \text{for all $x
\in \Rd, \quad f \in \lord$.}
\]
\end{lemma}
\begin{proof}
Write $k(x)=\rho(|x|)$ where $\rho : [0,\infty) \to \R$ is
nonnegative and decreasing. Assume $\rho$ is absolutely continuous,
for simplicity. We first establish a layer-cake decomposition of
$k$, like we did on the torus in Lemma~\ref{le:maj}:
\begin{align*}
k(y) = \rho(|y|) & = -\int_{|y|}^\infty \rho^\prime(\omega) \,
d\omega \qquad \text{since $\rho(\infty)=0$ by integrability of $k$} \\
& = -\int_0^\infty |B(\omega)| L_{1/\omega}(y)
\rho^\prime(\omega) \, d\omega ,
\end{align*}
because by \eqref{eq:Lform},
\[
L_{1/\omega}(y) =
\begin{cases}
1/|B(\omega)| & \text{if $\omega>|y|$,} \\
0 & \text{if $\omega \leq |y|$.}
\end{cases}
\]
Hence
\[
(k*f)(x) = \int_0^\infty |B(\omega)| (L_{1/\omega} * f)(x)
\big(-\rho^\prime(\omega)\big) \, d\omega
\]
and so
\begin{align*}
|(k*f)(x)| & \leq \int_0^\infty |B(\omega)|
\big(-\rho^\prime(\omega)\big)
\, d\omega \cdot (L^* f)(x) \qquad \qquad \text{since $\rho^\prime \leq 0$} \\
& = \int_0^\infty \int_0^\omega |\partial B(1)| r^{d-1} \, dr \,
\big(-\rho^\prime(\omega)\big) \, d\omega \cdot (L^* f)(x) \\
& \qquad \qquad \qquad \text{by spherical coordinates for $|B(\omega)|=\int_{B(\omega)} \, dy$} \\
& = \int_0^\infty |\partial B(1)| \omega^{d-1} \rho(\omega) \, d\omega \cdot (L^* f)(x) \\
& \qquad \text{by parts with respect to $\omega$ (why does the $\omega=\infty$ term vanish?)} \\
& = \int_\Rd k(y) \, dy \cdot (L^* f)(x)
\end{align*}
by using spherical coordinates again.
\end{proof}

\begin{theorem}[Lebesgue dominates Poisson and Gauss in all dimensions, and \Fejer in $1$
dimension] \label{th:lpgf}
\begin{align*}
F^* f & \leq \frac{4}{\pi} L^* |f| \qquad \text{(when $d=1$)} \\
P^* f & \leq L^* f \\
G^* f & \leq L^* f
\end{align*}
for all $f \in \lprd, 1 \leq p \leq \infty$.
\end{theorem}
\begin{proof}
$P^*f \leq L^*f$ by the Majorization Lemma~\ref{le:majr}, since
$P_\omega$ is nonnegative and radially symmetric decreasing, with
$\lv P_\omega \rv_\lord = 1$. Similarly $G^*f \leq L^*f$.

When $d=1$,
\begin{align*}
F_\omega(x) & = \frac{\omega}{2\pi} \left( \frac{\sin \big(
\frac{\omega}{2} x \big)}{\frac{\omega}{2} x} \right)^{\! \! \!
2} && \text{by \eqref{eq:FR3b}} \\
& \leq k(x) \overset{\text{def}}{=} \frac{\omega}{2\pi}
\begin{cases}
1 , & |x| \leq 2/\omega, \\
1/\big( \frac{\omega}{2} x \big)^2 , & |x|>2/\omega .
\end{cases}
\end{align*}
Note $k$ is nonnegative, even and decreasing, with $\lv k \rv_\lor =
4/\pi$. Hence $|F_\omega * f| \leq k*|f| \leq (4/\pi)L^*|f|$ by
Majorization Lemma~\ref{le:majr}.
\end{proof}

\begin{remark}\rm
The \Fejer kernel is \emph{not} majorized by a radially symmetric
decreasing integrable function, when $d \geq 2$. For example, taking
$\omega=2$ gives
\[
F_2(x) = \prod_{j=1}^d \frac{1}{\pi} \left( \frac{\sin x_j}{x_j}
\right)^{\! \! \! 2} ,
\]
which decays like $x_1^{-2}$ along the $x_1$-axis. Thus the best
possible radial bound would be $O(|x|^{-2})$, which is not
integrable at infinity in dimensions $d \geq 2$.
\end{remark}

\begin{corollary} \label{co:swr}
$F^*, P^*, G^*$ and $L^*$ are weak $(1,1)$ and strong $(p,p)$ on
$\lprd$, for $1<p \leq \infty$.
\end{corollary}
\begin{proof}
Combine Theorem~\ref{th:lpgf} and Lemma~\ref{le:lm} with the weak
and strong bounds on the Hardy--Littlewood maximal operator in
Chapter~\ref{ch:mf}.

For the \Fejer kernel in dimensions $d \geq 2$, see
\cite[Theorem~3.3.3]{G}.
\end{proof}

\begin{theorem}[Summability a.e.]
If $f \in \lprd, 1 \leq p \leq \infty$, then
\begin{align*}
F_\omega * f & \to f \ \text{a.e.\ as $\omega \to \infty$,} \\
P_\omega * f & \to f \ \text{a.e.\ as $\omega \to \infty$,} \\
G_\omega * f & \to f \ \text{a.e.\ as $\omega \to \infty$,} \\
L_\omega * f & \to f \ \text{a.e.\ as $\omega \to \infty$.}
\end{align*}
(The last statement is the Lebesgue differentiation theorem.)
\end{theorem}
\begin{proof}
Assume $1 \leq p < \infty$. $F^*$ is weak $(p,p)$ by
Corollary~\ref{co:swr}. Hence the Theorem in Chapter~\ref{ch:mf}
says
\[
{\mathcal C} = \{ f \in \lprd : \lim_{\omega \to \infty} F_\omega
* f = f \ \text{a.e.} \}
\]
is closed in $\lprd$. Obviously ${\mathcal C}$ contains every $f \in
C_c(\Rd)$, because $F_\omega * f \to f$ uniformly by
Theorem~\ref{th:srd}. Thus ${\mathcal C}$ is dense in $\lprd$ (using
here that $p<\infty$). Because ${\mathcal C}$ is closed, it must
equal $\lprd$, which proves the result.

When $p=\infty$, consider $f \in \lird$. For $m \in \N$, put
$g=\charfn_{B(m)} f$ and $h=f-g$. Then $g \in \lord$, and so
$F_\omega
* g \to g$ a.e., by the part of the theorem already proved. Hence
$F_\omega * g \to f$ a.e.\ on $B(m)$. Next $h \in \lird$ is
continuous on $B(m)$, with $h=0$ there, and so $F_\omega * h \to
h=0$ on $B(m)$ by Theorem~\ref{th:spr}. Since $f=g+h$ we deduce
$F_\omega * f \to f$ a.e.\ on $B(m)$. Letting $m \to \infty$ proves
the result.

Argue similarly for the other kernels.
\end{proof}

\chapter{Fourier integrals: norm convergence} \label{ch:fir}

\subsubsection*{Goal} Show norm convergence for $\lprd$ follows from
boundedness of the Hilbert transform on $\R$

\subsubsection*{Reference} I do not know a fully satisfactory reference for this material. Suggestions are welcome!

\vspace{18pt}

\begin{definition}\rm
Write
\[
S_\omega(f) = D_\omega * f
\]
where
\[
D_\omega(x)=\prod_{j=1}^d \omega D(\omega x_j) = \prod_{j=1}^d
\frac{\sin (\omega x_j)}{\pi x_j}
\]
is the Dirichlet kernel on $\Rd$ and $D(z)=(\sin z)/\pi z$ is the
Dirichlet function in $1$ dimension.
\end{definition}

$S_\omega$ is the ``partial sum'' operator for the Fourier
integral, because if $f \in \lprd,  1 \leq p \leq 2$, then
$S_\omega(f)=(\charfn_{[-\omega,\omega]^d} \hat{f} )\, \check{\
}$ by \eqref{eq:DR4} and Remark~\ref{re:bpp}. In particular,
\[
S_\omega : \ltrd \to \ltrd
\]
is bounded, by boundedness of the Fourier transform and its
inverse on $\lt$. Further, $\charfn_{[-\omega,\omega]^d} \fhat
\to \fhat$ and so $S_\omega(f) \to f$ in $\ltrd$, as $\omega \to
\infty$.

$S_\omega(f)$ is well defined whenever $f \in \lprd, 1 \leq p <
\infty$, because $D_\omega \in L^q(\Rd)$ for each $q>1$ and so
$D_\omega * f \in L^r(\Rd)$ for each $r \in (p,\infty]$, by the Generalized Young's
Theorem in Chapter~\ref{ch:ai}.

We will prove below that $S_\omega(f) \in \lprd$ when $f \in
\lprd, 1 < p < \infty$. But $S_\omega(f)$ need not belong to
$\lord$ when $f \in \lord$ (Exercise).

\vspace{6pt} Our goal in this Chapter is to improve the $\lp$
summability for Fourier integrals ($F_\omega * f \to f$ in
Theorem~\ref{th:srd}) to $\lp$ convergence ($D_\omega * f =
S_\omega(f) \to f$ in Theorem~\ref{th:ficn} below). As remarked
above, we have the result already for $p=2$.

First we reduce norm convergence to norm boundedness.

\begin{theorem} \label{th:red1}
Let $1 < p < \infty$ and suppose $\sup_\omega \lv S_\omega
\rv_{\lprd \to \lprd} < \infty$. Then Fourier integrals converge in
$\lprd$: $\lim_{\omega \to \infty} \lv S_\omega(f) - f \rv_\lprd =
0$ for each $f \in \lprd$.
\end{theorem}

\begin{proof}
Let ${\mathcal A} = \{ g \in L^1 \cap L^p(\Rd) : \text{$\ghat$ has
compact support} \}$. We claim ${\mathcal A}$ is dense in $\lprd$.
Indeed, if $f \in L^1 \cap L^p(\Rd)$ then $F_\omega * f \in L^1
\cap L^p(\Rd)$ and $\widehat{(F_\omega * f)} = \widehat{F_\omega}
\fhat$ has compact support by Table~\ref{ta:ftkernels}. Thus
$F_\omega * f \in {\mathcal A}$. Since $F_\omega * f \to f$ in
$\lprd$ by Theorem~\ref{th:srd}, and $L^1 \cap L^p$ is dense in
$L^p$, we see ${\mathcal A}$ is dense in $\lprd$.

We further show $S_\omega(g)=g$, when $g \in {\mathcal A}$,
provided $\omega$ is large enough that $[-\omega,\omega]^d$
contains the support of $\ghat$. To see this fact, note
$S_\omega(g) \in \ltrd$ because $D_\omega \in \ltrd$ and $g \in
\lord$; thus
\begin{align*}
\widehat{S_\omega(g)} & = \widehat{D_\omega} \widehat{g} \\
& = \charfn_{[-\omega,\omega]^d} \ghat && \text{by Table~\ref{ta:ftdirichlet}} \\
& = \ghat .
\end{align*}
Applying Fourier inversion in $\lt$ gives $S_\omega(g)=g$.

We conclude
\[
{\mathcal A} \subset \{ f \in \lprd : \lim_{\omega \to \infty}
S_\omega(f) = f \ \text{in $\lprd$} \} \overset{\text{def}}{=}
{\mathcal C} ,
\]
so that ${\mathcal C}$ is dense in $\lprd$. Because ${\mathcal C}$
is closed by Proposition~\ref{pr:red} (using the assumption that
$\sup_\omega \lv S_\omega \rv_{\lprd \to \lprd} < \infty$), we
conclude ${\mathcal C} = \lprd$, which proves the theorem.
\end{proof}

Next we reduce to norm boundedness in $1$ dimension. For the sake of
generality we allow different $\omega$-values in each coordinate
direction. (Thus our ``square partial sums'' for convergence of
Fourier integrals can be relaxed to ``rectangular partial sums'';
proof omitted.)

Given a vector $\vec{\omega} = (\omega_1, \ldots, \omega_d)$ of
positive numbers, define
\[
D_{\vec{\omega}}(x)=\prod_{j=1}^d \omega_j D(\omega_j x_j) .
\]
The Fourier multiplier
\[
\widehat{D_{\vec{\omega}}} = \charfn_{[-\omega_1,\omega_1] \times
\cdots \times [-\omega_d,\omega_d]}
\]
is the indicator function of a rectangular box.

Write
\[
C_{p,d} = \sup_{\vec{\omega}} \lv S_{\vec{\omega}} \rv_{\lprd \to
\lprd}
\]
for the norm bound on the partial sum operators. We have not yet
shown that this constant is finite.

\begin{theorem}[Reduction to $1$ dimension] \label{th:red2}
$C_{p,d} \leq (C_{p,1})^d$.
\end{theorem}

\begin{proof}
First observe that for $g \in \lpr$ and $\omega>0$,
\begin{align}
\int_\R \Big| \int_\R \omega D(\omega(x-y)) g(y) \, dy \Big|^p \,
dx & = \lv D_\omega * g \rv_\lpr^p \notag \\
& \leq C_{p,1}^p \lv g \rv_\lpr^p && \text{by definition of $C_{p,1}$} \notag \\
& = C_{p,1}^p \int_\R |g(y)|^p \, dy . \label{eq:red1}
\end{align}
Hence for $f \in L^p(\R^2)$ and
$\vec{\omega}=(\omega_1,\omega_2)$,
\begin{align*}
& \int_{\R^2} |(D_{\vec{\omega}}*f)(x_1,x_2)|^p \, dx_1 dx_2 \\
& = \int_\R \int_\R \Big| \int_\R \omega_1 D(\omega_1(x_1-y_1)) \int_\R \omega_2 D(\omega_2(x_2-y_2)) f(y_1,y_2) \, dy_2 dy_1 \Big|^p \, dx_1 dx_2 \\
& \leq C_{p,1}^p \int_\R \int_\R \big| \int_\R \omega_2 D(\omega_2(x_2-y_2)) f(y_1,y_2) \, dy_2 \big|^p \, dy_1 dx_2 \\
& \qquad \qquad \qquad \text{by \eqref{eq:red1} with $g(y_1)=\int_\R \omega_2 D(\omega_2(x_2-y_2)) f(y_1,y_2) \, dy_2$} \\
& \leq C_{p,1}^{2p} \int_\R \int_\R |f(y_1,y_2)|^p \, dy_2 dy_1 \\
& \qquad \qquad \qquad \text{by \eqref{eq:red1} with
$g(y_2)=f(y_1,y_2)$.}
\end{align*}
Taking $p$-th roots gives $\lv S_{\vec{\omega}} \rv_{L^p(\R^2) \to
L^p(\R^2)} \leq C_{p,1}^2$, which proves the theorem when $d=2$.

Argue similarly for $d \geq 3$.
\end{proof}

\noindent \emph{Aside.} The ``ball'' multiplier
$\charfn_{B(1)}(\xi)$ does \emph{not} yield a partial sum operator
with uniform norm bounds, when $p \neq 2$; see \cite[Section
10.1]{G}. Therefore Fourier integrals and series in higher
dimensions should be evaluated with ``rectangular'' partial sums,
and not ``spherical'' sums, when working in $\lp$ for $p \neq 2$.

\subsubsection*{Boundedness in $\lpr$}

\noindent 1. We shall prove (in Chapters~\ref{ch:hrt} and
\ref{ch:hpr}) the existence of a bounded linear operator
\[
H: \lpr \to \lpr , \qquad 1<p<\infty,
\]
called the \textbf{Hilbert transform} on $\R$, with the property
that
\[
\widehat{(Hf)}(\xi) = -i \sign(\xi) \fhatxi
\]
when $f \in L^p \cap L^2(\R)$. (Thus $H$ is a Fourier multiplier
operator.)

\vspace{6pt} \noindent 2. Then the \emph{Riesz projection} $P: \lpr
\to \lpr$ defined by
\[
Pf = \frac{1}{2}(f + iHf)
\]
is also bounded, when $1<p<\infty$.

Observe $P$ projects onto the positive frequencies:
\[
\widehat{(Pf)}(\xi) = \charfn_{(0,\infty)}(\xi) \fhatxi , \qquad f
\in \ltr ,
\]
since $i(-i \sign(\xi))=\sign(\xi)$.

\vspace{6pt} \noindent 3. The following formula expresses the
Fourier partial sum operator in terms of the Riesz projection and
some modulations: for $\omega>0$,
\begin{equation} \label{eq:rpr}
e^{-i\omega x}P(e^{i\omega x}f) - e^{i\omega x}P(e^{-i\omega
x}f) = S_\omega (f) , \qquad f \in \ltr .
\end{equation}
\noindent \Proof
\begin{align*}
[e^{i\omega x}f]\widehat{\ }(\xi) & = \fhat(\xi-\omega) \\
[P(e^{i\omega x}f)]\widehat{\ }(\xi) & = \charfn_{(0,\infty)}(\xi) \fhat(\xi-\omega) \\
[e^{-i\omega x} P(e^{i\omega x}f)]\widehat{\ }(\xi) & = \charfn_{(0,\infty)}(\omega+\xi) \fhatxi \\
& = \charfn_{(-\omega,\infty)}(\xi) \fhatxi \\
[e^{i\omega x} P(e^{-i\omega x}f)]\widehat{\ }(\xi) & =
\charfn_{(\omega,\infty)}(\xi) \fhatxi
\end{align*}
Subtracting the last two formulas gives
$\charfn_{(-\omega,\omega]} \fhat$, which equals
$\widehat{S_\omega(f)}$. Fourier inversion now completes the proof.

\vspace{6pt} \noindent 4. From \eqref{eq:rpr} applied to the dense
class of $f \in L^p \cap L^2(\R)$, and from boundedness of the Riesz
projection, it follows that
\[
C_{p,1} = \sup_\omega \lv S_\omega \rv_{\lpr \to \lpr} \leq 2 \lv
P \rv_{\lpr \to \lpr} < \infty
\]
when $1<p<\infty$. Hence from Theorems~\ref{th:red1} and
\ref{th:red2} we conclude:

\begin{theorem}[Fourier integrals converge in $\lprd$] \label{th:ficn} Let $1<p<\infty$.
Then
\[
\lim_{\omega \to \infty} \lv S_\omega(f) - f \rv_{\lprd} = 0
\qquad \text{for each $f \in \lprd$.}
\]
\end{theorem}

\vspace{6pt} It remains to prove $L^p$ boundedness of the Hilbert
transform on $\R$.

\chapter[Hilbert and Riesz transforms on $\ltrd$]{Hilbert and Riesz transforms on $\ltrd$} \label{ch:hrt}

\subsubsection*{Goal} Develop spatial and frequency representations of Hilbert and
Riesz transforms

\subsubsection*{Reference}

\cite{D} Section 4.3

\noindent \cite{G} Section 4.1

\vspace{18pt}

\begin{definition}\rm
The Riesz transforms on $\Rd$ are
\begin{align*}
R_j : \ltrd & \to \ltrd \\
f & \mapsto (-i (\xi_j/|\xi|) \fhat \, )\check{\,}
\end{align*}
for $j=1,\ldots,d$.

In dimension $d=1$, the Riesz transform equals the Hilbert
transform on $\R$, defined by
\begin{align*}
H : \ltr & \to \ltr \\
f & \mapsto (-i \sign(\xi) \fhat \, )\check{\,}
\end{align*}
because $\sign(\xi)=\xi/|\xi|$.
\end{definition}

$R_j$ is bounded since the Fourier multiplier $-i \xi_j/|\xi|$ is
a bounded function (in fact, bounded by $1$). Clearly
\begin{align*}
\lv R_j \rv_{\ltrd \to \ltrd} & \leq 1 && \text{by Plancherel,} \\
\sum_{j=1}^d R_j^2& = -I && \text{since $\sum_{j=1}^d (-i \xi_j/|\xi|)^2=-1$,} \\
R_j^* & = -R_j && \text{by Parseval.}
\end{align*}

\begin{proposition}[Spatial representation of Hilbert transform] \label{pr:htr}
If $f \in \ltr$ is $C^1$-smooth on an interval then
\begin{equation} \label{eq:hsr}
(Hf)(x) = \text{p.v.} \int_\R f(x-y) \frac{1}{\pi y} \, dy
\end{equation}
for almost every $x$ in the interval.
\end{proposition}

The proposition says formally that
\[
Hf = f * \frac{1}{\pi x}
\]
or
\[
\big( \text{p.v.}\ \frac{1}{\pi x} \big)\widehat{\ } =
-i\sign(\xi) .
\]
Later we will justify these formulas in terms of distributions.

The right side of \eqref{eq:hsr} is a \emph{singular integral},
since the convolution kernel $1/\pi y$ is not integrable.

\begin{proof}\
[This proof is similar to Proposition~\ref{pr:hst} on $\T$, and so
was skimmed only lightly in class.] For $\omega>0$,
\begin{align}
\pr \int_{[-\omega,\omega]} (-i) \sign(\xi) e^{i\xi
y} \, d\xi & = \frac{i}{2\pi} \int_{-\omega}^0 e^{i\xi y} \, d\xi
- \frac{i}{2\pi} \int_0^\omega e^{i\xi y} \, d\xi \notag \\
& = \frac{1-\cos(\omega y)}{\pi y} . \label{eq:hsr1}
\end{align}
If $f \in L^1 \cap L^2(\R)$ then
\begin{align*}
(\charfn_{[-\omega,\omega]} \widehat{Hf})\,\check{\,}(x) 
& = \pr \int_{[-\omega,\omega]} (-i) \sign(\xi)
\fhatxi e^{i\xi x} \, d\xi \\
& = \int_\R f(y) \pr \int_{[-\omega,\omega]} (-i)
\sign(\xi) e^{i\xi(x-y)} \, d\xi dy \qquad \text{by Fubini} \\
& = \int_\R f(x-y) \frac{1-\cos(\omega y)}{\pi y} \, dy \qquad
\text{by $y \mapsto x-y$ and \eqref{eq:hsr1}} \\
& = \int_{|y|<1} [f(x-y) - f(x)] \frac{1-\cos(\omega y)}{\pi y}
\, dy \\
& + \int_{|y|>1} f(x-y) \frac{1-\cos(\omega y)}{\pi y} \, dy
\end{align*}
by oddness of $\big( 1-\cos(\omega y) \big)/\pi y$. The second
integral converges to
\begin{equation} \label{eq:hr2}
\int_{|y|>1} f(x-y) \frac{1}{\pi y} \, dy
\end{equation}
as $\omega \to \infty$, by the Riemann--Lebesgue
Corollary~\ref{co:rlr}. The first integral similarly converges to
\begin{equation} \label{eq:hr1}
\int_{|y|<1} [f(x-y) - f(x)] \frac{1}{\pi y} \, dy ,
\end{equation}
assuming $f$ is $C^1$-smooth on a neighborhood of $x$ (which
ensures integrability of $y \mapsto [f(x-y) - f(x)]/\pi y$ on
$|y|<1$).

Meanwhile, $\charfn_{[-\omega,\omega]} \widehat{Hf}$ converges
to $\widehat{Hf}$ in $\ltr$ as $\omega \to \infty$, so that
$(\charfn_{[-\omega,\omega]} \widehat{Hf})\check{\,}$ converges
to $Hf$. Convergence holds a.e.\ for some subsequence of
$\omega$-values. Formula \eqref{eq:hsr} therefore follows from
\eqref{eq:hr2} and \eqref{eq:hr1}, since $\int_{\e<|y|<1} (1/\pi
y) \, dy = 0$.

Finally, one deduces \eqref{eq:hsr} in full generality by
approximating $f$ off a neighborhood of $x$ using functions in
$L^1 \cap L^2$. (Obviously $f$ belongs to $L^1 \cap L^2$ already
on each neighborhood of $x$.)
\end{proof}

\begin{proposition}[Spatial representation of Riesz transform] \label{pr:rt}
If $f \in \ltrd$ is $C^1$-smooth on an open set $U \subset \Rd$
then
\[
(R_j f)(x) = \text{p.v.} \int_\Rd f(x-y) \frac{c_d \,
y_j}{|y|^{d+1}} \, dy
\]
for almost every $x \in U$, for each $j=1,\ldots,d$.
\end{proposition}

Here $c_d = \Gamma \big( (d+1)/2 \big)/\pi^{(d+1)/2}>0$. For
example, $c_1=1/\pi$.

The proposition says formally that
\[
R_j f = f * \frac{c_d \, y_j}{|y|^{d+1}}
\]
or
\[
\big( \text{p.v.}\ \frac{c_d \, y_j}{|y|^{d+1}} \big)\widehat{\ }
= -i\frac{\xi_j}{|\xi|} .
\]

\begin{proof}
To motivate the following proof, observe
\begin{equation} \label{eq:mot0}
\frac{1}{|\xi|} = \int_0^\infty e^{-|\xi|z} \, dz
\end{equation}
and that $e^{-|\xi|z}$ is the Fourier transform of the Poisson
kernel $P_{1/z}$. Our proof will use a truncated version of this
identity:
\begin{equation} \label{eq:mot1}
\frac{e^{-|\xi|\delta} - e^{-|\xi|/\delta}}{|\xi|} =
\int_\delta^{1/\delta} e^{-|\xi|z} \, dz .
\end{equation}
In class we proceeded formally, skipping the rest of this proof and
using \eqref{eq:mot0} instead of \eqref{eq:mot1} in the proof of
Lemma~\ref{le:rtech} below.

For $f \in \ltrd$,
\begin{align*}
\widehat{(R_j f)}(\xi) & = -i \frac{\xi_j}{|\xi|} \fhatxi \\
& = \lim_{\delta \to 0} (-i \xi_j) \frac{e^{-|\xi|\delta} -
e^{-|\xi|/\delta}}{|\xi|} \fhatxi
\end{align*}
with convergence in $\ltrd$ (by dominated convergence). Applying
$L^2$ Fourier inversion yields
\[
(R_j f)(x) = \lim_{\delta \to 0} \left( -i \xi_j
\frac{e^{-|\xi|\delta} - e^{-|\xi|/\delta}}{|\xi|} \hat{f} \right)
\! \check{\ }(x)
\]
in $\ltrd$, and hence pointwise a.e.\ for some subsequence of
$\delta$ values. Thus the theorem is proved when $f \in L^1 \cap
L^2(\Rd)$, by Lemma~\ref{le:rtech} below.

Finally, one deduces the theorem for $f \in \ltrd$ by
approximating $f$ off a neighborhood of $x$ using functions in
$L^1 \cap L^2$. (Obviously $f$ belongs to $L^1 \cap L^2$ already
on each neighborhood of $x$.)
\end{proof}

\begin{lemma} \label{le:rtech}
If $f \in L^1 \cap L^2(\Rd)$ is $C^1$-smooth on an open set $U
\subset \Rd$, then
\begin{align}
& \lim_{\delta \to 0} \left( -i \xi_j \frac{e^{-|\xi|\delta} -
e^{-|\xi|/\delta}}{|\xi|} \hat{f} \right) \! \check{\ }(x) \notag \\
& = \int_{|y|<1} [f(x-y)-f(x)] \frac{c_d \, y_j}{|y|^{d+1}} \, dy
+ \int_{|y|>1} f(x-y) \frac{c_d \, y_j}{|y|^{d+1}} \, dy
\label{eq:tech1}
\end{align}
for almost every $x \in U$. Further, the first integral in
\eqref{eq:tech1} equals
\[
\lim_{\e \to 0} \int_{\e<|y|<1} f(x-y) \frac{c_d \,
y_j}{|y|^{d+1}} \, dy .
\]
\end{lemma}

\begin{proof}
First, $\xi_j/|\xi|$ is bounded by $1$, and the exponentials
$e^{-|\xi|\delta}$ and $e^{-|\xi|/\delta}$ are square integrable,
and so is $\fhat$. Thus their product is integrable, so that by the
$L^1$  Fourier Inversion Theorem~\ref{th:fi} (and the definition of
$\fhat$ for $f \in \lord$),
\begin{align*}
& \left( -i \xi_j \frac{e^{-|\xi|\delta} - e^{-|\xi|/\delta}}{|\xi|}
\hat{f} \right) \! \check{\ }(x) \\
& = - \frac{1}{(2\pi)^d} \int_\Rd
i \xi_j \frac{e^{-|\xi|\delta} -
e^{-|\xi|/\delta}}{|\xi|} \int_\Rd f(y) e^{-i\xi y} \, dy \, e^{i\xi x} \, d\xi \\
& = - \int_\Rd f(x-y) \, \frac{1}{(2\pi)^d} \int_\Rd i \xi_j
\frac{e^{-|\xi|\delta} - e^{-|\xi|/\delta}}{|\xi|} e^{i\xi y} \,
d\xi dy
\end{align*}
after changing $y \mapsto x-y$. To evaluate the inner integral, we
express it using Poisson kernels:
\begin{align*}
& \frac{1}{(2\pi)^d} \int_\Rd i \xi_j \frac{e^{-|\xi|\delta} -
e^{-|\xi|/\delta}}{|\xi|} e^{i\xi y} \, d\xi \\
& = \frac{\partial\ }{\partial y_j} \frac{1}{(2\pi)^d} \int_\Rd
\frac{e^{-|\xi|\delta} - e^{-|\xi|/\delta}}{|\xi|} e^{i\xi y}
\, d\xi \\
& = \int_\delta^{1/\delta} \frac{\partial\ }{\partial y_j}
\frac{1}{(2\pi)^d} \int_\Rd e^{-|\xi|z} e^{i\xi y} \,
d\xi dz && \text{by identity \eqref{eq:mot1}} \\
& = \int_\delta^{1/\delta} \frac{\partial\ }{\partial y_j} \, P_{1/z}(y) \, dz && \text{by \eqref{eq:PR3}} \\
& = \int_\delta^{1/\delta} c_d \, z \, \frac{\partial\ }{\partial y_j} \, \frac{1}{(|y|^2+z^2)^{(d+1)/2}} \, dz && \text{by \eqref{eq:PR3b}} \\
& = \int_\delta^{1/\delta} c_d \, y_j \, \frac{\partial\ }{\partial z} \, \frac{1}{(|y|^2+z^2)^{(d+1)/2}} \, dz && \text{(why?!)} \\
& = \left. \frac{c_d y_j}{(|y|^2+z^2)^{(d+1)/2}}
\right|_{z=\delta}^{z=1/\delta} .
\end{align*}
By substituting this expression into the above, we find
\begin{align}
& \left( -i \xi_j \frac{e^{-|\xi|\delta} -
e^{-|\xi|/\delta}}{|\xi|} \hat{f} \right) \! \check{\ }(x) \notag \\
& = - \int_\Rd f(x-y) \left. \frac{c_d \,
y_j}{(|y|^2+z^2)^{(d+1)/2}}
\right|_{z=\delta}^{z=1/\delta} \, dy \notag \\
& = - \int_{|y|<1} [f(x-y) - f(x)] \left. \frac{c_d \,
y_j}{(|y|^2+z^2)^{(d+1)/2}}
\right|_{z=\delta}^{z=1/\delta} \, dy \label{eq:tech2} \\
& \quad - \int_{|y|>1} f(x-y) \left. \frac{c_d \,
y_j}{(|y|^2+z^2)^{(d+1)/2}} \right|_{z=\delta}^{z=1/\delta} \, dy
\label{eq:tech3}
\end{align}
where we used the oddness of $y_j$ to insert $f(x)$ in
\eqref{eq:tech2}.

Now fix a point $x \in U$. As $\delta \to 0$, expression
\eqref{eq:tech2} converges to
\[
\int_{|y|<1} [f(x-y) - f(x)] \frac{c_d \, y_j}{|y|^{d+1}} \, dy
\]
by dominated convergence (noting the $C^1$-smoothness ensures the
integrand is $O(|y|) \cdot O(1/|y|^d) = O(1/|y|^{d-1})$ near the
origin, which is integrable). And as $\delta \to 0$, expression
\eqref{eq:tech3} converges to
\[
\int_{|y|>1} f(x-y) \frac{c_d \, y_j}{|y|^{d+1}} \, dy
\]
by dominated convergence (noting $f \in \ltrd$ and
$y_j/|y|^{d+1}=O(1/|y|^d)$ is square integrable for $|y|>1$).
(Exercise: explain why the terms with $z=1/\delta$ in
\eqref{eq:tech2} and \eqref{eq:tech3} vanish as $\delta \to 0$,
using dominated convergence.)

For the final claim in the lemma, write $\int_{|y|<1}=\lim_{\e \to
0} \int_{\e<|y|<1}$ and use the oddness of $y_j$ to remove the
term with $f(x)$.
\end{proof}

\subsubsection*{Connections to PDEs}

1. The Riesz transforms map the \emph{normal} derivative of a
harmonic function to its \emph{tangential} derivatives.

\vspace{6pt} \noindent \emph{Formal Proof.} Given a function $f$,
let
\[
u(x,x_{d+1})=(P_{1/x_{d+1}} * f)(x) , \qquad x \in \Rd, \quad
x_{d+1} > 0,
\]
so that $u$ is harmonic on the upper halfspace $\Rd \times
(0,\infty)$ with boundary value $u=f$ when $x_{d+1}=0$ (see
Chapter~\ref{ch:fis}). Put
\[
v(x) = \left. \frac{\partial \quad \ }{\partial x_{d+1}}
u(x,x_{d+1}) \right|_{x_{d+1}=0} = \text{normal derivative of $u$ at
the boundary.}
\]
Then
\[
R_j v = \frac{\partial f}{\partial x_j} , \qquad j=1,\ldots,d ,
\]
because
\begin{align*}
\widehat{(R_j v)}(\xi) & = -i \frac{\xi_j}{|\xi|} \widehat{v}(\xi) \\
& =  -i \frac{\xi_j}{|\xi|} \left. \frac{\partial \quad \ }{\partial
x_{d+1}} \widehat{u}(\xi,x_{d+1}) \right|_{x_{d+1}=0} \\
& =  -i \frac{\xi_j}{|\xi|} \left. \frac{\partial \quad \ }{\partial
x_{d+1}} \big( e^{-|\xi|x_{d+1}} \fhatxi \big) \right|_{x_{d+1}=0} \\
& =  -i \frac{\xi_j}{|\xi|} (-|\xi|) \fhatxi \\
& =  i \xi_j \fhatxi \\
& = \left( \frac{\partial f}{\partial x_j} \right)\widehat{\ }(\xi)
.
\end{align*}
Thus we have shown the $j$th Riesz transform maps the normal
derivative of $u$ to its $j$th tangential derivative, on the
boundary.

\vspace{6pt} \noindent 2. Mixed Riesz transforms map the Laplacian
to mixed partial derivatives.

\vspace{6pt} \noindent \emph{Formal Proof.}
\[
\left( \frac{\partial^2 f}{\partial x_j^2} \right)\widehat{\ }(\xi)
= (i\xi_j)^2 \fhatxi = -\xi_j^2 \fhatxi
\]
and so summing over $j$ gives
\[
(\Delta f)\widehat{\ }(\xi) = -|\xi|^2 \fhatxi .
\]
Hence
\begin{align*}
(R_j R_k \Delta f)\widehat{\ }(\xi) & = \frac{(-i \xi_j)}{|\xi|}
\frac{(-i \xi_k)}{|\xi|} (-|\xi|^2) \fhatxi \\
& = -(i\xi_j)(i\xi_k) \fhatxi \\
& = - \left( \frac{\partial^2 f}{\partial x_j \partial x_k}
\right)\widehat{\ }(\xi)
\end{align*}
so that
\[
R_j R_k \Delta f = - \frac{\partial^2 f}{\partial x_j \partial x_k}
.
\]
That is, mixed Riesz transforms map the Laplacian to mixed partial
derivatives.

The above formal derivation is rigorous if, for example, $f$ is
$C^2$-smooth with compact support.

\vspace{6pt} Consequently, the norm of a mixed second derivative is
controlled by the norms of the pure second derivatives in the
Laplacian, with
\[
\left\lv \frac{\partial^2 f}{\partial x_j \partial x_k}
\right\rv_\ltrd \leq \lv \Delta f \rv_\ltrd
\]
since each Riesz transform has norm $1$ on $\ltrd$. Similar
estimates hold on $\lprd, 1<p<\infty$, by the $L^p$ boundedness of
the Riesz transform proved in the next chapter.

\chapter[Hilbert and Riesz transforms on $\lprd$]{Hilbert and Riesz transforms on $\lprd$} \label{ch:hpr}

\subsubsection*{Goal} Prove weak $(1,1)$ for Riesz transform, and deduce strong $(p,p)$ by interpolation and
duality

\subsubsection*{Reference}

\cite{D} Section 5.1

\vspace{18pt}

\begin{theorem}[weak $(1,1)$ on $L^1 \cap L^2(\Rd)$]
\label{th:rw11} There exists $A>0$ such that
\[
| \{x \in \Rd : |(R_j f)(x)| > \omega \} | \leq \frac{A}{\omega}
\lv f\rv_\lord
\]
for all $\omega>0, j=1,\ldots,d$ and $f \in L^1 \cap L^2(\Rd)$.
\end{theorem}

\begin{proof}
Apply the Calder\'{o}n--Zygmund Theorem~\ref{th:cz} to get $f=g+b$.
Note $g \in L^1 \cap L^\infty(\Rd)$ and so $g \in \ltrd$, hence $R_j
g \in \ltrd$ by Chapter~\ref{ch:hrt}. And $b=f-g \in \ltrd$ so that
$R_j b \in \ltrd$.

Now proceed like in the proof of Theorem~\ref{th:hw11}, just
changing $\T$ to $\Rd$ and the interval $I(l)$ to the cube $Q(l)$.
To finish the proof, we want to show
\begin{equation} \label{eq:rlp1}
\sum_l \int_{\Rd \setminus 2\sqrt{d}Q(l)} |(R_j b_l)(x)| \, dx \leq
(\text{const.}) \lv f \rv_\lord .
\end{equation}
By Proposition~\ref{pr:rt} applied on the open set $U=\Rd \setminus
2\sqrt{d}\overline{Q(l)}$ (where $b_l=0$), we have
\begin{align*}
& \int_{\Rd \setminus 2\sqrt{d}Q(l)} |R_j b_l(x)| \, dx \\
& = \int_{\Rd \setminus 2\sqrt{d}Q(l)} \Big| \int_{Q(l)} b_l(y) \frac{c_d (x_j - y_j)}{|x-y|^{d+1}} \, dy \Big| \, dx \\
& \qquad \qquad \text{noting $x-y$ is bounded away from $0$,
since $y \in Q(l)$ and $x \notin 2\sqrt{d}Q(l)$,} \\
& = \int_{\Rd \setminus 2\sqrt{d}Q(l)} \big| \int_{Q(l)} b_l(y) \big[ \rho_j(x-y) - \rho_j(x-c(l)) \big] \, dy \big| \, dx \\
\end{align*}
where
\[
\rho_j(x) = c_d \frac{x_j}{|x|^{d+1}}
\]
is the $j$th Riesz kernel and $c(l)$ is the center of $Q(l)$; here
we used that $\int_{Q(l)} b_l(y) \, dy = 0$. Hence
\begin{align}
& \int_{\Rd \setminus 2\sqrt{d}Q(l)} |R_j b_l(x)| \, dx \notag \\
& \leq \int_{Q(l)} |b_l(y)| \int_{\Rd \setminus 2\sqrt{d}Q(l)} |\rho_j(x-y) - \rho_j(x-c(l))| \, dx dy \notag \\
& \leq (\text{const.}) \int_{Q(l)} |b_l(y)| \, dy \label{eq:rlp2}
\end{align}
by Lemma~\ref{le:hc} below; the hypotheses of that lemma are
satisfied here because
\[
|(\nabla \rho_j)(x)| \leq  \frac{(\text{const.})}{|x|^{d+1}}
\]
and if $x \in \Rd \setminus 2\sqrt{d}Q(l)$ and $y \in Q(l)$ then
\begin{align*}
|x-c(l)| & \geq \frac{1}{2} \, \text{side}\big( 2\sqrt{d}Q(l) \big) \\
& \geq 2|y-c(l)| .
\end{align*}
Now \eqref{eq:rlp1} follows by summing \eqref{eq:rlp2} over $l$ and
recalling that $\lv b \rv_\lord \leq 2 \lv f \rv_\lord$ by the
Calder\'{o}n--Zygmund Theorem~\ref{th:cz}.
\end{proof}

\begin{lemma}[H\"{o}rmander condition] \label{le:hc}
If $\rho \in C^1(\Rd \setminus \{ 0 \})$ with
\[
|(\nabla \rho)(x)| \leq  \frac{(\text{const.})}{|x|^{d+1}} , \qquad
x \in \Rd ,
\]
then
\[
\sup_{y,z \in \Rd} \int_{\{ x : |x-z| \geq 2|y-z| \}} |\rho(x-y) -
\rho(x-z)| \, dx < \infty .
\]
\end{lemma}

\begin{proof}
We can take $z=0$, by a translation. By the Fundamental Theorem,
\begin{align*}
\rho(x-y) - \rho(x) & = \int_0^1 \frac{\partial \ }{\partial s}
\rho(x-sy) \, ds \\
& = -\int_0^1 y \cdot (\nabla \rho)(x-sy) \, ds .
\end{align*}
Hence
\begin{align*}
& \int_{\{ x : |x| \geq 2|y| \}} |\rho(x-y) - \rho(x)| \, dx
\\
& \leq |y| \int_0^1 \int_{|x| \geq 2|y|} |(\nabla \rho)(x-sy)| \, dx ds \\
& \leq (\text{const.}) |y| \int_{|x| \geq 2|y|}
\frac{1}{(|x|/2)^{d+1}} \, dx \\
& \qquad \qquad \qquad \text{by using the hypothesis, since $|x-sy| \geq |x|-|y| \geq |x|/2$,} \\
& = (\text{const.}) |y| \int_{2|y|}^\infty \frac{1}{r^{d+1}} r^{d-1}
\, dr \\
& = (\text{const.})
\end{align*}
\end{proof}

Now we deduce strong $(p,p)$ estimates.

\begin{corollary} \label{co:rspp}
The Riesz transforms are strong $(p,p)$ for $1<p<\infty$.
\end{corollary}

\begin{proof}
$R_j$ is strong $(2,2)$ and linear, by definition in
Chapter~\ref{ch:hrt}, and $R_j$ is weak $(1,1)$ on $L^1 \cap
L^2(\Rd)$ (and hence on all simple functions with support of
finite measure) by Theorem~\ref{th:rw11}. So $R_j$ is strong
$(p,p)$ for $1<p<2$ by Remark~\ref{re:mi} after Marcinkiewicz
Interpolation (in Appendix~\ref{ap:ip}). That is, $R_j: \lprd \to
\lprd$ is bounded and linear for $1<p<2$.

For $2<p<\infty$ we use duality and anti-selfadjointness
$R_j^*=-R_j$ on $\ltrd$ to reduce to the case $1<p<2$, just like in
the proof of Corollary~\ref{co:hspp}.
\end{proof}

Alternatively, for singular integral kernels of the form
\[
\frac{{\mathcal O}(x/|x|)}{|x|^d}
\]
where ${\mathcal O}$ is an odd function on the unit sphere, one can
instead use the \emph{method of rotations} \cite[Section 4.2c]{G}.
The idea is to express convolution with this kernel as an average of
Hilbert transforms taken in all possible directions in $\Rd$.

The Riesz kernel $c_d(x_j/|x|)/|x|^d$ fits this form, since
${\mathcal O}(y)=y_j$ is odd.

The strong $(p,p)$ bound on the Riesz transform can be generalized
to a whole class of convolution-type singular integral operators
\cite[Section 5.1]{D}.

\part{Fourier series and integrals}

\chapter[Band limited functions]{Compactly supported Fourier transforms, and the sampling theorem} \label{ch:bl}

\subsubsection*{Goal}

Show band limited functions are holomorphic

\noindent Prove the Kotelnikov--Shannon--Whittaker sampling theorem

\subsubsection*{Reference}

\cite{K} Section VI.7

\vspace{18pt}

\begin{definition}
We say $f=\check{g}$ is \textbf{band limited} if $g \in \lord$ has
compact support.
\end{definition}

\begin{theorem}[Band limited functions are holomorphic] \label{th:blh}
Assume $g \in \lord$ is supported in a ball $B(R)$, and define
\[
f(z) = \check{g}(z) = \frac{1}{(2\pi)^d} \int_{B(R)} g(\xi) e^{i\xi
z} \, d\xi
\]
for $z=x+iy \in \Cd, x,y \in \Rd$. (Here $\xi z = \xi_1 z_1 + \cdots
+ \xi_d z_d$.)

Then $f$ is holomorphic, and $|f(z)| = O(e^{R|y|})$.

If in addition $g \in \ltrd$ then $|f(z)| = O(e^{R|y|}/\sqrt{|y|})$.
\end{theorem}

Thus once more, decay of the Fourier transform (here, compact
support) implies smoothness of the function (here, holomorphicity).
The theorem also bounds the rate of growth of the function in the
complex directions. (The function must vanish at infinity in the
real directions, by the Riemann--Lebesgue corollary, since $g$ is
integrable.)

For example, the Dirichlet kernel $D(x)=\sin(x)/\pi x =
(\charfn_{[-1,1]})\check{\ }(x)$ is band limited with $R=1$, in $1$
dimension. Taking $z=0+iy$, we calculate
\[
D(iy) = \frac{e^y-e^{-y}}{2\pi y} = O(e^{|y|}/|y|) ,
\]
which is better (by a factor of $\sqrt{|y|}$) than is guaranteed by
the theorem.

\begin{proof}
$f$ is well defined because $\xi \mapsto e^{i\xi z}$ is bounded on
$B(R)$, for each $z$. And $f$ is holomorphic because $e^{i\xi  z}$
is holomorphic and $f$ can be differentiated through the integral
with respect to the complex variable $z$. (Exercise. Justify these
claims in detail.)

Clearly
\begin{align*}
|f(z)| & \leq \frac{1}{(2\pi)^d} \int_{B(R)} |g(\xi)| e^{-\xi y} \,
d\xi && \text{since $e^{i\xi z} = e^{i\xi x} e^{-\xi y}$} \\
& \leq \frac{1}{(2\pi)^d} \lv g \rv_\lord e^{R|y|} .
\end{align*}
If in addition $g \in \ltrd$, then
\[
|f(z)| \leq \frac{1}{(2\pi)^d} \lv g\rv_\ltrd \big( \int_{B(R)}
e^{-2\xi y} \, d\xi \big)^{\! 1/2}
\]
and
\begin{align*}
\int_{B(R)} e^{-2\xi y} \, d\xi
& = \int_{B(R)} e^{-2\xi |y| e_1} \, d\xi \\
& \qquad \qquad \text{by $\xi \mapsto \xi A$ for some orthogonal matrix $A$ with $Ay=|y|e_1$} \\
& = \int_{B(R)} e^{-2\xi_1 |y|} \, d\xi \\
& \leq \int_{[-R,R]^d} e^{-2\xi_1 |y|} \, d\xi \\
& = (2R)^{d-1} \frac{e^{2R|y|}-e^{-2R|y|}}{2|y|} \\
& \leq \frac{(2R)^{d-1}}{2} \frac{e^{2R|y|}}{|y|} .
\end{align*}
Hence $|f(z)| \leq (\text{const.})e^{R|y|}/\sqrt{|y|}$.
\end{proof}

Holomorphic functions are known to be determined by their values
on lower dimensional sets in $\Cd$. For a band limited function,
that ``sampling set'' can be a lattice in $\Rd$.

\begin{theorem}[Sampling theorem for band limited functions]
\label{th:ss} Assume $f \in \ltrd$ is band limited, with $\fhat$
supported in the cube $[-\omega,\omega]^d$ for some $\omega>0$.

Then
\[
f(x) = \sum_{n \in \Zd} f \big( \frac{\pi}{\omega} n \big)
\prod_{j=1}^d \sinc(\omega x_j - \pi n_j)
\]
with the series converging in $\ltrd$, and also uniformly (in
$\lird$).
\end{theorem}

\begin{remark}\label{re:ss}\rm \

\noindent 1. The sampling rate $\omega/\pi$ is proportional to the
bandwidth $\omega$, that is, to the highest frequency contained in
the signal $f$. Intuitively, the sampling rate must be high when the
frequencies are high, because many samples are needed to determine a
highly oscillatory function.

\noindent 2. $\sinc(\omega x_j - \pi n_j)$ is centered at the
sampling location $(\pi/\omega)n_j$ and rescaled to have bandwidth
$\omega$. It vanishes at all the other sampling locations
$(\pi/\omega)m_j$, since
\[
\sinc \big( \omega(\pi/\omega)m_j - \pi n_j \big) = \sinc
\big(\pi(m_j-n_j) \big) = 0.
\]

\noindent 3. A graphical example of the sampling formula is shown in
Figure~\ref{samplefig}, for $\omega=2\pi$ and
\begin{align*}
f(x) & = -\sinc \big(2\pi (x + 1) \big) + 2 \sinc \big( 2\pi (x +
1/2) \big) + 3 \sinc \big( 2\pi x \big) \\
& \qquad \qquad \qquad + 2 \sinc \big( 2\pi (x - 1/2) \big) + 1
\sinc \big( 2\pi (x - 1) \big) .
\end{align*}
The figure shows $f$ with a solid curve, and $3 \sinc \big( 2\pi x
\big)$ and $2 \sinc \big( 2\pi (x - 1/2) \big)$ with dashed curves.
\begin{figure}
\begin{center}
  \includegraphics[scale=0.8]{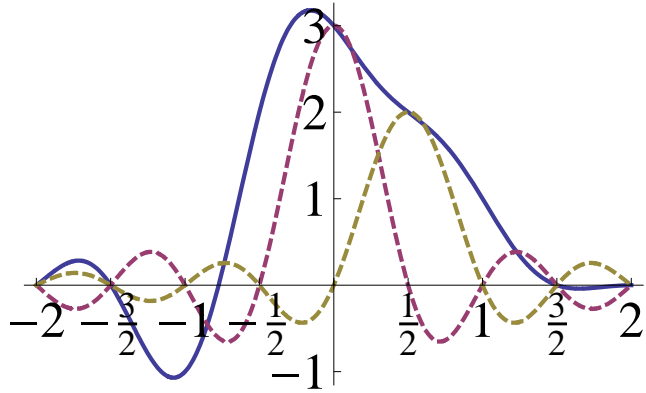}
  \caption{\label{samplefig}
    Example of sampling formula in Theorem~\ref{th:ss}, with $\omega=2\pi$ and sampling rate $\omega/\pi=2$.
    The dashed curves are two of the $\sinc$ functions making up the signal. See Remark~\ref{re:ss}.}
\end{center}
\end{figure}
\end{remark}

\begin{proof}[Proof of Sampling Theorem]
We can assume $\omega=\pi$, by replacing $x$ with $(\pi/\omega)x$
(Exercise).

Next, $\fhat$ is square integrable and compactly supported, and so
is integrable. Hence by $L^1$ Fourier inversion, $f$ is continuous
(after redefining it on some set of measure zero) with
\begin{equation} \label{eq:ss2}
f(x) = \frac{1}{(2\pi)^d} \int_\Rd \fhatxi e^{i\xi x} \, d\xi ,
\qquad x \in \Rd .
\end{equation}
Thus the pointwise sampled values $f\big( (\pi/\omega)n \big)$ in
the theorem are well defined.

We will prove
\begin{equation} \label{eq:ss1}
\fhatxi = \sum_{n \in \Zd} f(-n) e^{i\xi n} , \qquad \xi \in
[-\pi,\pi]^d ,
\end{equation}
with convergence in $\lt([-\pi,\pi]^d)$.  Indeed, if we regard
$\fhat$ as a square integrable function on the cube
$\Td=[-\pi,\pi]^d$, then its Fourier coefficients are
\begin{align*}
& \frac{1}{(2\pi)^d} \int_{[-\pi,\pi]^d} \fhatxi e^{-i\xi n} \,
d\xi
\\
& = \frac{1}{(2\pi)^d} \int_\Rd \fhatxi e^{-i\xi n} \, d\xi
&& \text{since $\fhat$ is supported in $[-\pi,\pi]^d$} \\
& = f(-n)
\end{align*}
by the inversion formula \eqref{eq:ss2}. Thus \eqref{eq:ss1}
simply expresses the Fourier series of $\fhat$ on the cube.

After changing $n \mapsto -n$ in \eqref{eq:ss1}, we have
\[
\fhatxi = \sum_{n \in \Zd} f(n) e^{-i\xi n}
\charfn_{[-\pi,\pi]^d}(\xi) , \qquad \xi \in \Rd ,
\]
with convergence in $\ltrd$ and in $\lord$. Applying $L^2$
inversion gives
\begin{align*}
f(x) & = \sum_{n \in \Zd} f(n) \big( e^{-i\xi n}
\charfn_{[-\pi,\pi]^d} \big)\check{\ }(x) \\
& = \sum_{n \in \Zd} f(n) \prod_{j=1}^d
\frac{\sin(\pi(x_j-n_j))}{\pi(x_j-n_j)}
\end{align*}
with convergence in $\ltrd$. Applying $L^1$ inversion gives
convergence in $L^\infty$.
\end{proof}

\subsubsection*{Paley--Wiener space}

For a deeper perspective on Sampling Theorem~\ref{th:ss}, consider
the \emph{Paley--Wiener space}
\[
PW(\omega) = \{ f \in \ltrd : \text{$\widehat{f}$ is supported in
$[-\omega,\omega]^d$ } \} .
\]
Clearly $PW(\omega)$ is a subspace of $\ltrd$, and it is a closed
subspace (since if $f = \lim_m f_m$ in $\ltrd$ and $\widehat{f_m}$
is supported in $[-\omega,\omega]^d$, then $\fhat = \lim_m
\widehat{f_m}$ is also supported in $[-\omega,\omega]^d$).

Hence $PW(\omega)$ is a Hilbert space with the $\lt$ inner product.
It is isometric, under the Fourier transform, to
$\lt([-\omega,\omega]^d)$ with inner product $(2\pi)^{-d} \langle
\cdot , \cdot \rangle_\lt$. That space has orthonormal Fourier basis
\[
\left\{ (\pi/\omega)^{d/2} \charfn_{[-\omega,\omega]^d}(\xi)
e^{-i\xi(\pi/\omega)n} \right\}_{n \in \Zd} ,
\]
where the indicator function simply reminds us that we are working
on the cube. Taking the inverse Fourier transform gives an
orthonormal basis of $\sinc$ functions for the Paley--Wiener space:
\[
\{ g_n \}_{n \in \Zd} = \big\{ (\omega/\pi)^{d/2} \prod_{j=1}^d
\sinc(\omega x_j - \pi n_j) \big\}_{n \in \Zd} .
\]
Using this orthonormal basis, we expand
\begin{equation} \label{eq:pw1}
f = \sum_{n \in \Zd} \la f , g_n \ra_\lt \, g_n , \qquad \text{for
all $f \in PW(\omega)$,}
\end{equation}
where the coefficient is
\begin{align*}
\la f , g_n \ra_\lt & = \frac{1}{(2\pi)^d} \la \fhat ,
(\pi/\omega)^{d/2} \charfn_{[-\omega,\omega]^d} e^{-i\xi
(\pi/\omega)n} \ra_\lt && \text{by
Parseval} \\
& = (\pi/\omega)^{d/2} f\big( (\pi/\omega)n \big)
\end{align*}
by Fourier inversion. Thus the orthonormal expansion \eqref{eq:pw1}
simply restates the Sampling Theorem~\ref{th:ss}.

Our calculations have, of course, essentially repeated the proof of
the Sampling Theorem.

\chapter{Periodization and Poisson summation} \label{ch:psf}

\subsubsection*{Goal}

Periodize functions on $\Rd$ to functions on $\Td$

\noindent Show the Fourier series of periodization gives the Poisson
summation formula

\subsubsection*{References}

\cite{F} Section 8.3

\noindent \cite{K} Section VI.1

\vspace{18pt}

\begin{definition}
Given $f \in \lord$, its \textbf{periodization} is the function
\[
\Pe(f)(x)=(2\pi)^d \sum_{n \in \Zd} f(x+2\pi n) , \qquad x \in \Rd .
\]
\end{definition}

\begin{example}
In $1$ dimension, if $f=\charfn_{[-\pi,2\pi)}$, then
$\Pe(f)=2\pi(2\charfn_{[-\pi,0)}+\charfn_{[0,\pi)})$ for $x \in
[-\pi,\pi)$, with $\Pe(f)$ extending $2\pi$-periodically to $\R$.
\end{example}

\begin{lemma} \label{le:per}
If $f \in \lord$ then the series for $\Pe(f)(x)$ converges
absolutely for almost every $x$, and $\Pe(f)$ is $2\pi
\Zd$-periodic. Further, $\Pe : \lord \to \lotd$ is bounded, with
\[
\lv \Pe(f) \rv_\lotd \leq \lv f \rv_\lord .
\]
The periodization has Fourier coefficients
\[
\widehat{\Pe(f)}(j) = \fhatj , \qquad j \in \Zd .
\]
\end{lemma}

That is, the $j$th Fourier coefficient of $\Pe(f)$ equals the
Fourier transform of $f$ at $j$.

\begin{proof}
See Problem~\ref{pr:pef} in Assignment 3.
\end{proof}

\begin{lemma}[Periodization of a convolution] \label{le:pcon}
If $f,g \in \lord$ then
\[
\Pe(f*g) = \Pe(f) * \Pe(g) .
\]
\end{lemma}

\begin{proof}
We have
\begin{align*}
\big( \Pe(f*g) \big)\widehat{\ }(j) & = (f*g)\widehat{\ }(j)
&& \text{by Lemma~\ref{le:per}} \\
& = \fhat(j) \, \ghat(j) \\
& = \widehat{\Pe(f)}(j) \, \widehat{\Pe(g)}(j) && \text{by Lemma~\ref{le:per} again} \\
& = \big( \Pe(f) * \Pe(g) \big)\widehat{\ }(j)
\end{align*}
and so $\Pe(f*g) = \Pe(f) * \Pe(g)$ by the uniqueness theorem for
Fourier series.

For a more direct proof, suppose $f$ and $g$ are bounded with
compact support, so that the sums in the following argument are all
finite rather than infinite. (Thus sums and integrals can be
interchanged, below.)

For each $x \in \Rd$,
\begin{align*}
& \Pe(f*g)(x) \\
& = (2\pi)^d \sum_{n \in \Zd} (f*g)(x+2\pi n) \\
& = (2\pi)^d \sum_{n \in \Zd} \int_\Rd f(x+2\pi n-y) g(y) \, dy \\
& = \int_\Rd \Pe(f)(x-y) g(y) \, dy && \text{by definition of $\Pe(f)$} \\
& = \sum_{m \in \Zd} \int_\Td \Pe(f)(x-y-2\pi m) g(y+2\pi m) \, dy && \text{since $\Rd = \bigcup_m (\Td + 2\pi m)$} \\
& = \sum_{m \in \Zd} \int_\Td \Pe(f)(x-y) g(y+2\pi m) \, dy && \text{using periodicity of $\Pe(f)$} \\
& = \frac{1}{(2\pi)^d} \int_\Td \Pe(f)(x-y) \Pe(g)(y) \, dy \\
& = \big( \Pe(f) * \Pe(g) \big)(x) ,
\end{align*}
remembering that our definition of convolution on $\Td$ has a
prefactor of $(2\pi)^{-d}$.

Finally, pass to the general case by a limiting argument, using that
if $f_m \to f$ in $\lord$ then $\Pe(f_m) \to \Pe(f)$ in $\lotd$ by
Lemma~\ref{le:per}.
\end{proof}

\begin{theorem}[Poisson summation formula] \label{th:psf}
Suppose $f \in \lord$ is continuous and decays in space and
frequency according to:
\begin{align}
|f(x)| & \leq \frac{C}{(1+|x|)^{d+\e}} , && x \in \Rd ,
\label{eq:psf1}
\\
|\fhatxi| & \leq \frac{C}{(1+|\xi|)^{d+\e}} , && \xi \in \Rd ,
\label{eq:psf2}
\end{align}
for some constants $C,\e>0$.

Then the periodization $\Pe(f)$ equals its Fourier series at every
point:
\[
(2\pi)^d \sum_{n \in \Zd} f(x+2\pi n) = \sum_{j \in \Zd} \fhatj
e^{ijx} , \qquad x \in \Rd .
\]
In particular, taking $x=0$ gives
\[
(2\pi)^d \sum_{n \in \Zd} f(2\pi n) = \sum_{j \in \Zd} \fhatj .
\]
\end{theorem}
This Poisson summation formula relates a lattice sum of values of
the function to a lattice sum of values of its Fourier transform.
\begin{proof}
$\Pe(f)$ has Fourier coefficients in $\ell^1(\Zd)$, since
\begin{align*}
\sum_{j \in \Zd} \big| \widehat{\Pe(f)}(j) \big| & = \sum_{j \in
\Zd} |\fhatj| && \text{by Lemma~\ref{le:per}} \\
& \leq \sum_{j \in \Zd} \frac{C}{(1+|j|)^{d+\e}} && \text{by \eqref{eq:psf2}} \\
& \leq \int_\Rd \frac{(\text{const.})}{(1+|\xi|)^{d+\e}} \, d\xi \\
& < \infty
\end{align*}
by spherical coordinates.

Hence the Fourier series of $\Pe(f)$ converges absolutely and
uniformly to a continuous function. That continuous function has the
same Fourier coefficients as $\Pe(f)$, and so it equals $\Pe(f)$
a.e.\ (just like in $1$ dimension; see Chapter~\ref{ch:l1}).

To complete the proof we will show $\Pe(f)$ is continuous, for then
$\Pe(f)$ equals its Fourier series everywhere (and not just almost
everywhere).

Notice that $\Pe(f)(x)=(2\pi)^d \sum_{n \in \Zd} f(x+2\pi n)$ is a
series of continuous functions. The series converges absolutely and
uniformly on each ball in $\Rd$ (by using \eqref{eq:psf1};
exercise), and so $\Pe(f)$ is continuous.
\end{proof}

\begin{example}[Periodizing the Poisson kernel]
The Poisson kernel $P_r$ on $\T$ equals the periodization of the
Poisson kernel $P_\omega$ on $\R$:
\begin{equation} \label{eq:pp}
\frac{1-r^2}{1-2r \cos x + r^2} = 2\pi \sum_{n \in \Z} \frac{1}{\pi}
\frac{\omega^{-1}}{(x+2\pi n)^2+\omega^{-2}} , \qquad x \in \R ,
\end{equation}
provided $r=e^{-1/\omega}$. Hence we obtain a series expansion for
the square of the cosecant:
\[
\frac{\pi^2}{\sin^2 \pi x} = \sum_{n \in \Z} \frac{1}{(x+n)^2} ,
\qquad x \in \R \setminus \Z .
\]
\noindent \Proof First, to partially motivate these results we note
$\Pe(P_\omega * f) = \Pe(P_\omega) * \Pe(f)$ by
Lemma~\ref{le:pcon}, so that it is plausible $P_\omega$ periodizes
to $P_r$ for some $r$.

To prove \eqref{eq:pp}, observe that $P_\omega$ satisfies decay
hypotheses \eqref{eq:psf1} and \eqref{eq:psf2} because
\begin{align*}
P_\omega(x) & = \frac{1}{\pi}
\frac{\omega^{-1}}{x^2+\omega^{-2}} , && x \in \R , \\
\widehat{P_\omega}(\xi) & = e^{-|\xi|/\omega} , && \xi \in \R ,
\end{align*}
by \eqref{eq:PR3b} and Table~\ref{ta:ftkernels}. Hence the Poisson
Summation Formula says that
\begin{align*}
\Pe(P_\omega)(x) & = \sum_{j \in \Z} \widehat{P_\omega}(j) e^{ijx} \\
& = \sum_{j \in \Z} e^{-|j|/\omega} e^{ijx} \\
& = \sum_{j \in \Z} r^{|j|} e^{ijx} && \text{since $r=e^{-1/\omega}$} \\
& = P_r(x)
\end{align*}
by \eqref{eq:P2}, which proves \eqref{eq:pp}.

Changing $x$ to $2\pi x$ in \eqref{eq:pp} gives
\[
\sum_{n \in \Z} \frac{1}{(x+n)^2+(2\pi \omega)^{-2}} = 2\pi^2
\omega \frac{1-r^2}{1-2r \cos(2\pi x)+r^2} .
\]
Since
\[
r=e^{-1/\omega}=1-\frac{1}{\omega}+O \big( \frac{1}{\omega^2}
\big) ,
\]
letting $\omega \to \infty$ implies that
\[
\sum_{n \in \Z} \frac{1}{(x+n)^2} = \frac{4\pi^2}{2-2\cos(2\pi x)} =
\frac{\pi^2}{(\sin \pi x)^2} ,
\]
where we used monotone convergence on the left side.
\end{example}

\begin{example}[Periodizing the Gauss kernel] \label{ex:pg}
The Gauss kernel $G_s(t)=\sum_{j \in \Z} e^{-j^2 s} e^{ijt}$ on $\T$
equals the periodization of the Gauss kernel $G_\omega$ on $\R$:
\begin{equation} \label{eq:pg}
\sum_{j \in \Z} e^{-j^2 s} e^{ijx} = 2\pi \sum_{n \in \Z}
\frac{\omega}{\sqrt{2\pi}} e^{-\omega^2 (x+2\pi n)^2/2} , \qquad x
\in \R ,
\end{equation}
provided $s>0$ and $\omega=1/\sqrt{2s}$. Hence
\[
\sum_{n \in \Z} e^{-n^2 \pi s} = s^{-1/2} \sum_{n \in \Z} e^{- n^2
\pi/s} , \qquad s > 0 .
\]
In terms of the \textbf{theta function} $\vartheta(s)=\sum_{n \in
\Z} e^{-n^2 \pi s}$, the last formula expresses the functional
equation
\[
\vartheta(s) = s^{-1/2} \vartheta(s^{-1}) .
\]
\noindent \Proof Decay hypotheses \eqref{eq:psf1} and
\eqref{eq:psf2} hold for $G_\omega$  because
\begin{align*}
G_\omega(x) & = \frac{\omega}{\sqrt{2\pi}} e^{-(\omega x)^2/2} , && x \in \R , \\
\widehat{G_\omega}(\xi) & = e^{-(\xi/\omega)^2/2} , && \xi \in \R
,
\end{align*}
by \eqref{eq:GR3b} and Table~\ref{ta:ftkernels}. Hence the Poisson
Summation Formula says that
\begin{align*}
\Pe(G_\omega)(x) & = \sum_{j \in \Z} \widehat{G_\omega}(j) e^{ijx} \\
& = \sum_{j \in \Z} e^{-(j/\omega)^2/2} e^{ijx} \\
& = \sum_{j \in \Z} e^{-j^2 s} e^{ijx} && \text{since $\omega=1/\sqrt{2s}$} \\
& = G_s(x) ,
\end{align*}
which proves \eqref{eq:pg}.

Taking $x=0$ in \eqref{eq:pg} and changing $s$ to $\pi s$ yields the
functional equation for the theta function.
\end{example}

\chapter{Uncertainty principles} \label{ch:up}

\subsubsection*{Goal}

Establish qualitative and quantitative uncertainty principles

\subsubsection*{References}

\cite{GM} Section 2

\noindent \cite{J} Section 1

\vspace{18pt}

Uncertainty principles say that $f$ and $\fhat$ cannot both be too
localized. Consequently, if $\fhat$ is well localized then $f$ is
not, and so we are ``uncertain'' of the value of $f$.

\begin{proposition}[Qualitative uncertainty principles]
\label{pr:qup}\

(a) If $f \in \ltt$ is continuous, $f$ has infinitely many zeros in
$\T$, and $\fhat$ is finitely supported, then $f \equiv 0$.

(b) If $f \in \ltrd$ is continuous, $f$ vanishes on some open set,
and $\fhat$ is compactly supported, then $f \equiv 0$.
\end{proposition}

\begin{proof}\

(a) $f$ is a trigonometric polynomial since it has only finitely
many nonzero Fourier coefficients. Thus part (a) says:
\begin{quotation}
a trigonometric polynomial that vanishes infinitely often in $\T$
must vanish identically.
\end{quotation}
To prove this claim, write $f(t)=\sum_{n=-N}^N a_n e^{int}$. Then
$f(t)=p(e^{it})/e^{iNt}$ where $p$ is the polynomial
\[
p(z) = \sum_{n=0}^{2N} a_{n-N} z^n , \qquad z \in \C .
\]
Since $f$ has infinitely many zeros $t \in \T$, we see $p$ has
infinitely many zeros $e^{it}$ on the unit circle. The Fundamental
Theorem of Algebra implies $p \equiv 0$.

(b) $f$ is band limited, and hence is holomorphic on $\Cd$ by
Theorem~\ref{th:blh}. In particular, $f$ is real analytic on $\Rd$.

Choose $x_0 \in \Rd$ such that $f \equiv 0$ on a neighborhood of
$x_0$; then the Taylor series of $f$ centered at $x_0$ is
identically zero. That Taylor series equals $f$ on $\Rd$, and so $f
\equiv 0$.
\end{proof}

\begin{theorem}[Benedicks' qualitative uncertainty principle]
\label{th:bqup} If $f \in \ltrd$ is continuous and $f$ and $\fhat$
are supported on sets of finite measure, then $f \equiv 0$.
\end{theorem}

In contrast to Proposition~\ref{pr:qup}, here the support of $\fhat$
need not be compact.

\begin{proof}
We prove only the $1$ dimensional case.

Let $A = \{ x \in \R : f(x) \neq 0 \}$ and $B = \{ \xi \in \R :
\fhatxi \neq 0 \}$. By dilating $f$ we can suppose $|A| < 2\pi$.
Then
\begin{align*}
& \quad \  \big| \{ x \in \T : \text{$f(x+2\pi n) \neq 0$ for some
$n \in \Z$} \} \big| \\
& = \big| \{ x \in \T : \sum_{n \in \Z} \charfn_A(x+2\pi n) \geq 1
\} \big| \\
& \leq \int_\T \sum_{n \in \Z} \charfn_A(x+2\pi n) \, dx \\
& = \int_\R \charfn_A(x) \, dx \\
& = |A| \\
& < |\T| = 2\pi .
\end{align*}
Therefore the complementary set
\[
E = \{ x \in \T : \text{$f(x+2\pi n) = 0$ for all $n \in \Z$} \}
\]
has positive measure.

Next,
\begin{align*}
\int_{[0,1)} \sum_{j \in \Z} \charfn_B(\xi+j) \, d\xi
& = \int_\R \charfn_B(\xi) \, d\xi \\
& = |B| \\
& < \infty ,
\end{align*}
so that $\sum_{j \in \Z} \charfn_B(\xi+j)$ is finite for almost
every $\xi \in [0,1)$, say for all $\xi \in F \subset [0,1)$ where
$F$ has full measure, $|[0,1) \setminus F|=0$. Hence when $\xi \in
F$, the set $\{ j \in \Z : \fhat(\xi+j) \neq 0 \}$ is finite.

Fix $\xi \in F$ and consider the periodization
\[
\Pe(f e^{-i\xi x})(x) = 2\pi \sum_{n \in \Z} f(x+2\pi n)
e^{-i\xi(x+2\pi n)} ,
\]
which is well defined since $f \in \lor$. The $j$th Fourier
coefficient of the periodization equals
\[
(fe^{-i\xi x})\widehat{\ }(j) = \fhat(\xi+j) ,
\]
which equals zero for but finitely many $j$, since $\xi \in F$. Thus
$\Pe(fe^{-i\xi x})$ equals some trigonometric polynomial $Q(x)$ a.e.
But $\Pe(fe^{-i\xi x})(x) = 0$ for all $x \in E$, and so $Q$
vanishes a.e.\ on $E$. In particular, $Q$ vanishes at infinitely
many points in $\T$ (using here that $E$ has positive measure).
Hence $Q \equiv 0$ by Proposition~\ref{pr:qup}(a). The Fourier
coefficient $\fhat(\xi+j)$ of $Q$ therefore vanishes for all $j$.

Since $\fhat(\xi+j)=0$ for all $j \in \Z$ and almost every $\xi \in
[0,1)$, we deduce $\fhatxi=0$ a.e., and so $f \equiv 0$.
\end{proof}

\begin{theorem}[Nazarov's quantitative uncertainty principle]
\label{th:nup} A constant $C_d>0$ exists such that
\[
\lv f \rv_\ltrd^2 = \lv \fhat \rv_\ltrd^2 \leq C_d^{\lvert A \rvert
\lvert B \rvert + 1} \Big( \int_{\Rd \setminus A} |f(x)|^2 \, dx +
\int_{\Rd \setminus B} |\fhatxi|^2 \, d\xi \Big)
\]
for all sets $A,B \subset \Rd$ of finite measure and all $f \in
\ltrd$.
\end{theorem}

We omit the proof.

Nazarov's theorem implies Benedicks' theorem, because if $f$ is
supported in $A$ and $\fhat$ is supported in $B$, then the right
side is zero and so $f \equiv 0$.

Next we develop an abstract commutator inequality that leads to the
Heisenberg Uncertainty Principle.

Let $H$ be a Hilbert space. Suppose $T$ is a linear operator from a
subspace ${\mathcal D}(T)$  into $H$. Write $T^*$ for its adjoint,
defined on a subspace ${\mathcal D}(T^*)$, meaning $T^*$ is linear
and
\[
\la Tf , g \ra = \la f , T^*g \ra \qquad \text{whenever} \quad f \in
{\mathcal D}(T), \quad g \in {\mathcal D}(T^*) .
\]
Define
\begin{align*}
\Delta_f(T) & = \min_{\alpha \in \C} \lv Tf - \alpha f \rv \\
& = \text{norm of component of $Tf$ perpendicular to $f$.}
\end{align*}
The minimum is attained for $\alpha = \la Tf , f \ra / \lv f \rv^2$.

\begin{theorem}[Commutator estimate] \label{th:comm}
Let $T$ and $U$ be linear operators like above. Then
\[
\big| \la [T,U]f , f \ra \big| \leq \Delta_f(T^*) \Delta_f(U) +
\Delta_f(T) \Delta_f(U^*)
\]
for all $f \in {\mathcal D} (TU) \cap {\mathcal D} (UT) \cap
{\mathcal D} (T^*) \cap {\mathcal D} (U^*)$.
\end{theorem}

Here $[T,U]=TU-UT$ is the \emph{commutator} of $T$ and $U$.

\begin{proof}
\begin{align}
\big| \la [T,U]f , f \ra \big| & = \big| \la TUf, f \ra - \la UTf , f \ra \big| \notag \\
& = \big| \la Uf , T^* f \ra - \la Tf , U^* f \ra \big| \notag \\
& \leq \lv Uf \rv \lv T^* f \rv + \lv Tf \rv \lv U^* f \rv .
\label{eq:comm1}
\end{align}
Let $\alpha , \beta \in \C$. Note that
\[
[T-\alpha I, U-\beta I] = [T,U] .
\]
Hence by replacing $T$ with $T-\alpha I$ and $U$ with $U-\beta I$ in
\eqref{eq:comm1} we find
\[
\big| \la [T,U]f , f \ra \big| \leq \lv Uf - \beta f \rv \lv T^* f -
\overline{\alpha} f \rv + \lv Tf - \alpha f \rv \lv U^* f -
\overline{\beta} f \rv .
\]
Minimizing over $\alpha$ and $\beta$ proves the theorem, noting for
the adjoints that
\[
\alpha = \frac{\la Tf , f \ra}{\lv f \rv^2} \quad
\Longleftrightarrow \quad \overline{\alpha} = \frac{\la T^*f , f
\ra}{\lv f \rv^2} .
\]
\end{proof}

\begin{example}[Heisenberg Uncertainty Principle]
Take $H=\ltr$,
\begin{align*}
(Tf)(x) = xf(x) & \qquad \text{with} \quad {\mathcal D}(T) = \{ f
\in \ltr : xf(x) \in \ltr \} , \\
(Uf)(x) = -if^\prime(x) & \qquad \text{with} \quad {\mathcal D}(U) =
\{ f \in \ltr : f^\prime \in \ltr \} .
\end{align*}
Here $T$ is the \emph{position} operator and $U$ is the
\emph{momentum} operator.

Observe $T^*=T, U^*=U$ and
\begin{align*}
[T,U]f & = TUf - UTf \\
& = x \cdot \big( -i \frac{d\ }{dx}f(x) \big) + i \frac{d\ }{dx}
\big( xf(x) \big) \\
& = if(x) .
\end{align*}
The Commutator Theorem~\ref{th:comm} implies
\begin{align*}
\lv f \rv_\ltr^2 & \leq 2 \Delta_f(T) \Delta_f(U) \\
& \leq 2 \lv xf - \alpha f \rv_\ltr \, \lv -if^\prime
- \beta f \rv_\ltr \\
& = 2 \lv (x-\alpha)f \rv_\ltr \, \frac{1}{\sqrt{2\pi}} \lv (\xi -
\beta) \fhat \rv_\ltr
\end{align*}
by Plancherel. Squaring yields the \textbf{Heisenberg Uncertainty
Principle}:
\begin{equation} \label{eq:hup1}
\frac{1}{4} \lv f \rv_\ltr^4 \leq \int_\R |x-\alpha|^2 |f(x)|^2 \,
dx \cdot \pr \int_\R |\xi-\beta|^2 |\fhatxi|^2 \, d\xi
\end{equation}
for all $\alpha, \beta \in \C$.

We interpret \eqref{eq:hup1} as restricting how localized $f$ and
$\fhat$ can be, around the locations $\alpha$ and $\beta$.

In quantum mechanics, we normalize $\lv f \rv_\ltr = 1$ and
interpret $|f(x)|^2$ as the probability density for the position $x$
of some particle, and regard $|\fhat|^2/2\pi$ as the probability
density for the momentum $\xi$. Thus the Heisenberg Uncertainty
Principle implies that the variance (or uncertainty) in position
multiplied by the variance in momentum is at least $1/4$.

Roughly, the Principle says that the more precisely one knows the
position of a quantum particle, the less precisely one knows its
momentum, and vice versa.
\end{example}

\begin{remark}\rm \

1. Equality holds in the Heisenberg Principle \eqref{eq:hup1} if and
only if $f(x) = Ce^{i\beta x} e^{-\gamma(x-\alpha)^2}$ is a
$\beta$-modulated Gaussian at $\alpha$ (with $C \in \C, \gamma>0$).

\vspace{3pt} 2. A more direct proof of \eqref{eq:hup1} can be given
by integrating by parts in
\[
\lv f \rv_\ltr^2 = \int_\R f(x) \overline{f(x)} (x-\alpha)^\prime \,
dx
\]
and then applying Cauchy--Schwarz.

\vspace{3pt} 3. The Heisenberg Uncertainty Principle extends
naturally to higher dimensions.

\vspace{3pt} 4. On $\T$, the analogous uncertainty principle says
\[
\frac{1}{4} m^2 \Big| \intp_\T e^{imt} |f(t)|^2 \, dt \Big|^2 \leq
\intp_\T |e^{imt}-\alpha|^2 |f(t)|^2 \, dt \cdot \sum_{n \in \Z}
|n-\beta|^2 |\fhatn|^2
\]
for all $\alpha, \beta \in \C, m \in \Z$ (exercise).

One considers here a quantum particle at position $e^{it}$ on the
unit circle, with momentum $n \in \Z$. When $\alpha=0$ we deduce a
lower bound on the localization of momentum, in terms of Fourier
coefficients of the position density $|f|^2$:
\[
\frac{1}{4} \lv f \rv_\ltt^{-2} \sup_{m \in \Z} m^2 \big|
(|f|^2)\widehat{\ }(m) \big|^2 \leq \sum_{n \in \Z} |n-\beta|^2
|\fhatn|^2 .
\]
\end{remark}

\part{Problems}

\chapter*{Assignment 1} \label{ch:a1}

\begin{problem}
Do the following problems, but do not hand them in:

\cite{K} Ex.\ 1.1.2, 1.1.4.
\end{problem}

\begin{problem} (\cite{K} Ex.\ 1.1.5: downsampling)

Let $f \in \lot, m \in \N$, and define
\[
f_{(m)}(t) = f(mt) .
\]

\noindent (a) Prove that $\widehat{f_{(m)}}(n)=\widehat{f}(n/m)$ if
$m | n$ and $\widehat{f_{(m)}}(n)=0 $ otherwise. Use only the
definition of the Fourier coefficients, and elementary
manipulations.

\noindent (b) Then give a quick, formal (nonrigorous) proof using
the Fourier \text{series of $f$}.
\end{problem}

\begin{problem} (\cite{K} Ex.\ 1.2.8: Fej\'{e}r's Lemma)

Let $f\in \lpt$ and $g \in L^q(\T)$, where $1 < p \leq \infty, 1
\leq q < \infty$ and $\frac{1}{p}+\frac{1}{q}=1$. Prove that
\[
\lim_{m \to \infty} \pr \int_\T f(mt) \overline{g(t)} \,
dt = \widehat{f}(0) \overline{\widehat{g}(0)} .
\] \emph{Hint.} Use that
trigonometric polynomials are dense in $L^q(\T)$.
\end{problem}

\begin{problem} (Weak convergence and oscillation)

Let $H$ be a Hilbert space. We say $u_n$ converges weakly to $u$,
written $u_n \rightharpoonup u$ weakly, if $\la u_n , v \ra \to \la
u , v \ra$ as $n \to \infty$, for each $v \in H$. Clearly if $u_n
\to u$ in norm (meaning $\lVert u_n - u \rVert \to 0$) then $u_n
\rightharpoonup u$ weakly.

(a) Show that $e^{imt} \rightharpoonup 0$ weakly in $\ltt$, as $m
\to \infty$.

(b) Let $f \in \ltt$. Show
\[
f_{(m)} \rightharpoonup \widehat{f}(0)=\text{(mean value of $f$)}
\]
weakly in $\ltt$, as $m \to \infty$.

\noindent \emph{Remark.} Thus rapid oscillation yields weak
convergence to the mean.
\end{problem}

\begin{problem} (Smoothness of $f$ implies rate of decay
of $\widehat{f}$)

(a) Show that if $f$ has bounded variation, then
$\widehat{f}(n)=O(|n|^{-1})$.

(b) Show that if $f$ is absolutely continuous and $f^\prime$ has
bounded variation, then $\widehat{f}(n)=O(|n|^{-2})$.

\noindent \emph{Remark.} These results cover most of the functions
encountered in elementary courses. For example, functions that are
smooth expect for finitely many jumps (such as the sawtooth $f(t)=t,
t \in (-\pi,\pi]$) have bounded variation. And functions that are
smooth except for finitely many corners (such as the triangular wave
$f(t)=|t|, t \in (-\pi,\pi]$) have first derivative with bounded
variation. That is why one encounters so many functions with Fourier
coefficients decaying like $1/n$ or $1/n^2$.
\end{problem}

\begin{problem} (\cite{K} Ex.\ 1.3.2: rate of uniform
summability)

Assume $f$ is H\"{o}lder continuous, with $f \in C^\alpha(\T)$ for
some $0<\alpha<1$. Prove there exists $C>0$ (depending on the
H\"{o}lder constant of $f$) such that
\[
\lVert \sigma_N(f) - f \rVert_\li \leq \frac{C}{1-\alpha} \,
\frac{1}{N^\alpha} , \qquad N \in \N .
\]
\noindent \emph{Remark.} Thus the ``smoother'' $f$ is, the faster
$\sigma_N(f)$ converges to $f$ as $N \to \infty$.
\end{problem}

\begin{problem} (\cite{K} Ex.\ 1.5.4)

Let $f$ be absolutely continuous on $\T$ with $f^\prime\in \ltt$. In
other words, $f \in W^{1,2}(\T)$.

(a) Prove that
\[
\lVert \widehat{f} \rVert_{\ell^1(\Z)} \leq \lVert f \rVert_\lot +
\left( 2 \sum_{n=1}^\infty \frac{1}{n^2} \right)^{\! \! 1/2} \lVert
f^\prime \rVert_\lt .
\]
\emph{Hint.} First evaluate $\lVert f^\prime \rVert_\lt^2$.

(b) Deduce that $f \in A(\T)$.

\noindent \emph{Remark.} Hence the Fourier series of $f$ converges
uniformly by Chapter~\ref{ch:l1}, so that $S_n(f) \to f$ in $\lit$.
In particular, if $f$ is smooth except for finite many corners (such
as the triangular wave $f(t)=|t|$ for $t \in (-\pi,\pi]$), then the
Fourier series converges uniformly to $f$.
\end{problem}

\begin{problem} (A lacunary series)

Assume $0<\alpha<1$.

(a) Suppose that $f$ is continuous on $\T$ and that
\[
\sum_{2^n \leq |j| < 2^{n+1}} |\widehat{f}(j)| \leq C 2^{-n\alpha}
\]
for each $n \geq 0$. Prove $f \in A(\T)$, and then $f \in
C^\alpha(\T)$.

(b) Let $f(t)=\sum_{n=0}^\infty 2^{-n\alpha} e^{i 2^n t}$. Show $f
\in C^\alpha(\T)$. Deduce that the rate of decay $\widehat{f}(n) =
O(|n|^{-\alpha})$ proved in Theorem~\ref{th:rd1} for $C^\alpha(\T)$
is sharp. (That is, show $\widehat{f}(n) = O(|n|^{-\beta})$ fails
for some $f \in C^\alpha(\T)$, when $\beta>\alpha$.)
\end{problem}

\begin{problem}\label{pr:llogl} (Maximal function when $p=1$)

Define $L \log L (\Rd)$ to be the class of measurable functions for
which $\int_\Rd |f(x)| \log \big( 1 + |f(x)| \big) \, dx < \infty$.
Prove that
\[
f \in L \log L (\Rd) \quad \implies \quad Mf \in L^1_{loc}(\Rd) .
\]
\noindent \emph{Remark.} Thus if the singularities of $f$ are
``logarithmically better than $\lo$'' then the Hardy--Littlewood
maximal function belongs to $\lo$ (at least locally).
\end{problem}

\begin{problem} Enjoyable reading (nothing to hand in).

Read Chapter~8 ``Compass and Tides'' from \cite{Ko}, which shows how
sums of Fourier series having different underlying periods can be
used to model the heights of tides.

Sums of periodic functions having different periods are called
\emph{almost periodic functions}. Their theory was developed by the
Danish mathematician Harald Bohr, brother of physicist Niels Bohr.
Harald Bohr won a silver medal at the 1908 Olympics, in soccer.
\end{problem}

\chapter*{Assignment 2} \label{ch:a2}

\begin{problem}[Hilbert transform of indicator function]\

(a) Evaluate $(H \charfn_{[a,b]})(t)$, where $[a,b] \subset
(-\pi,\pi)$ is a closed interval. Sketch the graph, for $t \in
[-\pi,\pi]$.

(b) Conclude that the Hilbert transform on $\T$ is not strong
$(\infty,\infty)$.
\end{problem}

\begin{problem}[Fourier synthesis on $\ell^p$] Let $1 \leq p \leq 2$.

Prove that the Fourier synthesis operator $T$, defined by
\[
(T \{ c_n \})(t) = \sum_{n \in \Z} c_n e^{int} ,
\]
is bounded from $\ell^p(\Z)$ to $L^{p^\prime}(\T)$. Estimate the
norm of $T$.

\emph{Extra credit.} Show the series converges unconditionally, in
$L^{p^\prime}(\T)$.
\end{problem}

\begin{problem}[Parseval on $L^p$] Do part (a) or part (b). You may do both
parts if you wish.

(a) Let $1 \leq p \leq 2$. Take $f \in \lpt$ and $g \in \lot$ with
$\{ \ghatn \} \in \ell^p (\Z)$. Prove that $g \in L^{p^\prime}(\T)$,
and establish the Parseval identity
\[
\intp_\T f(t) \overline{g(t)} \, dt = \sum_{n \in \Z} \fhatn
\overline{\ghatn} .
\]
(In your solution, explain why the integral and sum are absolutely
convergent.)

(b) Let $1<p<\infty$. Take $f \in \lpt$ and $g \in
L^{p^\prime}(\T)$. Prove the Parseval identity
\[
\intp_\T f(t) \overline{g(t)} \, dt = \lim_{N \to \infty} \sum_{|n|
\leq N} \fhatn \overline{\ghatn} .
\]
\end{problem}

\begin{problem}[Fourier analysis into a weighted space] Let $1 < p \leq 2$.

(a) Show
\[
\left( \sum_{n \neq 0} |\fhatn|^p |n|^{p-2} \right)^{\! \! \! 1/p}
\leq C_p \lv f \rv_\lpt \qquad \text{for all $f \in \lpt$.}
\]
\emph{Hint.} $Y=\Z \setminus \{ 0 \}$ with $\nu=$ counting measure
weighted by $n^{-2}$.

(b) Show that combining the H\"{o}lder and Hausdorff--Young
inequalities in the obvious way does \emph{not} prove part (a).
\end{problem}

\begin{problem}[Poisson extension] \label{pr:pe}
Recall $P_r$ denotes the Poisson kernel on $\T$, and write $\D$ for
the open unit disk in the complex plane. Suppose $f \in C(\T)$ and
define
\[
v(re^{it})= \begin{cases} (P_r * f)(t) & \text{for $0 \leq r
< 1 , \ t \in \T$,} \\
f(t) & \text{for $r=1, \ t \in \T$,}
\end{cases}
\]
so that $v$ is defined on the closed disk $\overline{\D}$.

(a) Show $v$ is $C^\infty$ smooth and harmonic ($\Delta v = 0$) in
$\D$.

(b) Show $v$ is continuous on $\overline{\D}$.

(c) [Optional; no credit] Assume $f \in C^\infty(\T)$ and show $v
\in C^\infty(\overline{\D})$. (Parts (a) and (b) show $v$ is smooth
on $\D$ and continuous on $\overline{\D}$. Thus the task is to prove
each partial derivative of $v$ on $\D$ extends continuously to
$\overline{\D}$.).

\noindent \emph{Aside.} $(P_r * f)(t)$ is called the \emph{harmonic
extension} to the disk of the boundary function $f$.
\end{problem}

\begin{problem}[Boundary values lose half a derivative]
Assume $u$ is a smooth, real-valued function on a neighborhood of
$\overline{\D}$, and define
\[
f(t) = u(e^{it})
\]
for the boundary value function of $u$. Hence $f \in C^\infty(\T)$,
and so the Poisson extension $v$ belongs to
$C^\infty(\overline{\D})$ by Problem~\ref{pr:pe}(c).

(a) Prove
\[
\pr \int_\D |\nabla v|^2 \, dA = \sum_{n \in \Z} |n|
|\fhatn|^2 .
\]
\emph{Hint.} Use one of Green's formulas, and remember
$v=\overline{v}$ since $f$ and $v$ are real-valued.

(b) Prove
\[
\int_\D |\nabla v|^2 \, dA \leq \int_\D |\nabla u|^2 \, dA .
\]
\emph{Hint.} Write $u=v+(u-v)$ and use one of Green's formulas.

\emph{Aside.} This result is known as ``Dirichlet's principle''. It
asserts that among all functions having the same boundary values,
the harmonic function has smallest Dirichlet integral. As your proof
reveals, this result holds on arbitrary domains.

(c) Conclude
\[
\sum_{n \in \Z} |n| |\fhatn|^2 \leq \pr  \int_\D |\nabla
u|^2 \, dA .
\]

\noindent \emph{Discussion.} We say $f$ has ``half a derivative'' in
$L^2$, since $\{ |n|^{1/2} \fhatn \} \in \ell^2(\Z)$. Justification:
if $f$ has zero derivatives ($f \in \ltt$) then $\{ \fhatn \} \in
\ell^2(\Z)$, and if $f$ has one derivative ($f^\prime \in \ltt$)
then $\{ n \fhatn \} \in \ell^2(\Z)$. Halfway inbetween lies the
condition $\{ |n|^{1/2} \fhatn \} \in \ell^2(\Z)$.

\emph{Boundary trace} inequalities like in part (c) are important
for partial differential equations and Sobolev space theory. The
inequality says, basically, that if a function $u$ has one
derivative $\nabla u$ belonging to $L^2$ on a domain, then $u$ has
half a derivative in $L^2$ on the boundary. Thus the boundary value
loses half a derivative, compared to the original function.

Note that in this problem, $f \in C^\infty(\T)$ and so certainly
$f^\prime \in \ltt$, which implies $\{ n \fhatn \} \in \ell^2(\Z)$.
You might wonder, then, why you should bother proving the weaker
result $\{ |n|^{1/2} \fhatn \} \in \ell^2(\Z)$ in part (c). But
actually you prove more in part (c): you obtain a \emph{norm
estimate} on $\{ |n|^{1/2} \fhatn \} \in \ell^2(\Z)$ in terms of the
$L^2$ norm of $\nabla u$. (We do not have such a norm estimate on
$\{ n \fhatn \}$.) This norm estimate means that the restriction map
\begin{align*}
H^1(\D) & \to H^{1/2}(\partial \D) \\
u & \mapsto f
\end{align*}
is bounded from the Sobolev space $H^1(\D)$ on the disk with one
derivative in $L^2$ to the Sobolev space $H^{1/2}(\partial \D)$ on
the boundary circle with half a derivative in $L^2$.

\vspace{6pt} \noindent \emph{Aside.} The notion of fractional
derivatives defined via Fourier coefficients can be extended to
fractional derivatives in $\Rd$, by using Fourier transforms.
\end{problem}

\begin{problem}[Measuring diameters of stars]\

Enjoyable reading; nothing to hand in.

Read Chapter~95 ``The Diameter of Stars'' from \cite{Ko}, which
shows how the diameters of stars can be estimated using Fourier
transforms of radial functions, and convolutions.
\end{problem}

\chapter*{Assignment 3} \label{ch:a3}

\begin{problem}[Adjoint of Fourier transform]\

Find the adjoint of the Fourier transform on $\ltrd$.
\end{problem}

\begin{problem}[Periodization, and Fourier coefficients and
transforms]\label{pr:pef} \

\vspace{3pt} Suppose $f \in \lord$.

\vspace{6pt} (a) Prove that the periodization
\[
\Pe(f)(x)=(2\pi)^d \sum_{n \in \Zd} f(x+2\pi n)
\]
of $f$ satisfies
\[
\lv \Pe(f) \rv_\lotd \leq \lv f \rv_\lord .
\]

\vspace{6pt} (b) Deduce from your argument that the series for
$\Pe(f)(x)$ converges absolutely for almost every $x$, and that
$\Pe(f)$ is $2\pi \Zd$-periodic.

\vspace{6pt} (c) Show that the $j$th Fourier coefficient of $\Pe(f)$
equals the Fourier transform of $f$ at $j$:
\[
\widehat{\Pe(f)}(j) = \fhatj , \qquad j \in \Zd
\]

\end{problem}

\begin{problem}[Course summary]\

Write a one page description of the most important and memorable
results and general techniques from this course. Be brief, but
thoughtful; explain how these main results fit together.

You need not state the results technically --- intuition is more
helpful than rigor, at this stage.
\end{problem}

\part{Appendices}

\appendix

\chapter{Minkowski's integral inequality} \label{ap:mii}

\subsubsection*{Goal} State Minkowski's integral inequality, and apply it to norms of convolutions

\vspace{18pt} Minkowski's inequality on a measure space $(X,\mu)$ is
simply the triangle inequality for $L^p(X)$, saying that the norm of
a sum is bounded by the sum of the norms:
\[
\big\lv \sum_j f_j \nu_j \big\rv_{L^p(X)} \leq \sum_j \lv f_j \rv_{L^p(X)} \nu_j
\]
whenever $f_j \in L^p(X)$ and the constants $\nu_j$ are nonnegative.
Similarly, the norm of an integral is bounded by the integral of the
norms:

\begin{theorem} \label{th:mii}
Suppose $(X,\mu)$ and $(Y,\nu)$ are $\sigma$-finite measure spaces, and that $f(x,y)$ is measurable on the product space $X \times Y$. If $1 \leq p \leq \infty$ then
\[
\Big\lv \int_Y f(x,y) \, d\nu(y) \Big\rv_{L^p(X)} \leq \int_Y \lv f(x,y) \rv_{L^p(X)} \, d\nu(y)
\]
whenever the right side is finite.
\end{theorem}

\begin{proof}
Take $q$ to be the conjugate exponent, with $\frac{1}{p}+\frac{1}{q}=1$. Then for all $g \in L^q(X)$,
\begin{align*}
& \Big| \int_X \Big( \int_Y f(x,y) \, d\nu(y) \Big) g(x) \, d\mu(x) \Big| \\
& \leq \int_Y \int_X \lvert f(x,y) \rvert \lvert g(x) \rvert \, d\mu(x) d\nu(y) \\
& \leq \int_Y \Big( \int_X |f(x,y)|^p \, d\mu(x) \Big)^{\! \! 1/p} \lv g \rv_{L^q(X)} \, d\nu(y) && \text{by H\"{o}lder} \\
& = \int_Y \lv f(x,y) \rv_{L^p(X)} \, d\nu(y) \cdot \lv g \rv_{L^q(X)} .
\end{align*}
Now the theorem follows from the dual characterization of the norm
on $\lp(X)$ (see \cite[Theorem 6.14]{F}).
\end{proof}

\begin{definition} Define the \emph{convolution} of functions $f$ and $g$ on $\T$ by
\[
(f*g)(t) = \intp_\T f(t-\tau) g(\tau) \, d\tau , \qquad t \in \T ,
\]
whenever the integral makes sense.

Define the \emph{convolution} of functions $f$ and $g$ on $\Rd$ by
\[
(f*g)(x) = \int_\Rd f(x-y) g(y) \, dy , \qquad x \in \Rd ,
\]
whenever the integral makes sense.
\end{definition}

\begin{theorem}[Young's theorem] \label{th:yt} Fix $1 \leq p \leq \infty$. Then
\begin{align*}
\lv f*g \rv_\lpt & \leq \lv f \rv_\lpt \lv g \rv_\lot , \\
\lv f*g \rv_\lprd & \leq \lv f \rv_\lprd \lv g \rv_\lord ,
\end{align*}
whenever the right sides are finite. In particular, the convolution $f*g$ is well defined a.e.\ whenever $f \in \lp$ and $g \in \lo$.
\end{theorem}

\begin{proof}
\begin{align*}
\lv f*g \rv_\lprd & = \Big\lv \int_\Rd f(\cdot-y) g(y) \, dy \Big\rv_\lprd \\
& \leq \int_\Rd \lv f(\cdot-y) \rv_\lprd |g(y)| \, dy \\
& \qquad \qquad \qquad \text{by Minkowski's integral inequality, Theorem~\ref{th:mii},} \\
& = \lv f \rv_\lprd \lv g \rv_\lord .
\end{align*}
Argue similarly on $\T$.
\end{proof}

\chapter{$L^p$ norms and the distribution function} \label{ap:df}

\subsubsection*{Goal}

Express $L^p$-norms in terms of the distribution function

\vspace{18pt} Given a $\sigma$-finite measure space $(X,\mu)$ and a measurable
function $f$ on $X$, write
\[
E(\lambda) = \{ x \in X : |f(x)| > \lambda \}
\]
for the level set of $f$ above level $\lambda$. The
\emph{distribution function} of $f$ is $\mu(E(\lambda))$.

\begin{lemma} \label{le:df}
Let $\alpha > 0$.

If $-\infty < r < p < \infty$ then
\begin{equation}
\label{le:df1} \int_0^\infty \lambda^{p-r-1}
\int_{E(\lambda/\alpha)} |f(x)|^r \, d\mu(x) d\lambda =
\frac{\alpha^{p-r}}{p-r} \int_X |f(x)|^p \, d\mu(x) .
\end{equation}

If $-\infty < p < r < \infty$ then
\begin{equation}
\label{le:df2} \int_0^\infty \lambda^{p-r-1}
\int_{E(\lambda/\alpha)^c} |f(x)|^r \, d\mu(x) d\lambda =
\frac{\alpha^{p-r}}{r-p} \int_X |f(x)|^p \, d\mu(x) .
\end{equation}

In particular, when $r = 0 < p < \infty$ and $\alpha=1$, formula
\eqref{le:df1} expresses the $L^p$-norm in terms of the distribution
function:
\begin{equation}
\label{le:df3} \int_0^\infty p \lambda^{p-1} \mu(E(\lambda)) \,
d\lambda = \int_X |f(x)|^p \, d\mu(x) = \lv f \rv_{L^p(X)}^p .
\end{equation}
\end{lemma}
\begin{proof}
We can assume $\alpha=1$ without loss of generality, by changing
variable with $\lambda \mapsto \alpha \lambda$.

Write $E=\{ (x,\lambda) \in X \times (0,\infty) : |f(x)|>\lambda
\}$, so that $(x,\lambda) \in E \Leftrightarrow x \in E(\lambda)$.
Then the left hand side of \eqref{le:df1} equals
\begin{align*}
& \int_0^\infty \lambda^{p-r-1} \int_X \charfn_E(x,\lambda) |f(x)|^r
\, d\mu(x) d\lambda \\
& = \int_X |f(x)|^r \int_0^\infty \lambda^{p-r-1}
\charfn_E(x,\lambda) \, d\lambda d\mu(x) && \text{by Fubini} \\
& = \int_X |f(x)|^r \int_0^{|f(x)|} \lambda^{p-r-1}
\, d\lambda d\mu(x) && \text{since $\lambda < |f(x)|$ on $E$} \\
& = \int_X |f(x)|^r \frac{1}{p-r} |f(x)|^{p-r} \, d\mu(x)
\end{align*}
since $p-r>0$. Thus we have proved \eqref{le:df1}, and
\eqref{le:df2} is similar.
\end{proof}

\chapter{Interpolation} \label{ap:ip}

\subsubsection*{Goal}

Interpolation of operators on $L^p$ spaces, assuming either weak
endpoint bounds (Marcinkiewicz) or strong endpoint bounds
(Riesz--Thorin)

\subsubsection*{References}

\cite{F} Chapter 6

\noindent \cite{G} Section 1.3

\vspace{18pt}

\begin{definition}
An operator is \emph{sublinear} if
\[
|T(f+g)(y)| \leq |(Tf)(y)|+|(Tg)(y)|
\]
\[
|T(cf)(y)|=|c| |(Tf)(y)|
\]
for all $f, g$ in the domain of $T$, all $c \in \C$, and all $y$ in
the underlying set.
\end{definition}

\begin{theorem}[Marcinkiewicz Interpolation] \label{th:mi}
Let $1 \leq p_0 < p_1 \leq \infty$ and suppose $(X, \mu)$ and $(Y,
\nu)$ are measure spaces. Assume
\[
T: L^{p_0} + L^{p_1}(X) \to \{\text{measurable functions on} \ Y \}
\]
is sublinear. If $T$ is weak $(p_0, p_0)$ and weak $(p_1, p_1)$,
then $T$ is strong $(p,p)$ whenever $p_0 < p < p_1$.
\end{theorem}

\begin{proof}
Write $A_0, A_1$ for the constants in the weak $(p_0, p_0)$ and
$(p_1, p_1)$ estimates. Let $\alpha > 0$. Consider $f \in \lp(X),
\lambda >0$. Split $f$ into ``large'' and ``small'' parts:
\[
g = f \charfn_{\{x : |f(x)| > \lambda/\alpha \}} \quad \text{and}
\quad h = f \charfn_{\{x : |f(x)| \leq \lambda/\alpha \}} .
\]
Notice that
\begin{align*}
& \text{$g \in L^{p_0}(X)$ since $|g|^{p_0} \leq |f|^p (\lambda /
\alpha)^{p_0-p}$,} \\
& \text{$h \in L^{p_1}(X)$ since $|h|^{p_1} \leq |f|^p (\lambda /
\alpha)^{p_1-p}$.} \end{align*}
Hence $f = g + h \in L^{p_0} + L^{p_1}(X)$. By sublinearity, $|Tf|
\leq |Tg| + |Th|$.

\subparagraph*{Case 1.} Assume $p_1 < \infty$. Then
\begin{align}
& \nu \left( \left\{ y \in Y : |Tf(y)| > \lambda \right\} \right) \notag \\
& \leq \nu \left( \left\{ y \in Y : |Tg(y)| > \lambda/2 \right\}
\right) + \nu \left( \left\{ y \in Y : |Th(y)|>\lambda/2 \right\}
\right) \quad \text{by sublinearity} \notag \\
& \leq \left(\frac{A_0}{\lambda/2} \lv g \rv_{L^{p_0}(X)}
\right)^{\! \! p_0} + \left( \frac{A_1}{\lambda/2} \lv h
\rv_{L^{p_1}(X)} \right)^{\! \! p_1} \qquad \text{by the weak estimates on $T$} \notag \\
& = (2A_0)^{p_0} \lambda^{-p_0} \int_{\{x : |f(x)| > \lambda/\alpha\}} |f(x)|^{p_0} \, d\mu(x) \notag \\
& + (2A_1)^{p_1}\lambda^{-p_1} \int_{\{x : |f(x)| \leq
\lambda/\alpha\}} |f(x)|^{p_1} \, d\mu(x) . \label{eq:mi1}
\end{align}
Therefore
\begin{align*}
\lv Tf \rv_{\lp(Y)}^p & = \int_0^{\infty} p \lambda^{p-1} \nu \left( \left\{ y \in Y : |Tf(y)| > \lambda \right\} \right) d\lambda \\
& \leq p(2A_0)^{p_0} \frac{\alpha^{p-p_0}}{p-p_0} \lv f
\rv^p_{\lp(X)} + p(2A_1)^{p_1} \frac{\alpha^{p-p_1}}{p_1-p} \lv f
\rv^p_{\lp(X)}
\end{align*}
by \eqref{eq:mi1} and formulas \eqref{le:df1}, \eqref{le:df2}. We
have proved strong $(p,p)$.

Choosing $\alpha = 2 A_1^{p_1/(p_1-p_0)}/A_0^{p_0/(p_1-p_0)}$ gives
simple constants:
\begin{equation} \label{eq:mi2}
\lv Tf \rv_{\lp(Y)} \leq 2p^{1/p}\left( \frac{1}{p-p_0} +
\frac{1}{p_1-p} \right)^{\! \! \! 1/p} A_0^{1-\theta} A_1^{\theta}
\lv f \rv_{\lp(X)}
\end{equation}
where $0<\theta<1$ is determined by expressing $\frac{1}{p}$ as a
convex combination of $\frac{1}{p_0}$ and
$\frac{1}{p_1}$:
\[
\frac{1}{p}=\frac{1-\theta}{p_0}+\frac{\theta}{p_1} .
\]
Note the estimate in \eqref{eq:mi2} blows up as $p$ approaches $p_0$
or $p_1$.

\subparagraph{Case 2.} Assume $p_1=\infty$. Let $\alpha=2A_1$. We
have $\lv Th \rv_{\li(Y)} \leq A_1 \lv h \rv_{\li(X)}$, because weak
$(\infty,\infty)$ is defined to mean strong $(\infty,\infty)$, and
so
\[
\lv Th \rv_{\li(Y)} \leq A_1 \frac{\lambda}{\alpha} =
\frac{\lambda}{2}
\]
by definitions of $h$ and $\alpha$. Hence
\[
\nu \left( \{ y \in Y : |Tf(y)|>\lambda \} \right) \leq \nu \left(
\{ y \in Y : |Tg(y)|> \lambda/2 \} \right)
\]
because $|Tf|\leq |Tg|+|Th|$. Now argue like in Case~1 to get strong
$(p,p)$.
\end{proof}

Next we weaken the hypotheses of Marcinkiewicz Interpolation.

\begin{definition}
Given a measure space $(X, \mu)$, write
\[
\Sigma(X)=\{\text{simple functions on $X$ with support of finite
measure}\}.
\]
That is, $f \in \Sigma(X)$ provided $f=\sum \alpha_j
\charfn_{F_j}$ where the sum is finite, $\alpha_j \in \C \setminus
\{ 0 \}$, and the sets $F_j$ have finite measure and are disjoint.
\end{definition}

\begin{remark}[Linear Operators]\rm \label{re:mi}
Suppose
\[
T : \Sigma(X) \to \{ \text{measurable functions on $Y$} \}
\]
is linear. Then Marcinkiewicz Interpolation still holds: if
$T$ is weak $(p_0,p_0)$ and weak $(p_1,p_1)$ on the simple
functions in $\Sigma(X)$, then $T$ is strong $(p,p)$ on $\lp(X)$
whenever $p_0 < p < p_1$.

\begin{proof}
If $f$ is simple with support of finite measure, then so are
$g=f\charfn_{\{ |f|> \lambda/ \alpha\} }$ and $h=f\charfn_{\{ |f|
\leq \lambda/ \alpha\} }$. And $Tf=Tg+Th$ by linearity. Hence the
proof of Marcinkiewicz Interpolation gives a strong $(p,p)$ bound
for $T$ on $\Sigma(X)$. By density of $\Sigma(X)$ in $\lp(X)$
(using here that $p<p_1$ implies $p < \infty$), we deduce $T$ has
a unique extension to a bounded linear operator on $\lp(X)$. (This
extension step uses linearity of $T$.)
\end{proof}
\end{remark}

Our next interpolation result needs:
\begin{lemma}[Hadamard's Three Lines] \label{le:htl}
Assume $H(z)$ is holomorphic on $U=\{z \in \mathbb{C} : 0< \Re(z)
<1\}$ and continuous and bounded on $\overline{U}=\{z \in \mathbb{C}
: 0 \leq \Re(z) \leq 1\}$. Let $B_0= \sup_{\Re(z)=0} |H(z)|$ and
$B_1=\sup_{\Re(z)=1}|H(z)|$.

Then $|H(z)| \leq B_0^{1-\theta} B_1^{\theta}$ whenever
$\Re(z)=\theta \in [0,1]$.
\end{lemma}

\noindent (Exercise. Let $B_{\theta}=\sup_{\Re(z)=\theta}|H(z)|$, so
that $B_{\theta}\leq B_0^{1-\theta} B_1^{\theta}$ by the Lemma. Show
that $\theta \mapsto \log B_{\theta}$ is convex.)

\begin{proof}
Assume $B_0>0$ and $B_1>0$. Then
\[
G(z)=\frac{H(z)}{B_0^{1-z}B_1^z}
\]
is holomorphic on $U$ and bounded on $\overline{U}$, since $H$ is
bounded and $|B_0^{1-z}B_1^{z}|=B_0^{1-\Re(z)}B_1^{\Re(z)} \geq
\min(B_0, B_1) > 0$. Let $G_m=G(z)e^{(z^2-1)/m}, m>0$. Then $G_m$ is
holomorphic on $U$ with
\begin{align*}
|G_m(z)| &= |G(z)|e^{-(y^2+1)/m}e^{(x^2-1)/m} && \text{where $z=x+iy$} \\
& \leq |G(z)|e^{-(y^2+1)/m} && \text{since $x^2 \leq 1$ on
$\overline{U}$.}
\end{align*}
Hence $G_m \to 0$ as $|z| \to \infty$ in $U$. Therefore the Maximum
Principle applied to $G_m$ in $U$ says
\begin{align*}
\sup_{z \in U} |G_m(z)| & \leq \sup_{\partial U \cup \{\infty\}} |G_m| \\
& = \sup_{\partial U}|G_m| \\
& \leq \sup_{\partial U} |G| \\
& \leq 1 ,
\end{align*}
since $|H| \leq B_0$ on $\{\Re(z)=0\}$ and $|H|\leq B_1$ on
$\{\Re(z)=1\}$. Letting $m \to \infty$ gives $|G(z)| \leq 1$, which
proves the lemma.

If $B_0=0$ or $B_1=0$, then add $\e$ to $H$ and argue as above. Let
$\e \to 0$.
\end{proof}

\begin{theorem}[Riesz--Thorin Interpolation] \label{th:rti}
Let $1 \leq p_0, p_1, q_0, q_1 \leq \infty$, and $(X, \mu)$ and
$(Y, \nu)$ be measure spaces. (If $q_0=q_1=\infty$, then further
assume $\nu$ is semi-finite.) Suppose
\[
T : L^{p_0}+L^{p_1}(X) \to L^{q_0}+L^{q_1}(Y)
\]
is linear.

If $T$ is strong $(p_0,q_0)$ and $(p_1,q_1)$, then $T$ is strong
$(p,q)$ whenever
\[
\Big( \frac{1}{p}, \frac{1}{q} \Big) = (1-\theta) \Big(
\frac{1}{p_0},\frac{1}{q_0} \Big) + \theta \Big(
\frac{1}{p_1},\frac{1}{q_1} \Big)
\]
for some $0<\theta<1$. Specifically,
\[
\lv T \rv_{\lp(X)\to L^q(Y)} \leq \lv T
\rv^{1-\theta}_{L^{p_0}(X)\to L^{q_0}(Y)} \lv T
\rv^{\theta}_{L^{p_1}(X) \to L^{q_1}(Y)} .
\]
\end{theorem}

\begin{remark}\rm \

\noindent 1. The relationship between the $p$ and $q$ parameters is
shown in Figure~\ref{RTfig}. In particular, if $\theta =0$ then
$(p,q)=(p_0,q_0)$, and if $\theta =1$ then $(p,q)=(p_1,q_1)$.
\begin{figure}
\begin{center}
  \includegraphics[scale=0.6]{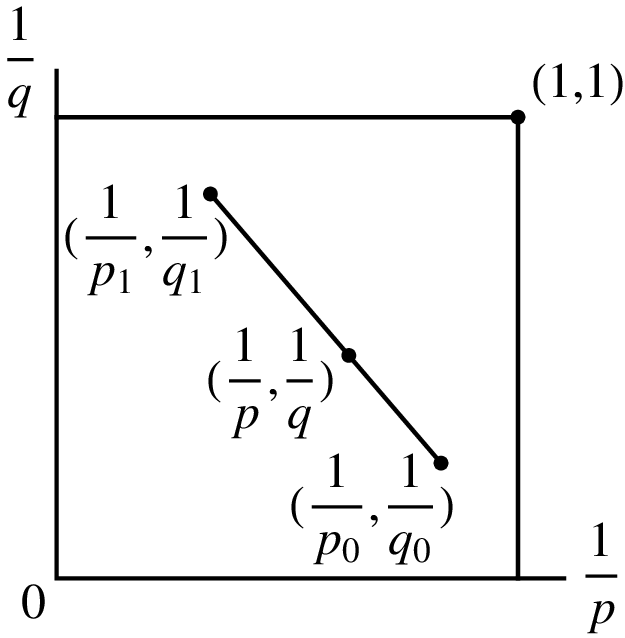}
  \caption{\label{RTfig}
    Parameters in the Riesz--Thorin theorem.}
\end{center}
\end{figure}

\vspace{6pt} \noindent 2. The space
\[
L^{p_0} + L^{p_1}(X) = \{ f_0+f_1 : f_0 \in L^{p_0}(X), f_1 \in
L^{p_1}(X) \}
\]
consists of all sums of functions in $L^{p_0}$ and $L^{p_1}$.
Recall from measure theory that
\[
\lp \subset L^{p_0} + L^{p_1} ,
\]
by splitting $f \in \lp$ into large and small parts.

A subtle aspect of the theorem is that when we assume $T$ maps
$L^{p_0} + L^{p_1}(X)$ into $L^{q_0} + L^{q_1}(Y)$, we need the
value of $Tf$ to be independent of the choice of decomposition
$f=f_0+f_1$.

In applications of the theorem, usually one has $T$ defined and
linear on $L^{p_0}(X)$ and $L^{p_1}(X)$, and the two definitions
agree on the intersection $L^{p_0} \cap L^{p_1}(X)$. Then one
defines $T$ on $f=f_0+f_1 \in L^{p_0}+L^{p_1}$ by $Tf=Tf_0+Tf_1$.
This definition is independent of the decomposition, as follows.
For suppose $f=\widetilde{f_0}+\widetilde{f_1}$. Then
\[
f_0-\widetilde{f_0} = \widetilde{f_1}-f_1 \in L^{p_0} \cap
L^{p_1}(X)
\]
and so $T \big( f_0-\widetilde{f_0} \big) = T \big(
\widetilde{f_1}-f_1 \big)$, where on the left side we use $T$
defined on $L^{p_0}(X)$ and on the right side we use $T$ on
$L^{p_1}(X)$. Linearity of $T$ now yields
$Tf_0+Tf_1=T\widetilde{f_0}+T\widetilde{f_1}$ so that the
definition of $Tf$ is independent of the decomposition of $f$.

\vspace{6pt} \noindent 3. When $T=\text{identity}$, Riesz--Thorin
says that
\[
L^{p_0} \cap L^{p_1} \subset \lp
\]
with
\begin{equation} \label{eq:api1}
\lv f \rv_{L^p(X)} \leq \lv f \rv_{L^{p_0}(X)}^{1-\theta} \lv f
\rv_{L^{p_1}(X)}^{\theta}
\end{equation}
where $\frac{1}{p}=\frac{1-\theta}{p_0}+\frac{\theta}{p_1}$. Here is
a direct proof:
\begin{align*}
\lv f \rv_{L^p(X)}^p & = \int_X |f|^p \, d\mu \\
& = \int_X |f|^{p(1-\theta)} |f|^{p\theta} \, d\mu \\
& \leq \Big( \int_X |f|^{p(1-\theta) \cdot p_0/p(1-\theta)} \,
d\mu \Big)^{\! p(1-\theta)/p_0} \Big( \int_X |f|^{p\theta \cdot p_1/p\theta} \, d\mu \Big)^{\! p\theta/p_1} && \text{by H\"{o}lder} \\
& = \lv f \rv^{p(1-\theta)}_{L^{p_0}(X)} \lv f
\rv^{p\theta}_{L^{p_1}(X)}
\end{align*}
\end{remark}

\begin{proof}[Proof of Riesz--Thorin Interpolation]
First suppose $p_0=p_1$, so that $p_0=p_1=p$. Then
\[
\lv Tf \rv_{L^q(Y)} \leq \lv Tf \rv_{L^{q_0}(Y)}^{1-\theta} \lv Tf
\rv_{L^{q_1}(Y)}^\theta
\]
by \eqref{eq:api1} applied to $Tf$ on $Y$. Now the $(p_0,q_0)$ and
$(p_1,q_1)$ bounds can be applied directly to give the $(p,q)$
bound.

Next suppose $p_0 \neq p_1$, so that $p<\infty$.

We will prove an $\lp \to L^q$ bound on $Tf$ for $f \in
\Sigma(X)$. Then at the end we will prove the bound for $ f \in
\lp(X)$.

Let $f \in \Sigma(X)$ and $g \in \Sigma(Y)$, say $f=\sum \alpha_j
\charfn_{F_j}$ and $g=\sum \beta_j \charfn_{G_j}$. Fix $\theta \in
(0,1)$, which fixes $p$ and $q$. For $z \in \C$, define
\begin{align*}
P(z) & =\frac{p}{p_0}(1-z)+\frac{p}{p_1}z , \\
Q^\prime(z) & =
\frac{q^\prime}{q_0^\prime}(1-z)+\frac{q^\prime}{q_1^\prime}z
\end{align*}
where $\frac{1}{q}+\frac{1}{q^\prime}=1$. (The ${}^\prime$ in
$Q^\prime$ does not denote a derivative, here.) Let
\begin{align*}
f_z(x) & = |f(x)|^{P(z)}e^{i \arg f(x)}, && x \in X ,\\
g_z(y) & = |g(y)|^{Q^\prime(z)}e^{i \arg g(y)}, && y \in Y .
\end{align*}
Note $f_{\theta}=f$ and $g_{\theta}=g$, since $P(\theta)=1$ and
$Q^\prime(\theta)=1$. Clearly
\begin{equation} \label{eq:rti1}
    g_z=\sum |\beta_j|^{Q^\prime(z)}e^{i \arg \beta_j}
    \charfn_{G_j}.
\end{equation}
Therefore $g_z(y)$ is bounded for $y \in Y, z \in \overline{U}$, and
it has support (independent of $z$) with finite measure. Similarly,
\begin{equation} \label{eq:rti2}
    Tf_z=\sum |\alpha_j|^{P(z)}e^{i \arg \alpha_j} (T\charfn_{F_j})
\end{equation}
by linearity, so that
\begin{align*}
|Tf_z| & \leq \sum |\alpha_j|^{\Re P(z)}|T\charfn_{F_j}| \\
& \leq (\text{const.}) \sum |T\charfn_{F_j}| && \text{for $z \in
\overline{U}$.}
\end{align*}
The right side belongs to $L^{q_0} \cap L^{q_1}(Y)$ by the strong
$(p_0,q_0)$ and $(p_1,q_1)$ bounds, since $\charfn_{F_j} \in L^{p_0}
\cap L^{p_1}(X)$. Hence the function
\begin{equation} \label{eq:rti3}
    H(z)=\int_Y (Tf_z)(y) g_z(y) \, d\nu(y)
\end{equation}
is well-defined and bounded for $z \in \overline{U}$, by H\"{o}lder.
And $H$ is holomorphic, as one sees by substituting \eqref{eq:rti1}
and \eqref{eq:rti2} into \eqref{eq:rti3} and taking the sums outside
the integral. Next,
\begin{align*}
\Re(z)=0 \quad \Rightarrow & \quad \Re P(z)=\frac{p}{p_0}, \qquad
\Re Q^\prime(z)=\frac{q^\prime}{q_0^\prime} \\
\Rightarrow & \quad |f_z|^{p_0}=|f|^{p_0 \Re P(z)}=|f|^p \\
& \quad |g_z|^{q_0^\prime}=|g|^{q_0^\prime \Re Q^\prime(z)}=|g|^{q^\prime} \\
\Rightarrow & \quad \lv f_z \rv_{L^{p_0}(X)}=\lv f
\rv^{p/p_0}_{L^p (X)} \\
& \quad \lv g_z \rv_{L^{q_0^\prime}(Y)}=\lv g \rv^{q^\prime/q_0^\prime}_{L^{q^\prime}(Y)} \\
& \qquad \qquad \qquad \text{(valid even when $p_0=\infty$ or $q_0=\infty$)} \\
\Rightarrow & \quad |H(z)| \leq \lv Tf_z \rv_{L^{q_0}(Y)}\lv g_z \rv_{L^{q_0^\prime}(Y)} \qquad \text{by H\"{o}lder} \\
& \quad |H(z)| \leq \lv T \rv_{L^{p_0}(X) \to L^{q_0}(Y)}\lv f
\rv_{\lp(X)}^{p/p_0} \lv g
\rv_{L^{q^\prime}(Y)}^{q^\prime/q_0^\prime} .
\end{align*}
Similarly,
\[
\Re(z)=1 \quad \Rightarrow \quad |H(z)| \leq \lv T \rv_{L^{p_1}(X)
\to L^{q_1}(Y)}\lv f \rv_{\lp(X)}^{p/p_1} \lv g
\rv_{L^{q^\prime}(Y)}^{q^\prime/q_1^\prime} .
\]
Hence by the Hadamard Three Lines Lemma~\ref{le:htl} and a short
calculation, if $z=\theta$ then
\begin{align*}
|\langle Tf,\overline{g} \rangle| & = |H(\theta)| \\
& \leq \lv T \rv^{1-\theta}_{L^{p_0}(X)\to L^{q_0}(Y)} \lv T
\rv^{\theta}_{L^{p_1}(X) \to L^{q_1}(Y)} \lv f \rv_{\lp(X)} \lv g
\rv_{L^{q^\prime}(Y)} .
\end{align*}
Now the dual characterization of the norm on $L^q(Y)$ implies
\begin{equation} \label{eq:rti4}
\lv Tf \rv_{L^q(Y)} \leq \lv T \rv_{L^{p_0}(X) \to
L^{q_0}(Y)}^{1-\theta} \lv T \rv_{L^{p_1}(X) \to
L^{q_1}(Y)}^\theta \lv f \rv_{L^p(X)} ,
\end{equation}
which is the desired strong $(p,q)$ bound. (See \cite[Theorem
6.14]{F} for the dual characterization of the norm, which uses
semi-finiteness of $\nu$ when $q=\infty$.)

We must extend this bound \eqref{eq:rti4} from $f \in \Sigma(X)$
to $f \in \lp(X)$. So fix $f \in \lp(X)$ and let $E=\{ x :
|f(x)|>1 \}$. Choose a sequence of simple functions $f_n \in
\Sigma(X)$ with $|f_n| \leq |f|$ and $f_n \to f$ at every point,
and with $f_n \to f$ uniformly on $X \setminus E$; such a sequence
exists by \cite[Theorem~2.10]{F}. Define
\[
g=f \charfn_E, \qquad g_n = f_n \charfn_E ,
\]
and
\[
h=f \charfn_{X \setminus E}, \qquad h_n = f_n \charfn_{X \setminus
E} ,
\]
so that $f=g+h, f_n=g_n+h_n$, and $|g_n| \leq |g|, |h_n| \leq
|h|$. Suppose $p_0 < p_1$, by swapping $p_0$ and $p_1$ if
necessary. Then $g_n \to g$ in $L^{p_0}(X)$ by dominated
convergence, and so $Tg_n \to Tg$ in $L^{q_0}(Y)$. By passing to a
subsequence we can further suppose $Tg_n \to Tg$ pointwise a.e.

Also $h_n \to h$ in $L^{p_1}(X)$ by dominated convergence (or, if
$p_1=\infty$, by the uniform convergence $f_n \to f$ on $X
\setminus E$). Hence $Th_n \to Th$ in $L^{q_1}(Y)$. By passing to
a subsequence we can suppose $Th_n \to Th$ a.e.

Therefore by linearity of $T$, we have $Tf_n \to Tf$ pointwise
a.e.\ and so
\begin{align*}
\lv Tf \rv_{L^q(Y)} & \leq \liminf_n \lv Tf_n \rv_{L^q(Y)} \qquad
\text{by Fatou's lemma} \\
& \leq \lv T \rv_{L^{p_0}(X) \to
L^{q_0}(Y)}^{1-\theta} \lv T \rv_{L^{p_1}(X) \to L^{q_1}(Y)}^\theta \liminf_n \lv f_n \rv_{\lp(X)} \\
& \qquad \qquad \text{by \eqref{eq:rti4}, the strong $(p,q)$ bound
on the simple functions,} \\
& = \lv T \rv_{L^{p_0}(X) \to L^{q_0}(Y)}^{1-\theta} \lv T
\rv_{L^{p_1}(X) \to L^{q_1}(Y)}^\theta \lv f \rv_{\lp(X)}
\end{align*}
since $f_n \to f$ in $\lp(X)$ by dominated convergence.

We have proved the desired strong $(p,q)$ bound for all $f \in
\lp(X)$, and so the proof is complete.
\end{proof}

\end{document}